\documentclass[final,onefignum,onetabnum]{siamart171218}

\input{ex_shared}

\begin{document}

\maketitle
\begin{abstract}
This paper describes a new algorithm for computing a low-Tucker-rank approximation of a tensor.
The method applies a randomized linear map to the tensor to obtain a \emph{sketch}
that captures the important directions within each mode, as well as the interactions among the modes.
The sketch can be extracted from streaming or distributed data or with a single pass over the tensor,
and it uses storage proportional to the degrees of freedom in the output Tucker approximation.
The algorithm does not require a second pass over the tensor, although it can exploit another view
to compute a superior approximation. The paper provides a rigorous theoretical guarantee
on the approximation error. Extensive numerical experiments show that that the algorithm
produces useful results that improve on the state-of-the-art for streaming Tucker decomposition.
\end{abstract}




\section{Introduction}

Large-scale datasets with natural tensor (multidimensional array) structure
arise in a wide variety of applications including
computer vision \cite{vasilescu2002multilinear},
neuroscience \cite{cichocki2013tensor},
scientific simulation \cite{austin2016parallel},
sensor networks \cite{sun2008incremental},
and data mining \cite{kolda2008scalable}.
In many cases, these tensors are too large to manipulate, to transmit,
or even to store in a single machine.
Luckily, tensors often exhibit a low-rank structure,
and can be approximated by a low-rank tensor factorization,
such as CANDECOMP/PARAFAC (CP), tensor train, or Tucker factorization~\cite{kolda2009tensor}.
These factorizations reduce the storage costs
by exposing the latent structure.
Sufficiently low rank tensors can be compressed by several orders of magnitude
with negligible loss.
However, computing these factorizations can require
substantial computational resources.
One challenge is that these large tensors
may not fit in the main memory on our computer.

In this paper, we develop a new algorithm to compute a low-rank Tucker
approximation of a tensor from streaming data,
using working storage proportional to the degrees of freedom in the output Tucker approximation.
The algorithm forms a linear sketch of the tensor,
and it operates on the sketch to compute a low-rank Tucker approximation.
The main computational work is all performed on
a small tensor whose size is proportional to the core tensor in the Tucker factorization.
We derive detailed probabilistic error bounds on the quality of the approximation
in terms of the tail energy of any matricization of the target tensor.

This algorithm is useful in at least three concrete problem settings:
\begin{enumerate}
\item{\bf Streaming:} Data about the tensor is received sequentially.
At each time, we observe
a low-dimensional slice,
an individual entry,
or an additive update to the tensor
(the ``turnstile'' model \cite{muthukrishnan2005data}).
For example, each slice of the tensor may represent one 
time step in a computer simulation or the measurements from a sensor array at a particular time.
In the streaming setting, the complete tensor is not stored; indeed, it may be much larger than
available computing resources.

Our algorithm can approximate a tensor, presented as a data stream, 
by sketching the updates and storing the sketch.
The linearity of the sketching operation guarantees that sketching commutes with
slice, entrywise, or additive updates.
Our method forms an approximation of the tensor
only after all the data has been observed,
rather than approximating the tensor-observed-so-far at any time.
This protocol allows for offline data analysis,
including many scientific applications.
Conversely, this protocol is not suitable for real-time monitoring.

\item{\bf Limited memory:} Data describing the tensor is stored on the hard disk
of a computer with much smaller RAM.
This setting reduces to the streaming setting by streaming the data
from disk.
\item {\bf Distributed:} Data describing the tensor may
be stored on many different machines.
Communicating data among these machines may be costly due to low
network bandwidth or high latency.
Our algorithm can approximate tensors stored in a distributed computing environment
by sketching the data on each slave machine and transmitting the sketch to a master,
which computes the sum of the sketches.
Linearity of the sketch guarantees that
the sum of the sketches is the sketch of the full tensor.
\end{enumerate}
In the streaming setting, the tensor is not stored, so
we require an algorithm that can compute an approximation from a single pass over the data.
In contrast, multiple passes over the data are possible in
the memory-limited or distributed settings.

This paper presents algorithms for all these settings,
among other contributions:
\begin{itemize}
	\item We present a new linear sketch for higher order tensors that we call the \emph{Tucker sketch}.
	This sketch captures the principal subspace of the tensor along each mode
	(corresponding to factor matrices in a Tucker decomposition)
	and the action of the tensor that links these subspaces (corresponding to the core).
	The sketch is linear, so it naturally handles streaming or distributed data.
	The Tucker sketch can be constructed from any dimension reduction map,
	and it can be used directly to, \eg, cluster the fibers of the tensor along some mode.
	It also can be used to approximate the original tensor.

	\item We develop a practical algorithm to compute a low-rank Tucker approximation from the Tucker sketch.
	This algorithm requires a single pass over the data to form the sketch,
	and does not require further data access.
	A variant of this algorithm, using the truncated QR decomposition,
	yields a quasi-optimal method for tensor approximation
	that (in expectation) matches the guarantees for HOSVD or ST-HOSVD up to constants.

	\item We show how to efficiently compress the output of our low-rank Tucker approximation
	to any fixed rank, without further data access.
	This method exploits the spectral decay of the original tensor,
	and it often produces results that are superior to truncated QR.
	It can also be used to adaptively choose the final size of the Tucker decomposition
	sufficient to achieve a desired approximation quality.

	\item We propose a two-pass algorithm that uses  additional data access to  improve on the one-pass method.
	This two-pass algorithm was also proposed in the simultaneous work \cite{minster2019randomized}.
	Both the one-pass and two-pass methods are appropriate for limited memory or distributed data settings.

	\item We develop provable probabilistic guarantees on the performance of both the
	one-pass and two-pass algorithms
	when the tensor sketch is composed of Gaussian dimension reduction maps.

	\item We exhibit several random maps that can be used to sketch the tensor.
	Compared to the Gaussian map,
	these alternatives are cheaper to store, easier to apply, and experimentally
	deliver similar performance as measured by the tensor approximation error.
	In particular, we demonstrate the benefits of a
	Khatri--Rao product of random matrices, which we call the tensor random projection (TRP),
	which uses exceedingly low storage.

	\item We perform a comprehensive simulation study with synthetic data,
	and we consider applications to several real datasets.
	These results demonstrate the practical performance of our method.
	In comparison to the only existing one-pass Tucker approximation algorithm \cite{malik2018low},
	our methods reduce the approximation error by more than an order of magnitude given the same storage budget.

	\item We have developed and released an open-source package in Python,
	available at \url{https://github.com/udellgroup/tensorsketch}, that implements our algorithms.
\end{itemize}

\section{Background and Related Work}

We begin with a short review of tensor notation
and some related work on low-rank matrix and tensor approximation.

\subsection{Notation}
Our paper follows the notation of \cite{kolda2009tensor}.
We denote \textit{scalar}, \textit{vector}, \textit{matrix}, and \textit{tensor} variables,
respectively, by lowercase letters ($x$), boldface lowercase letters ($\mathbf{x}$),
boldface capital letters  ($\mathbf{X}$),  and boldface Euler script letters ($\T{X}$).
For two vectors $\mathbf{x}$ and $\mathbf{y}$,
we write $\mathbf{x} \succ \mathbf{y}$ if $\mathbf{x}$ is greater than $\mathbf{y}$ elementwise.

Define $[N] := \{1,\dots, N\}$.
For a matrix $\mathbf{X} \in \mathbb{R}^{m \times n}$,
we respectively denote its $i$th row, $j$th column,
and $(i,j)$th element
by $\mathbf{X}(i,.)$, $\mathbf{X}(.,j)$, and $\mathbf{X}(i,j)$
for each $i \in [m]$, $j \in [n]$.
We write $\mathbf{X}^\dag \in \mathbb{R}^{n \times m}$ for the
\textit{Moore--Penrose pseudoinverse} of a matrix $\mathbf{X} \in \mathbb{R}^{m \times n}$.
In particular, $\mathbf{X}^\dag = (\mathbf{X}^\top \mathbf{X})^{-1}\mathbf{X}^T$
if $m \geq n$ and $\mathbf{X}$ has full column rank;
$\mathbf{X}^\dag = \mathbf{X}^T(\mathbf{XX}^T)^{-1}$,
if $m < n$ and $\mathbf{X}$ has full row rank.

\subsubsection{Kronecker and Khatri--Rao product}
For two matrices $\mathbf{A} \in \mathbb{R}^{I \times J}$ and $\mathbf{B} \in \mathbb{R}^{K \times L}$,
we define the \textit{Kronecker product}
$\mathbf{A} \otimes \mathbf{B} \in \mathbb{R}^{IK \times JL}$ as
\begin{equation}\label{kronecker}
\mathbf{A} \otimes \mathbf{B} = \left[
\begin{array}{ccc}
\mathbf{A}(1,1)\mathbf{B}  & \cdots & \mathbf{A}(1,J)\mathbf{B} \\
\vdots & \ddots & \vdots \\
\mathbf{A}(I,1)\mathbf{B} & \cdots &   \mathbf{A}(I,J)\mathbf{B}
\end{array}
\right].
\end{equation}
For $J = L$, we define the \textit{Khatri--Rao product} as $\mathbf{A} \odot \mathbf{B}$, that is, the ``matching columnwise'' Kronecker product.
The resulting matrix of size $(IK) \times J$ is defined as
\[\mathbf{A} \odot \mathbf{B} = [\mathbf{A}(\cdot,1) \otimes \mathbf{B}(\cdot,1) \cdots \mathbf{A}(\cdot,J) \otimes \mathbf{B}(\cdot,J) ] \]


\subsubsection{Tensor basics}
For a tensor $\T{X} \in \mathbb{R}^{I_1 \times \cdots \times I_N}$,
its \textit{mode} or \textit{order} is the number $N$ of dimensions.
If $I = I_1 = \cdots = I_N$, we denote $\mathbb{R}^{I_1 \times \cdots \times I_N}$ as $\mathbb{R}^{I^N}$.
The inner product of two tensors $\T{X}, \T{Y}$ is defined as
$\langle\T{X}, \T{Y}\rangle = \sum_{i_1=1}^{I_1}\cdots \sum_{i_N=1}^{I_N}\T{X}_{i_1\dots i_N}\T{Y}_{i_1\dots i_N}$.
The \textit{Frobenius norm} of $\T{X}$ is
$\|\T{X}\|_F = \sqrt{\langle\T{X},\T{X}\rangle}$.

\subsubsection{Tensor unfoldings}
Let $\bar{I} = \Pi_{j = 1}^N I_j $ and $I_{(-n)} = \Pi_{j \neq n} I_j $,
and let $\vc(\T{X})$ denote the vectorization of $\T{X}$.
The \textit{mode-$n$ unfolding} of $\T{X}$ is the matrix
$\mathbf{X}^{(n)} \in \mathbb{R}^{I_n \times I_{(-n)}}$.
The inner product for tensors matches that of any mode-$n$ unfolding:
\begin{equation}
\label{eq:F_norm_equivalent}
\langle \T{X}, \T{Y}\rangle = \langle \mathbf{X}^{(n)}, \mathbf{Y}^{(n)}\rangle = \rm{Tr}((\mathbf{X}^{(n)})^\top \mathbf{Y}^{(n)}).
\end{equation}

\subsubsection{Mode $n$-Rank of A Tensor}
The \textit{mode-$n$ rank} is the rank of the mode-$n$ unfolding.
We say a tensor $\T{X}$ has (multilinear) rank $\mathbf{r}(\T{X}) = (r_1,\dots, r_N)$
if its \textit{mode-n rank} is $r_n$ for each $n\in [N]$.

\subsubsection{Tensor contractions}
Write $\T{G} =\T{X} \times_n \mathbf{U}$ for the \textit{mode-$n$ (matrix) product}
of $\T{X}$ with $\mathbf{U} \in \mathbb{R}^{J \times I_n}$.
That is, $\T{G} =\T{X} \times_n \mathbf{U} \; \iff \; \mathbf{G}^{(n)} = \mathbf{U}\mathbf{X}^{(n)}$.
The tensor $\T{G}$ has dimension $I_1 \times \cdots \times I_{n-1} \times J \times I_{n+1} \times \cdots \times I_N$.
Mode products with respect to different modes commute:
for $\mathbf{U} \in \mathbb{R}^{J_1 \times I_n}$, $\mathbf{V} \in \mathbb{R}^{J_2 \times I_m}$,
\[
\T{X} \times_n \mathbf{U} \times_m \mathbf{V} = \T{X} \times_m \mathbf{V} \times_n \mathbf{U}
\quad \text{if} \quad n \ne m.
\]
Mode products obey the associative rule.
This rule simplifies mode products with matrices along the same mode:
for $\M{A} \in \mathbb{R}^{J_1 \times I_n}$, $\M{B} \in \mathbb{R}^{J_2 \times J_1}$,
\[
\label{eq: tensor_product_association}
\T{X}\times_n \mathbf{A} \times_n \mathbf{B} = \T{X}\times_n (\mathbf{BA}).
\]

\subsubsection{Tail energy}
To state our results, we will need a tensor equivalent for the decay in the
spectrum of a matrix.
For each unfolding $\mathbf{X}^{(n)}$,
define the $\rho$\textit{th tail energy}
\begin{equation}
(\tau_\rho^{(n)})^2 := \sum_{k>\rho}^{\min(I_n,I_{(-n)})} \sigma_{k}^2(\mathbf{X}^{(n)}), \nonumber
\end{equation}
where $\sigma_{k}(\mathbf{X}^{(n)})$ is the $k$th largest singular value of $\mathbf{X}^{(n)}$.

\subsection{Tucker Approximation}
Given a tensor $\T{X}\in\mathbb{R}^{I_1\times \dots \times I_N}$
and target rank $\mathbf{r}=(r_1, \ldots, r_N)$,
the goal of multilinear approximation is to approximate $\T{X}$ by
a tensor of multilinear rank $\mathbf{r}$.
Concretely, we search over Tucker decompositions of the approximating tensor
with \emph{core tensor}
$\T{G}\in \mathbb{R}^{r_1 \times \cdots \times r_N}$
and \emph{factor matrices} $\mathbf{U}_n \in \mathbb{R}^{I_n \times r_n}$ for $n\in [N]$
with each $\M{U}_n$ satisfying $\M{U}_n^\top \M{U}_n = \M{I}$.
For brevity, we define
$\llbracket\T{G}; \mathbf{U}_1, \ldots, \mathbf{U}_N \rrbracket
= \T{G}\times_1 \mathbf{U}_1\times_2 \cdots\times_N \mathbf{U}_N$.
Any best rank-$\mathbf{r}$ Tucker approximation is of the form
$\llbracket\T{G}^\star; \mathbf{U}_1^\star, \ldots, \mathbf{U}_N^\star \rrbracket$,
where $\T{G}^\star, \mathbf{U}_n^\star$ solve the \emph{Tucker approximation} problem
\begin{equation}
\begin{array}{ll}
\label{eq:tucker_optimization}
\mbox{minimize} & \|\T{X} -\T{G}\times_1 \times \cdots \mathbf{U}_{n+1} \times_N \mathbf{U}_N\|_F^2 \qquad \\
\mbox{subject to} & \mathbf{U}_n^\top \mathbf{U}_n = \mathbf{I}.
\end{array}
\end{equation}
The problem \cref{eq:tucker_optimization} is a challenging nonconvex optimization problem.
Moreover, the solution is not unique \cite{kolda2009tensor}.
We use the notation $\llbracket\T{X} \rrbracket_\mathbf{r}$
to represent a best rank-$\mathbf{r}$ Tucker approximation of the tensor $\T{X}$,
which in general we cannot compute.

\subsubsection{HOSVD}
The standard approach to computing a rank $\mathbf{r} = (r_1, \ldots, r_N)$ Tucker approximation for a tensor $\mathscr{X}$
begins with the higher order singular value decomposition (HOSVD) \cite{de2000multilinear,tucker1966some}
(\cref{alg:hosvd}). 
\begin{algorithm}[ht]
  \caption{Higher order singular value decomposition (HOSVD)
	\cite{de2000multilinear,tucker1966some}
	\label{alg:hosvd}}
  \textbf{Given:} tensor $\T{X}$, target rank $\V{r} = (r_1, \ldots, r_N)$
  \begin{enumerate}
  \item \emph{Factors.} For $n \in [N]$, compute the top $r_n$ left singular vectors $\M{U}_n$
  of $\M{X}^{(n)}$.
  \item \emph{Core.} \label{stage:hosvd-factor} Contract these with $\mathscr{X}$ to form the core
  \[
  \T{G} = \T{X} \times_1 \M{U}_1^T \cdots \times_N \M{U}_N^T.
  \]
  \end{enumerate}
  \textbf{Return:} Tucker approximation $\T{X}_{\rm{HOSVD}} = \llbracket\T{G}; \M{U}_1, \ldots, \M{U}_N \rrbracket$
\end{algorithm}

The HOSVD can be computed in two passes over the tensor
\cite{zhou2014decomposition, battaglino2019faster}.
We describe this method briefly here, and in more detail in the next section.
In the first pass, sketch each matricization $\M{X}^{(n)}$, $n \in [N]$,
and use randomized linear algebra
(\eg, the randomized range finder of \cite{halko2011finding})
to (approximately) recover its range $\mathbf{U}_n$.
To form the core $\T{X} \times_1 \M{U}_1^T \cdots \times_N \M{U}_N^T$
requires a second pass over $\T{X}$, since the factor matrices $\M{U}_n$
depend on $\T{X}$.
The main algorithmic contribution of this paper is to develop a method to
approximate both the factor matrices and the core in just one pass over $\T{X}$.

It is possible to improve the accuracy of the resulting approximation.
The higher order orthogonal iteration (HOOI) \cite{de2000best}, 
for example,
uses the HOSVD to initialize an alternating minimization method,
and sequentially minimizes over each of the factor matrices and the core tensor.
However, this method is rarely used in practice due to the memory and computation required.

\subsubsection{ST-HOSVD}
The sequentially truncated higher order singular value decomposition (ST-HOSVD)
modifies the HOSVD to reduce the computational burden \cite{vannieuwenhoven2012new}.
This method compresses the target tensor after extracting each factor matrix.
The resulting algorithm can be accelerated using randomized matrix approximations
\cite{minster2019randomized},
but seems to require $N$ passes over the tensor.
Hence the method is difficult to implement when the data is too large to store locally.

\begin{algorithm}[ht]
	\caption{Sequentially truncated HOSVD (ST-HOSVD)
	\cite{vannieuwenhoven2012new}
	\label{alg:st-hosvd}}
	\textbf{Given:} tensor $\T{X}$, target rank $\V{r} = (r_1, \ldots, r_N)$
	\begin{enumerate}
	\item $\T{G} = \T{X}$
	\item \text{For $n = 1$ to $N$}
	\bit
	\item Compute a best rank $r_n$ approximation of the mode $n$ unfolding of $\T{G}$:
	\[
	\M{U}_n, \M{\Sigma}_n, \M{V}_n =  \text{TruncatedSVD}( \T{G}^{(n)}, \V{r}_{n}).
	\]
	\item Form the updated tensor $\T{G}$ from its mode $n$ unfolding $\M{G}^{(n)} \gets \M{\Sigma}_n \M{V}_n^\top$.
	\eit
	\end{enumerate}
	\textbf{Return:} Tucker approximation $\T{X}_{\rm{ST-HOSVD}} = \llbracket\T{G}; \M{U}_1, \ldots, \M{U}_N \rrbracket$
\end{algorithm}

\subsubsection{Quasi-optimality}
A method for tensor approximation is called \emph{quasi-optimal} if
the error of the resulting approximation is comparable to the best possible:
more precisely, we say an approximation method is quasi-optimal with factor $d$ if
for any $\T{X}$ and any multilinear rank $\V{r}$,
the rank-$\V{r}$ approximation $\T{\hat{X}}$ produced by the method satisfies
\[
\| \T{X} - \T{\hat{X}} \|_F \leq d \| \T{X} - \llbracket \T{X} \rrbracket_{\V{r}} \|_F.
\]
We call a randomized tensor approximation method quasi-optimal if this inequality holds in expectation.
This definition shows the advantage of a quasi-optimal approximation method:
the method finds a good  approximation of the tensor whenever a good rank-$\V{r}$ approximation exists.
Moreover, it exactly recovers a rank-$\V{r}$ decomposition of a tensor that is
exactly rank $\V{r}$.

Both the HOSVD and the ST-HOSVD are quasi-optimal with factor $\sqrt{N}$
\cite{vannieuwenhoven2012new,hackbusch2012tensor,grasedyck2010hierarchical}.
This paper demonstrates the first known quasi-optimal streaming Tucker approximations.



\subsection{Previous Work}\label{sec: previous_work}



The only previous work on streaming Tucker approximation is \cite{malik2018low},
which develops a streaming method called Tucker TensorSketch (T.-TS) \cite[Algorithm 2]{malik2018low}.
T.-TS improves on the HOOI by sketching the data matrix in the least squares problems.
However, the success of the approach depends on the quality of the initial
core and factor matrices, and the alternating least squares algorithm takes
several iterations to converge.

In contrast, our work is motivated by the HOSVD (not HOOI)
and requires no initialization or iteration.
We treat the tensor as a \emph{multilinear} operator.
The sketch identifies a low-dimensional subspace \emph{for each mode of the tensor}
\mnote{for each input argument: what does this mean?}
that captures the action of the operator along that mode.
The reconstruction produces a low-Tucker-rank multilinear operator
with the same action on this low-dimensional tensor product space.
This linear algebraic view allows us to develop the
first guarantees on approximation error for this class of problems.\footnote{
The guarantees in \cite{malik2018low} hold only when a new sketch is applied
for each subsequent least squares solve;
the resulting algorithm cannot be used in a streaming setting.
In contrast, the practical streaming method T.-TS fixes the sketch for each mode,
and so has no known guarantees.
Interestingly, experiments in \cite{malik2018low} show that the method achieves
lower error using a fixed sketch (with no guarantees) than using fresh sketches at each iteration.
}
Moreover, we show numerically that
our algorithm achieves a better approximation of the original tensor given the same storage budget.

More generally, there is a large literature on randomized algorithms
for matrix factorizations and for solving optimization problems;
for example, see the review articles \cite{halko2011finding, woodruff2014sketching}.
In particular, our method is strongly motivated by the recent papers \cite{tropp2018more, tropp2019streaming},
which provide methods for one-pass matrix approximation.
The novelty of this paper is in our design of a core sketch (and reconstruction) for
the Tucker decomposition,
together with provable performance guarantees.
The proof requires a careful accounting of the errors resulting from
the factor sketches and from the core sketch.
The structure of the Tucker sketch guarantees that these errors are independent.

Many researchers have used randomized algorithms to compute tensor decompositions.
For example, the authors of \cite{wang2015fast, battaglino2018practical} apply sketching techniques to the CP decomposition,
while the author of \cite{tsourakakis2010mach} suggests sparsifying the tensor. Several papers aim to make Tucker decomposition efficient in the limited-memory or distributed settings \cite{baskaran2012efficient, zhou2014decomposition,
austin2016parallel, kaya2016high, li2015input, battaglino2019faster}.

\section{Dimension Reduction Maps}
In this section, we first introduce some commonly used
randomized dimension reduction maps together with some mathematical background,
and we explain how to calculate and update sketches.

\subsection{Dimension Reduction Map} Dimension reduction maps (DRMs)
take a collection of high-dimensional objects to a lower-dimensional space
while maintaining certain geometric properties \cite{oymak2015universality}.
For example, we may wish to preserve the pairwise distances between vectors,
or to preserve the column space of matrices.
We call the output of a DRM on an object $x$ a \emph{sketch} of $x$.

Common DRMs include matrices with i.i.d.~Gaussian entries
or i.i.d.~Rademacher entries (uniform on $\{\pm 1\}$).
The scrambled subsampled randomized fourier transform (SSRFT) \cite{woolfe2008fast}
and sparse random projections \cite{achlioptas2003database, li2006very}
can achieve similar performance with fewer computational and storage requirements;
\ifdefined \issupplement
see supplement for details.
\else
see \cref{appendix: ssrft} for details.
\fi

Our theoretical bounds rely on properties of the Gaussian DRM.
However, our numerical experiments indicate that many other DRMs
yield qualitatively similar results;
\ifdefined \issupplement
see supplement for details.
\else
see, \eg, \cref{fig:vary-k-600}, \cref{fig:vary-k-200-app}
and \cref{fig:vary-k-400-app})
in \cref{appendix:more_result}.
\fi

\subsection{Tensor Random Projection}\label{s-trp}
Here we present a strategy for reducing the storage of the random map
that makes use of the tensor random projection (TRP),
an extremely low storage structured dimension reduction map
proposed in \cite{sun2018tensor}.
The \emph{tensor random projection (TRP)}
$\mathbf{\Omega}: \prod_{n=1}^N I_n \to \mathbb{R}^k$ is defined as
the iterated Khatri--Rao product of DRMs
$\mathbf{A}_n \in \mathbb{R}^{I_n \times k}$, $n \in [N]$:
\begin{equation}
\label{eq:TRP}
\mathbf{\Omega} = \mathbf{A}_1 \odot \cdots \odot \mathbf{A}_N.
\end{equation}
Each $\mathbf{A}_n \in \mathbb{R}^{I_n \times k}$
can be a Gaussian map, a Rademacher matrix, an SSRFT, etc.
The number of constituent maps $N$ and their dimensions $I_n$ for $n \in [N]$
are parameters of the TRP,
and control the quality of the map; see \cite{sun2018tensor} for details.
The TRP map is a row-product random matrix,
which behaves like a Gaussian map in many respects \cite{rudelson2012row}.
Our experimental results confirm this behavior.

For simplicity, suppose $I_n$ is the same for each $n \in [N]$.
Then the TRP can be formed (and stored) using only $kNI$ random variables,
while standard dimension reduction maps use randomness (and storage)
that grows as $I^N$ when applied to a generic (dense) tensor.
\jatnote{Might be good to quantify the cost of applying the maps to a sparse or rank-one tensor?}
\cref{tbl: random_map} compares the computational and storage costs
for different DRMs.
\begin{table}[ht]
	\centering
	\begin{tabular}{|c|c|c|}
		\hline
		DRM & Storage & Computation \\
		\hline
		Gaussian                                                            & $kI^N$       & $kI^N$           \\ \hline
		Sparse                                                              & $\mu kI^N$    & $\mu kI^N$        \\ \hline
		SSRFT                                                               & $I^N$        & $\log(k) I^N$     \\ \hline
		TRP & $kNI$        & $kI^N$ \\
		\hline
	\end{tabular}
	\caption{Performance of different dimension reduction maps:
  We compare the storage and the computational cost of applying a DRM
	mapping $\mathbb{R}^{I^N}$ to $\mathbb{R}^k$
  to a dense tensor in $\mathbb{R}^{I^N}$. Here $\mu$ is the fraction of nonzero entries in the sparse DRMs.
	The TRP considered here is composed of Gaussian DRMs.
	}\label{tbl: random_map}
\end{table}

We do not need to explicitly form or store the TRP map $\mathbf{\Omega}$.
Instead,
we can store the constituent DRMs $\mathbf{A}_1, \dots, \mathbf{A}_N$
and compute the action of the map on the matricized tensor
using the definition of the TRP.
The additional computation required is minimal, and it empirically incurs almost no performance loss.

\section{Algorithms for Tucker approximation}
In this section, we present our proposed tensor sketch and our
algorithms for one- and two-pass Tucker approximation,
and we discuss the computational complexity and storage required
for both sparse and dense input tensors.
We present guarantees for these methods in \cref{sec:theory}.

\subsection{Tensor compression via sketching}\label{sec:sketch}


Our Tucker sketch generalizes the matrix sketch of \cite{tropp2019streaming} to higher order tensors.
To compute a Tucker sketch for tensor $\T{X} \in \reals^{I_1 \times \cdots \times I_N}$
with sketch size parameters $\V{k}$ and $\V{s}$,
draw independent, random DRMs
\begin{equation}\label{sketches}
\mathbf{\Omega}_1, \mathbf{\Omega}_2, \dots, \mathbf{\Omega}_N \quad \text{and} \quad \mathbf{\Phi}_1, \mathbf{\Phi}_2, \dots, \mathbf{\Phi}_N,
\end{equation}
with $\mathbf{\Omega}_n \in \mathbb{R}^{I_{(-n)} \times k_n}$ and
$\mathbf{\Phi}_n \in \mathbb{R}^{I_n \times s_n}$ for $n \in [N]$.
Use these DRMs to compute
\begin{align*}
\label{eq:sketchy_matrix}
\M{V}_n  &= \M{X}^{(n)}\M{\Omega}_n &\in &\reals^{I_n \times k_n}, \quad n \in [N], \\
\T{H}    &= \T{X} \times_1 \M{\Phi}_1^\top \cdots \times_N \M{\Phi}_N^\top &\in &\reals^{s_1 \times \cdots \times s_N}.
\end{align*}
The \emph{factor sketch} $\mathbf{V}_n$
captures the span of the mode-$n$ fibers of $\T{X}$ for each $n \in [N]$,
while the \emph{core sketch} $\T{H}$ contains information about
the interaction between different modes.
See \cref{alg:tensor_sketch} for pseudocode.

To produce a rank $\V{r} = \{r_1, \ldots, r_N\}$ Tucker approximation of $\T{X}$,
choose sketch size parameters $\V{k} = (k_1, \dots, k_N) \succeq \V{r}$
and $\V{s} = (s_1, \dots, s_N) \succeq \V{k}$.
(Vector inequalities hold elementwise.)
Our approximation guarantees depend closely on the parameters $\V{k}$ and $\V{s}$.
As a rule of thumb, we suggest selecting $\mathbf{s} = 2 \mathbf{k}+1$,
as the theory requires $\mathbf{s} \succ 2\mathbf{k}$,
and choosing $\mathbf{k}$ as large as possible given storage limitations.

The sketches $\mathbf{V}_n$ and $\T{H}$ are linear functions of the original tensor $\T{X}$
and so can be computed in a single pass over $\T{X}$.
Linearity enables easy computation of the sketch even
\ifdefined \issupplement
in the streaming model or distributed model (see supplement for details).
\else
in the streaming model (\cref{alg:linear_update}) or distributed model (\cref{alg:sketch_distributed}).
\fi
Storing the sketches 
requires memory $\sum_{n=1}^N I_n\cdot k_n + \Pi_{i = 1}^N s_n $:
much less than the full tensor.


\begin{algorithm}[htb]
\caption{Tucker Sketch}\label{alg:tensor_sketch}
\textbf{Given:} RDRM (a function that generates a random DRM)
  \begin{algorithmic}[1]
  \Function{TuckerSketch}{$\T{X}; \V{k}, \V{s}$}
  \State Form DRMs $\M{\Omega}_n = \text{RDRM}(I_{(-n)}, k_n)$
  and $\M{\Phi}_n = \text{RDRM}(I_n, s_n)$, $n \in [N]$
  \State Compute factor sketches $\M{V}_n= \M{X}^{(n)}\M{\Omega}_n $, $n \in [N]$
  \State Compute core sketch $\T{H} =\T{X}\times_1 \mathbf{\Phi}_1^\top \times \dots \times_N  \mathbf{\Phi}_N^\top $
  \State \Return $(\T{H}, \mathbf{V}_1,\dots,\mathbf{V}_N, \{\M{\Phi}_n, \M{\Omega}_n\}_{n \in [N]})$
  \EndFunction
\end{algorithmic}
\end{algorithm}

\begin{remark}
The DRMs $\mathbf{\Omega}_n \in \mathbb{R}^{I_{(-n)} \times k_n}$
are large---much larger than the size of the Tucker factorization we seek!
Even using a low memory mapping such as the SSRFT and sparse random map,
the storage required grows as $\mathcal{O}(I_{(-n)})$.
However, we do not need to store these matrices.
Instead, we can generate (and regenerate) them as needed using a (stored) random seed.\footnote{
Our theory assumes the DRMs are random, whereas our experiments use
pseudorandom numbers. In fact, for many pseudorandom number generators
it is \textsf{NP}-hard to determine whether the output is
random or pseudorandom \cite{arora2009computational}.
In particular, we expect both to perform similarly for tensor approximation.
}
\end{remark}

\begin{remark}
Alternatively, the TRP (\cref{s-trp}) can be used to limit the storage required for $\mathbf{\Omega}_n$.
The Khatri--Rao structure in the sketch need not match the structure in the matricized tensor.
However, we can take advantage of the structure of our problem to reduce storage even further.
We generate DRMs $\mathbf{A}_n \in \reals^{I_n \times k}$ for $n \in [N]$
and define
$\M{\Omega}_n = \M{A}_1 \odot \cdots \M{A}_{n-1} \odot \M{A}_{n+1} \odot \cdots \odot \M{A}_N$ for each $n \in [N]$.
Hence we need not store the maps $\M{\Omega}_n$, but only the
small matrices $\mathbf{A}_n$.
The storage required is thereby reduced from $\mathcal{O}(N(\prod_{n=1}^N I_n) k)$
to $\mathcal{O}((\sum_{n=1}^N I_n) k)$,
while the approximation error is essentially unchanged.
We use this method in our experiments.
\end{remark}

\subsection{Low-Rank Approximation}
Now we explain how to construct a Tucker decomposition of $\T{X}$
with target multilinear rank $\mathbf{k}$
from the factor and core sketches.

We first present a simple two-pass algorithm, \cref{alg:two_pass_low_rank_appro},
that uses only the factor sketches
by projecting the unfolded matrix of original tensor $\T{X}$ to the column space of each factor sketch.
(Notice that \cref{alg:two_pass_low_rank_appro} does not use the core sketch,
so the core sketch parameter $\V{s}$ of the Tucker sketch is set to 0.)

To project to the column space of each factor matrix, we calculate the QR decomposition of each factor sketch:
\begin{equation} \label{eqn:qr}
\mathbf{V}_n = \mathbf{Q}_n\mathbf{R}_n \quad\text{for $n \in [N]$},
\end{equation}
where $\mathbf{Q}_n \in \mathbb{R}^{I_n \times k_n}$ has orthonormal columns
and $\mathbf{R}_n \in \mathbb{R}^{k_n\times k_n}$ is upper triangular.
Consider the tensor approximation
\begin{equation} \label{eq:x_tilde}
\begin{aligned}
\T{\hat{X}}_2 &=\T{X}\times_1 \mathbf{Q}_1\mathbf{Q}_1^\top \times_2 \cdots \times_N \mathbf{Q}_N\mathbf{Q}_N^\top.
\end{aligned}
\end{equation}
This approximation admits the guarantees stated in \cref{thm:low_rank_err_two_pass}.
Using the commutativity of the mode product between different modes,
we can rewrite $\tilde{\T{X}}$ as
\begin{equation}
\label{eq: two-pass}
\hat{\T{X}}_2= \underbrace{\left[\T{X}\times \mathbf{Q}_1^\top \times_2 \cdots \times_N \mathbf{Q}_N^\top\right]}_{\T{W}_2} \times_1 \mathbf{Q}_1 \times_2 \cdots \times_N \mathbf{Q}_N
= \llbracket \T{W}_2; \mathbf{Q}_1,\ldots,\mathbf{Q}_N \rrbracket,
\end{equation}
which gives an explicit Tucker approximation $\tilde{\T{X}}$ of our original tensor.
The core approximation $\T{W}_2 \in \mathbb{R}^{k_1 \times \dots \times k_N}$
is much smaller than the original tensor $\T{X}$.
To compute this approximation, we need access to $\T{X}$ twice:
once to compute $\mathbf{Q}_1,\ldots,\mathbf{Q}_N$,
and again to apply them to $\T{X}$ in order to form $\T{W}_2$.

\begin{algorithm}[H]
  \caption{Two-Pass Sketch and Low-Rank Recovery \label{alg:two_pass_low_rank_appro}}
  \textbf{Given:} tensor $\T{X}$, sketch parameters $\V{k}$
  \ben
	\item \emph{Sketch.}
  $\left(\T{H}, \mathbf{V}_1, \ldots, \mathbf{V}_N, \{\M{\Phi}_n, \M{\Omega}_n\}_{n \in [N]}\right) =
  \textproc{TuckerSketch}\left(\T{X}; \V{k}, 0\right)$
  \item \emph{Recover factor matrices.} For $n \in [N]$,
  $
  (\mathbf{Q}_n, \sim) = \rm{QR}(\mathbf{V}_n)
  $
  \item \emph{Recover core.}
  $
  \T{W}_2 = \T{X} \times_1 \M{Q}_1 \cdots \times_N \M{Q}_N
  $
  \een
  \textbf{Return:} Tucker approximation
  $\hat{\T{X}}_2 = \llbracket\T{W}_2; \M{Q}_1, \ldots, \M{Q}_N \rrbracket$
  with rank $\preceq \V{k}$
\end{algorithm}

This two-pass algorithm, \cref{alg:two_pass_low_rank_appro}, has a simple motivation.
In the first step of the HOSVD, \cref{alg:hosvd},
we approximately compute the top $r_n$ eigenvectors of each matricization $\M{X}^{(n)}$
using the randomized SVD \cite{halko2011finding}.
Indeed, the same idea was proposed in concurrent work \cite{minster2019randomized},
which extends the idea to the ST-HOSVD and provides an error analysis.
The error analyses of the two papers essentially coincide for \cref{alg:two_pass_low_rank_appro}.
One major difference is that the authors of \cite{minster2019randomized} focus on
the computational benefits gained using the randomized SVD,
while here we focus primarily on the benefits due to reduced storage.

To find an algorithm for streaming data, when it is impossible to store the full tensor,
we require a one-pass algorithm.


\paragraph{One-Pass Approximation}
To develop an one-pass method, we must use the core sketch $\T{H}$
(the compression of $\T{X}$ using the random projections $\M{\Phi}_n$)
to approximate $\T{W}_2$
(the compression of $\T{X}$ using the random projections $\M{Q}_n$).
%
To develop intuition, consider the following calculation:
if the factor matrix approximations $\M{Q}_n$ capture the range of $\T{X}$ well, then projection onto their ranges in each mode approximately preserves the action of $\T{X}$:
\[
\T{X} \approx \T{X} \times_1 \M{Q}_1 \M{Q}_1^\top \times \cdots \times_N \M{Q}_N \M{Q}_N^\top.
\]
Recall that for tensor $\T{A}$ and matrices $\mathbf{B}$ and $\mathbf{C}$ with compatible sizes,
$\T{A} \times_n (\mathbf{B} \mathbf{C}) = (\T{A} \times_n \mathbf{C}) \times_n \mathbf{B}$.
Use this rule to recognize the two-pass core approximation $\T{W}_2$:
\[
\T{X}
      \approx \left( \T{X} \times_1 \M{Q}_1^\top \times \cdots \times_N \M{Q}_N^\top \right) \times_1 \M{Q}_1 \dots \times_N \M{Q}_N
      = \T{W}_2 \times_1 \M{Q}_1 \dots \times_N \M{Q}_N
\]
Now contract both sides of this approximate equality with the DRMs $\M{\Phi}_n$
to identify the core sketch $\T{H}$:
\[
\T{H} := \T{X} \times_1 \M{\Phi}_1^\top \dots \times_N \M{\Phi}_N^\top
\approx \T{W}_2 \times_1 \M{\Phi}_1^\top \M{Q}_1 \times \cdots \times_N \M{\Phi}_N^\top \M{Q}_N.
\]
We have chosen $\V{s} \succ \V{k}$ so that each $\M{\Phi}_n^\top \M{Q}_n$ has a left inverse
with high probability.  Hence, we can solve the approximate equality for $\T{W}_2$: 
\[
\T{W}_2 \approx \T{H} \times_1 (\M{\Phi}_1^\top \M{Q}_1)^\dagger \times \cdots \times_N (\M{\Phi}_N^\top \M{Q}_N)^\dagger =: \T{W}_1.
\]
The right-hand side of the approximation defines the one-pass core approximation $\T{W}_1$. \cref{lemma:err_core_sketch} controls the error in this approximation.

\cref{alg:one_pass_low_rank_appro} summarizes
the resulting one-pass algorithm.
One (streaming) pass over the tensor can be used to sketch the tensor;
to recover the tensor, we only access the sketches.
\cref{thm:low_rank_err} (below) bounds
the overall quality of the approximation.

\begin{algorithm}[H]
  \caption{One-Pass Sketch and Low-Rank Recovery \label{alg:one_pass_low_rank_appro}}
  \textbf{Given:} tensor $\T{X}$, sketch parameters
  $\V{k}$ and $\V{s} \succ \V{k}$
  \ben
  \item \emph{Sketch.}
  $\left(\T{H}, \mathbf{V}_1, \ldots, \mathbf{V}_N, \{\M{\Phi}_n, \M{\Omega}_n\}_{n \in [N]}\right) =
  \textproc{TuckerSketch}\left(\T{X}; \V{k}, \V{s}\right)$
  \item \emph{Recover factor matrices.} For $n \in [N]$,
  $
  (\mathbf{Q}_n, \sim) = \rm{QR}(\mathbf{V}_n)
  $
  \item \emph{Recover core.}
  $
  \T{W}_1 = \T{H} \times_1 (\M{\Phi}_1^\top \M{Q}_1)^\dagger \times \cdots \times_N (\M{\Phi}_N^\top \M{Q}_N)^\dagger
  $
  \een
  \textbf{Return:} Tucker approximation
  $\hat{\T{X}}_1 = \llbracket\T{W}_1; \M{Q}_1, \ldots, \M{Q}_N \rrbracket$
  with rank $\leq \V{k}$
\end{algorithm}

The computational complexity and storage required by \cref{alg:one_pass_low_rank_appro} is presented in \cref{tbl: time-complexity}.
These requirements compare favorably to the
only previous method for streaming Tucker approximation \cite{malik2018low};
\ifdefined \issupplement
see the supplement for details.
\else
see \cref{appendix: time-complexity} for details.
\fi

\begin{table*}[h!]
	\centering
	\begin{tabular}{|c|c|c|}
		\hline
		Stage     & Computation                                & Storage \\ \hline
		Sketching & $\mathcal{O}(((1-(s/I)^N)/(1-(s/I))+Nk)I^N)$          &              \\
		Recovery  & $\mathcal{O}((k^2s^N(1-(k/s)^N))/(1-k/s) + k^2NI)$ & $kNI + s^N$  \\
		Total     & $\mathcal{O}(((s(1-(s/I)^N))/(1-s/I)+Nk)I^N)$          &              \\ \hline
	\end{tabular}
	\caption{\label{tbl: time-complexity}
	Computational complexity of one-pass approximation (\cref{alg:one_pass_low_rank_appro})
	on tensor $\T{X} \in \mathbb{R}^{I \times \dots \times I}$ with parameters $(k,s)$,
	using a TRP composed of Gaussian DRMs inside the Tucker sketch.
	Most of the time is spent sketching the tensor $\T{X}$.
	}
\end{table*}

\subsection{Fixed-Rank Approximation}\label{sec:fixed_rank}
The low-rank approximation methods
\cref{alg:two_pass_low_rank_appro,alg:one_pass_low_rank_appro}
of the previous section produce
approximations with rank no more than $\V{k}$.
It is often valuable to truncate this approximation to a
user-specified target rank $\mathbf{r} \leq \mathbf{k}$
\cite[Figure 4]{tropp2019streaming}.
Increasing $\V{k}$ relative to $\V{r}$ can improve the quality of the final approximation
by using more intermediate storage,
without changing the storage requirements of the final approximation to $\T{X}$.
In this section, we introduce a few methods to compute fixed-rank approximations
with rank no more than $\V{r}$ by way of a sketch with parameter $\V{k} \succeq \V{r}$.

\subsubsection{Truncated QR}\label{sec:truncatedQR}
One simple fix to \cref{alg:one_pass_low_rank_appro,alg:two_pass_low_rank_appro}
results in a final approximation with rank $\V{r}$ rather than $\V{k}$:
simply replace the QR decomposition with a truncated QR decomposition \cite{gu1996efficient}.
Indeed, we will show that this simple change results in a one-pass algorithm that
achieves quasi-optimality with factor $2\sqrt{N}$, nearly matching the guarantee
for the HOSVD and ST-HOSVD.
This approach is best for tensors with many modes that are (almost) exactly rank $\V{r}$.
However, for tensors with few modes and slower spectral decay,
a more sophisticated fixed-rank approximation method outperforms this naive approach.

\subsubsection{Optimal Fixed-Rank Approximation.}
For tensors with few modes,
we recommend computing a rank-$\V{r}$ approximation to $\T{X}$
by forming an initial approximation with rank $\V{k} \succeq \V{r}$ using
a randomized method such as
\cref{alg:one_pass_low_rank_appro} or \cref{alg:two_pass_low_rank_appro}
and then truncating it to rank $\V{r}$ using a deterministic method
such as ST-HOSVD.
In previous work, we have found that rank truncation is essential to ensure
that the final approximation is fully reliable:
for matrices, we find that the top singular values and vectors of
the approximation are accurate
when $\V{k} \succapprox 4\V{r}$ \cite{tropp2019streaming}.
Rank truncation can also be used to choose the final size of the Tucker decomposition
adaptively to achieve a desired approximation quality.

For moderate $\V{k}$, it is computationally easy to truncate
an initial rank-$\V{k}$ approximation to rank $\V{r}$,
thanks to the following lemma.
\begin{lem}[Core truncation]
\label{lemma: equivalance_one_pass}
Let $\T{W}\in \mathbb{R}^{k_1 \times \cdots \times k_N}$ be a tensor with $\V{k} \succeq \V{r}$,
and let $\mathbf{Q}_n \in \mathbb{R}^{I_n\times k_n}$ be orthogonal matrices for each $n \in [N]$.
Then
\begin{equation}
\llbracket \T{W}\times_1 \mathbf{Q}_1 \cdots \times_N \mathbf{Q}_N \rrbracket_\mathbf{r} =
\llbracket \T{W} \rrbracket_\mathbf{r} \times_1 \mathbf{Q}_1 \cdots \times_N \mathbf{Q}_N. \nonumber
	\end{equation}
\end{lem}
\cref{lemma: equivalance_one_pass} shows that we can compute
the optimal rank-$\V{r}$ approximation of
the (large) tensor $\hat{\T{X}}= \T{W}\times_1 \mathbf{Q}_1 \cdots \times_N \mathbf{Q}_N$
by calculating the optimal rank-$\V{r}$ approximation of
the (small) core $\T{W}$.
Interestingly, the same result holds if we replace the best rank-$\V{r}$ Tucker approximation
$\llbracket \cdot \rrbracket$ by the HOOI \cite[Lemma A.1]{breiding2018riemannian}.

\begin{proof}[Proof of \cref{lemma: equivalance_one_pass}]
The target tensor to be approximated, $\T{W}\times_1 \mathbf{Q}_1 \cdots \times_N \mathbf{Q}_N$,
lies in the subspace spanned by the $\M{Q}_n$,
$\{ \T{X} : \T{X}^{(n)} \in \range(\mathbf{Q}_n) \}$.
By the Pythagorean theorem, any optimal Tucker decomposition also lies in this subspace.

Suppose $\llbracket \T{W'}; \mathbf{V}_1, \ldots, \mathbf{V}_N  \rrbracket$
is an optimal Tucker decomposition.
Since its $n$th unfolding is in $\range(\mathbf{Q}_n)$,
each $\mathbf{V}_n$ can be written as $\mathbf{Q}_n \mathbf{U}_n$ for some orthogonal $\mathbf{U}_n\in \reals^{k_n\times r_n}$.
Then, using the orthogonal invariance of the Frobenius norm,
\begin{align*}
&\| \T{W} \times_1 \mathbf{Q}_1 \times \cdots \times_N \mathbf{Q}_N - \T{W'} \times_1  \mathbf{Q}_1\mathbf{U}_1 \times \cdots \times_N   \mathbf{Q}_N\mathbf{U}_N\|_F  \\
 & = \|\T{W}-\T{W'} \times_1 \mathbf{U}_1\times \cdots \times_N \mathbf{U}_N\|_F \ge \|\T{W} -\llbracket \T{W} \rrbracket_\mathbf{r}\|_F\\
 & =  \|\T{W}\times_1\mathbf{Q}_1\times \cdots \times_N \mathbf{Q}_N -\llbracket \T{W} \rrbracket_\mathbf{r} \times_1 \mathbf{Q}_1 \times \cdots \times_N\mathbf{Q}_N\|_F.
\end{align*}
\end{proof}

Motivated by this lemma, to produce a fixed rank-$\V{r}$ approximation of $\T{X}$,
we compress the core tensor approximation from
\cref{alg:two_pass_low_rank_appro} or \cref{alg:one_pass_low_rank_appro}
to rank $\mathbf{r}$.
This compression is cheap because the core approximation
$\T{W} \in \mathbb{R}^{k_1 \times \dots \times k_N}$
is small.

We present this method (using ST-HOSVD as the compression algorithm)
as \cref{alg:fix_rank_appro}.
One convenient aspect of this scheme is that the rank-$\V{k}$
approximation can be stored in memory, allowing users to experiment with
different desired final ranks $\V{r}$ and with different algorithms $\mathcal A$
to compress the core to rank $\V{r}$.
Reasonable choices to compress the core include the HOSVD, 
the ST-HOSVD,
or TTHRESH \cite{ballester2019tthresh}.
It is possible to use these strategies to adaptively compute a core approximation
that achieves some target approximation error.
For example, for the HOSVD, one can successively truncate the core of the HOSVD (using the ordering property \cite[Theorem 2]{de2000multilinear});
for the ST-HOSVD, one can use the error tolerance strategy of \cite{vannieuwenhoven2012new};
and the iterative strategy of TTHRESH naturally terminates upon reaching a target error approximation.
Both HOSVD and ST-HOSVD are quasi-optimal \cite{vannieuwenhoven2012new}, while
ST-HOSVD requires less storage and computation.

\begin{algorithm}[H]
  \caption{Fixed-rank approximation \label{alg:fix_rank_appro}}
  \textbf{Given:} Tucker approximation
  $\llbracket\T{W}; \mathbf{Q}_1, \ldots, \mathbf{Q}_N\rrbracket$ of tensor $\T{X}$,
  rank target $\V{r}$,
	algorithm $\mathcal A( \T{W}, \V{r})$ to compute rank $\V{r}$ approximation to $\T{W}$ (\eg, ST-HOSVD).
  \ben
  \item \label{core_approx}\emph{Approximate core with fixed rank.}
  $\T{G}, \mathbf{U}_1, \ldots, \mathbf{U}_N =
  \mathcal{A}( \T{W},\mathbf{r})$ \label{line:core-decom}
  \item \emph{Compute factor matrices.} For $n \in [N]$,
  $\mathbf{P}_n = \mathbf{Q}_n \mathbf{U}_n$
  \een
  \textbf{Return:} Tucker approximation
  $\hat{\T{X}}_\V{r} = \llbracket\T{G}; \M{P}_1, \ldots, \M{P}_N \rrbracket$
  with rank $\leq \V{r}$
\end{algorithm}

\section{Guarantees}
\label{sec:theory}
In this section, we present probabilistic guarantees on the preceding algorithms.
We show that the approximation error for the one-pass algorithm is the sum of
the error from the two-pass algorithm and the error resulting from the core approximation.
We present most of the proofs in this section,
and defer some more technical parts to the appendix.

\subsection{Low-rank approximation}
\cref{thm:low_rank_err_two_pass} guarantees the performance of the two-pass method (\cref{alg:two_pass_low_rank_appro}).
\begin{theorem}[Two-pass low-rank approximation]
\label{thm:low_rank_err_two_pass}
Sketch the tensor $\T{X}$ using a Tucker sketch with parameters $\V{k}$
using DRMs 
with i.i.d. standard normal entries.
Then the approximation $\hat{\T{X}}_2$ computed by the two-pass method (\cref{alg:two_pass_low_rank_appro})
satisfies
\begin{equation*}
	\mathbb{E}\| \T{X} - \hat{\T{X}}_2 \|_F^2  \le  \min_{1\preceq \rho \preceq k-1}
	\sum_{n=1}^N \left(1+\frac{\rho_n}{k_n-\rho_n-1}\right)(\tau^{(n)}_{\rho_n})^2.
	\end{equation*}
\end{theorem}
This theorem shows that the proposed randomized tensor approximation
works best for tensors whose unfoldings exhibit spectral decay.
As a simple consequence, we see that the two-pass method
with $\V{k} \succ \V{r}+1$
perfectly recovers a tensor with exact (multilinear) rank $\V{r}$, since
in that case $\tau^{(n)}_{r_n} = 0$ for each $n \in [N]$.
However, the theorem states a stronger bound:
the method exploits decay in the spectrum,
wherever (in the first $k_n$ singular values of each mode $n$ unfolding)
it occurs.

Another useful consequence shows that the rank-$\V{k}$ approximation
computed with this two-pass method competes with the best rank-$\V{r}$ approximation.
\begin{corollary}\label{cor-two-pass}
Suppose $\V{k} \succeq 2\V{r} +1$.
Then the approximation $\hat{\T{X}}_2$ computed by the two-pass method (\cref{alg:two_pass_low_rank_appro}) satisfies
\begin{equation*}
	\mathbb{E}\| \T{X} - \hat{\T{X}}_2 \|_F^2  \le  2
	\sum_{n=1}^N (\tau^{(n)}_{r_n})^2 \le 2 N \| \T{X} - \llbracket \T{X} \rrbracket_r \|_F^2.
	\end{equation*}
\end{corollary}

\begin{proof}[Proof of \cref{thm:low_rank_err_two_pass}]
	Suppose $\T{\hat{X}}_2$ is the low-rank approximation from \cref{alg:two_pass_low_rank_appro}.
	Use the definition of the mode-$n$ product and the commutativity of the mode product between different modes to see that
	\begin{equation*}
	\begin{aligned}
	\T{\hat{X}}_2 &=  \left[\T{X}\times_1 \mathbf{Q}_1^\top \times_2 \cdots \times_N \mathbf{Q}_N^\top\right] \times_1 \mathbf{Q}_1\times_1 \cdots\times_N \mathbf{Q}_N\\
	&= \T{X}\times_1 \mathbf{Q}_1\mathbf{Q}_1^\top \times_2 \cdots \times_N \mathbf{Q}_N\mathbf{Q}_N^\top.
	\end{aligned}
	\end{equation*}
Here we see that $\T{\hat{X}}_2$
is a \emph{multilinear orthogonal projection}
of $\T{X}$ onto the subspace spanned by the $\M{Q}_n$,
$\{ \T{X} : \T{X}^{(n)} \in \range(\mathbf{Q}_n) \}
= \{\llbracket \T{W}; \M{Q}_1,\ldots, \M{Q}_N \rrbracket : \T{W} \in \reals^{k_1 \times \cdots \times k_N}\}$.
(See \cite{de2008tensor} for more background on multilinear orthogonal projection.)
As in \cite[Theorem 5.1]{vannieuwenhoven2012new},
we sequentially apply the Pythagorean theorem to each mode to show that
	\begin{equation}
	\|\hat{\T{X}}_2 - \T{X}\|_F^2 \le \sum_{n=1}^N  \left \| (\mathbf{I} - \mathbf{Q}_n\mathbf{Q}_n^\top) \mathbf{X}^{(n)} \right\|_F^2 .
	\end{equation}
	We then take the expectation over $\mathbf{Q}_n$ for each term in the sum
	and use \cref{cor:err-decreasing}
	to show this expectation is bounded by the corresponding tail energy,
	\[
	\mathbb{E} \left \| (\mathbf{I} - \mathbf{Q}_n\mathbf{Q}_n^\top) \mathbf{X}^{(n)} \right\|_F^2
	\leq
	\min_{1\le \rho_n \le k_n -1} \left(1+\frac{\rho_n}{k_n-\rho_n-1}\right)(\tau^{(n)}_{\rho_n})^2.
	\]

\end{proof}

\cref{thm:low_rank_err} guarantees the performance of the one-pass method \cref{alg:one_pass_low_rank_appro}.

\begin{theorem}[One-pass low-rank approximation]
\label{thm:low_rank_err}
Sketch the tensor $\T{X}$ using a Tucker sketch with parameters $\V{k}$ and $\V{s} \succeq 2\V{k}$
using DRMs 
with i.i.d. standard normal entries.
Then the approximation $\hat{\T{X}}_1$ computed with the one-pass method (\cref{alg:one_pass_low_rank_appro})
satisfies the bound
\begin{equation*}
\E \| \T{X} - \hat{\T{X}}_1 \|_F^2  \le (1+\Delta) \min_{1\preceq \rho \preceq k-1}
	\sum_{n=1}^N \left(1+\frac{\rho_n}{k_n-\rho_n-1}\right)(\tau^{(n)}_{\rho_n})^2,
\end{equation*}
where $\Delta := \max_{n=1}^N k_n / (s_n-k_n-1)$.
\end{theorem}

The resulting rank-$\V{k}$ approximation, computed in a single pass, is nearly optimal.
\begin{corollary}\label{cor-one-pass}
Suppose $\V{k} \succeq 2\V{r} +1$ and $\V{s} \succeq 2 \V{k}$.
Then the approximation $\hat{\T{X}}_2$ computed by the one-pass method (\cref{alg:one_pass_low_rank_appro}) satisfies
\begin{equation*}
	\mathbb{E}\| \T{X} - \hat{\T{X}}_1 \|_F^2  \le  4
	\sum_{n=1}^N (\tau^{(n)}_{r_n})^2 \le 4 N \| \T{X} - \llbracket \T{X} \rrbracket_r \|_F^2.
	\end{equation*}
\end{corollary}

\begin{proof}[Proof of \cref{thm:low_rank_err}]
	We decompose the approximation error into
	the error due to the factor matrix approximations
	and the error due to the core approximation.
	Recall that $\T{\hat{X}}_1$ is the one-pass approximation from \cref{alg:one_pass_low_rank_appro} and
\[
	\T{\hat{X}}_2 = \T{X}\times_1 \mathbf{Q}_1\mathbf{Q}_1^\top \times_2 \cdots \times_N \mathbf{Q}_N\mathbf{Q}_N^\top,
\]
	is the two-pass approximation from \cref{alg:two_pass_low_rank_appro}.
	The one-pass and two-pass approximations differ only in the core approximation:
	\begin{equation} \label{eq:core-error}
	\begin{aligned}
	\T{\hat{X}}_1-\hat{\T{X}}_2= (\T{W}-\T{X}\times_1 \mathbf{Q}_1^\top \times_2 \cdots \times_N \mathbf{Q}_n^\top)  \times_1 \mathbf{Q}_1 \dots \times_N \mathbf{Q}_N.
	\end{aligned}
	\end{equation}
	Thus $\T{\hat{X}}_1-\hat{\T{X}}_2$ is in the subspace spanned by the $\M{Q}_n$,
	$\{ \T{X} : \T{X}^{(n)} \in \range(\mathbf{Q}_n) \}$,
	while $\hat{\T{X}}_2 - \T{X}$ is orthogonal to that subspace. Therefore,
\[
	\langle \hat{\T{X}}_1 - \hat{\T{X}}_2, \hat{\T{X}}_2 - \T{X} \rangle = 0.
\]
	Now, use the Pythagorean theorem
	to bound the error of the one-pass approximation:
	\begin{equation}
	\label{eq:error_decom}
\| \hat{\T{X}}_1- \T{X} \|_F^2 = \| \hat{\T{X}}_1 - \hat{\T{X}}_2\|_F^2 + \|\hat{\T{X}}_2 - \T{X} \|_F^2.
	\end{equation}
	Consider the first term. Using \cref{eq:core-error}, we see that
	\begin{align*}
	\|\hat{\T{X}}_1 - \hat{\T{X}}_2\|^2_F &=
	\|(\T{W}_1 - \T{X}\times_1 \mathbf{Q}_1^\top \cdots \times_N \mathbf{Q}^\top_N)\times_1 \mathbf{Q}_1\cdots \times_N \mathbf{Q}_N \|^2_F\\
	& = \|(\T{W}_1- \T{X}\times_1 \mathbf{Q}_1^\top \cdots \times_N \mathbf{Q}^\top_N)\|_F^2,
	\end{align*}
	where we use the orthogonal invariance of the Frobenius norm for the second equality.
	Next, we use \cref{lemma:err_core_sketch} to bound the expected error from the core approximation as
\[
	\mathbb{E} \|\hat{\T{X}}_1- \hat{\T{X}}_2\|^2_F \le \Delta \|\T{X} - \hat{\T{X}}_2\|.
\]
Taking the expectation of \cref{eq:error_decom} and using this bound on the core error,
we find that
\[
	\mathbb{E} \|\T{X} - \hat{\T{X}}_1\|_F^2 \le (1+\Delta) \mathbb{E}\|\T{X} - \hat{\T{X}}_2\|_F^2
\]
	Finally, we use the two-pass approximation error bound \cref{thm:low_rank_err_two_pass}:
\[
	\mathbb{E}\|\hat{\T{X}}_2 - \T{X} \|_F^2 \le \min_{1\le \rho_n<k_n-1}\left[ \sum_{n=1}^N \left(1+\frac{\rho_n}{k_n-\rho_n-1}\right)(\tau^{(n)}_{\rho_n})^2\right].
\]
\end{proof}

We see that the additional error due to sketching the core
is a multiplicative factor $\Delta$ more than the error due to sketching
the factor matrices. This factor $\Delta$ decreases as the size of the
core sketch $\V{s}$ increases.

\cref{thm:low_rank_err} also offers guidance on how to select
the sketch size parameters $\mathbf{s}$ and $\mathbf{k}$.  In particular,
suppose that the mode-$n$ unfolding has a good rank $r_n$ approximation
for each mode $n$.  Then the choices $k_n = 2r_n + 1$ and $s_n = 2k_n + 1$
ensure that
\begin{equation*}
\mathbb{E}\| \T{X} - \hat{\T{X}}_1 \|_F^2
\leq 4 \sum_{n=1}^N (\tau_{r_n}^{(n)})^2.
\end{equation*}
More generally, as $k_n/r_n$ and $s_n/k_n$ increase,
the leading constant in the approximation error tends to one.

\subsection{Fixed-rank approximation}

We now present bounds on the error of the fixed rank-$\V{r}$ approximations
produced either using the truncated QR method in \cref{alg:one_pass_low_rank_appro,alg:two_pass_low_rank_appro}
or by using the fixed-rank approximation on the output of
\cref{alg:one_pass_low_rank_appro,alg:two_pass_low_rank_appro}.
The former method produces algorithms that are quasi-optimal
with factor $\sqrt{2N}$ (two-pass) or $2\sqrt{N}$ (one-pass),
matching the rate of the HOSVD and the ST-HOSVD.
The latter method produces algorithms that are quasi-optimal
with factor that grows linearly in the number of modes $N$,
but which adapts better to spectral decay.
For tensors with few modes and high dimension, such as those that appear in our
numerical experiments, the latter methods substantially outperform the former.

The resulting bounds show that the best rank-$\V{r}$
approximation of the output from the one- or two-pass algorithms
is comparable in quality to a true best rank-$\V{r}$ approximation of the input tensor.
An important insight 
is that the sketch size parameters $\mathbf{s}$ and $\mathbf{k}$
that guarantee a good low-rank approximation
also guarantee a good fixed-rank approximation:
the error due to sketching depends only
on the sketch size parameters $\V{k}$ and $\V{s}$,
and not on the target rank $\V{r}$.

\subsubsection{Truncated QR}
We can modify the argument in the proofs of 
\cref{thm:low_rank_err_two_pass,thm:low_rank_err}
to provide an error bound for the rank-$\V{r}$ approximation
obtained by truncating the QR decomposition in \cref{alg:one_pass_low_rank_appro,alg:two_pass_low_rank_appro}
to rank $\V{r}$ as in \cref{sec:truncatedQR}.
This bound will allow us to show quasi-optimality of the resulting algorithm.
\begin{theorem}[Fixed-rank approximation via truncated QR] 
\label{thm:fix_rank_err_one_time}
Sketch the tensor $\T{X}$ using a Tucker sketch with parameters $\V{k}$
using DRMs with i.i.d. standard normal entries.
The rank-$\V{r}$ approximation $\hat{\T{X}}_2$ computed with the two-pass method (\cref{alg:two_pass_low_rank_appro}),
using a rank-$\V{r}$ truncated QR in step 2 of the algorithm,
satisfies
\begin{equation*}
\mathbb{E}\| \T{X} - \hat{\T{X}}_{2} \|_F^2  \le
\sum_{n=1}^N \left(1+\frac{r_n}{k_n-r_n-1}\right)(\tau^{(n)}_{r_n})^2.
\end{equation*}
Similarly, the rank-$\V{r}$ approximation $\hat{\T{X}}_1$ computed with the one-pass method (\cref{alg:one_pass_low_rank_appro}),
using a rank-$\V{r}$ truncated QR in step 2 of the algorithm,
satisfies
\begin{equation*}
\mathbb{E}\| \T{X} - \hat{\T{X}}_{1} \|_F^2  \le  (1+\Delta)
\sum_{n=1}^N \left(1+\frac{r_n}{k_n-r_n-1}\right)(\tau^{(n)}_{r_n})^2,
\end{equation*}
where $\Delta := \max_{n=1}^N r_n / (s_n-r_n-1)$.
\end{theorem}

\begin{proof}
For the two-pass error, in the proof of \cref{thm:low_rank_err_two_pass} 
use the tail bound from \cref{lemma:sketchy_column_space_err}
to bound the error when $\M{Q}_n \in \reals^{I_n \times r_n}$ is chosen by the truncated QR algorithm \cite{gu1996efficient}.
For the one-pass error,
use \cref{lemma:err_core_sketch} to
show that the error can be no more than a factor $(1+\Delta)$ times the error bound for the two-pass approximation with truncated QR,
where $\Delta := \max_{n=1}^N r_n / (s_n-r_n-1)$.
\end{proof}

%

\begin{corollary}[Quasi-optimality with truncated QR] 
For a given target rank $\V{r}$, choose the sketch size parameters $\V{k} = 2\V{r}+1$ and $\V{s} = 2\V{r}+1$.
Replace step 2 of \cref{alg:one_pass_low_rank_appro,alg:two_pass_low_rank_appro}
by a rank-$\V{r}$ truncated QR.
The resulting algorithms produce quasi-optimal rank-$\V{r}$ approximations.
Specifically, the two-pass approximation (\cref{alg:two_pass_low_rank_appro}) with truncated QR
is quasi-optimal with factor $\sqrt{2N}$
and the one-pass approximation (\cref{alg:one_pass_low_rank_appro}) with truncated QR
is quasi-optimal with factor $2\sqrt{N}$.
\end{corollary}


In simultaneous work, the authors of \cite{minster2019randomized} also prove the two-pass
approximation with truncated QR is quasi-optimal.

\subsubsection{Optimal fixed-rank approximation}
We now bound the error induced by applying the fixed-rank approximation method
\cref{alg:fix_rank_appro} to a given (random) low-rank approximation.

Recall that $\llbracket \T{X} \rrbracket_r$ returns a best rank-$\V{r}$ approximation to $\T{X}$.

\begin{lemma}
\label{thm:fix_rank_err}
For any tensor $\T{X}$ and random approximation $\hat{\T{X}}$ of the same size,
\begin{equation*} \label{eq-fixed-rank-bound}
\mathbb{E}\|\T{X} - \llbracket \hat{\T{X}} \rrbracket_{\V{r}} \|_F \le
\|\T{X} - \llbracket \T{X}\rrbracket_\mathbf{r}\|_F
+2\sqrt{\mathbb{E}\|\T{X}-\T{\hat{X}}\|_F^2}.
\end{equation*}
\end{lemma}

\begin{proof}[Proof of \cref{thm:fix_rank_err}]
	Our argument follows the proof of \cite[Proposition 6.1]{tropp2017practical}:
	\begin{align*}
	\|\T{X} - \llbracket \hat{\T{X}} \rrbracket_\mathbf{r}\|_F
	&  \le \|\T{X} -  \hat{\T{X}}\|_F+\|\hat{\T{X}} -  \llbracket\hat{\T{X}}\rrbracket_\mathbf{r}\|_F\\
	&\le \|\T{X} -  \hat{\T{X}}\|_F+\|\hat{\T{X}} -  \llbracket \T{X}\rrbracket_\mathbf{r}\|_F \\
	& \le \|\T{X} -  \hat{\T{X}}\|_F+\|\hat{\T{X}} - \T{X}  + \T{X} - \llbracket \T{X}\rrbracket_\mathbf{r}\|_F \\
	&\le 2\|\T{X} - \hat{\T{X}} \|_F + \|\T{X} -  \llbracket \T{X} \rrbracket_\mathbf{r}\|_F.
	\end{align*}
	The first and the third lines use the triangle inequality,
	and the second line follows from the definition of the best rank-$r$ approximation.
	Take the expectation and use Lyapunov's inequality to finish the proof.
\end{proof}

\begin{corollary}[Quasi-optimality with truncated core] 
Suppose $\V{k} \succeq 2\V{r} +1$ and $\V{s} \succeq 2 \V{k}$
and	the core approximation $\mathcal A$ in \cref{alg:fix_rank_appro}
computes an optimal rank-$\V{r}$ approximation to its input $\T{W}$.
That is, $\mathcal A( \T{W}, \V{r}) = \llbracket \T{W} \rrbracket_r$.
Then the rank-$\V{r}$ approximation algorithm produced by composing
the two-pass approximation (\cref{alg:two_pass_low_rank_appro}),
resp.~the one-pass approximation (\cref{alg:one_pass_low_rank_appro}),
with \cref{alg:fix_rank_appro}
is quasi-optimal with factor $\sqrt{2N}$, resp.~$2\sqrt{N}$ for one pass.
\end{corollary}
\begin{proof}
	Use \cref{thm:fix_rank_err} together with \cref{cor-two-pass} or \cref{cor-one-pass}.
\end{proof}


Of course, in general it is not possible to compute an optimal rank-$\V{r}$ approximation for the core.
We can still bound the error of the resulting approximation if the
approximation algorithm $\mathcal{A}$ is quasi-optimal
using the following lemma.
\begin{lemma}
Suppose the tensor $\T{X}$ and random approximation $\hat{\T{X}}$ satisfy
	\[
	\mathbb{E} \|\T{X} - \hat{\T{X}} \|_F \le C(N) \| \T{X} - \llbracket \T{X} \rrbracket_{\V{r}}\|_F.
	\]
Further suppose algorithm $\mathcal A( \T{W}, \V{r})$ computes a quasi-optimal
rank-$\V{r}$ approximation to $\T{W}$ with factor $C'(N)$. 
Then
\begin{equation}
\mathbb{E} \|\T{X} - \mathcal A(\hat{\T{X}}, \V{r}) \|_F \le
(C(N)C'(N)+C(N)+C'(N)) \|\T{X} - \llbracket \T{X} \rrbracket_{\mathbf{r}}\|_F.
\end{equation}
\end{lemma}

\begin{proof}
We calculate that
	\begin{align*}
	\mathbb{E} \|\T{X} - \mathcal A(\hat{\T{X}}, \V{r}) \|_F
	& \le \mathbb{E}\left[ \|\T{X} - \T{\hat{X}}\|_F + \|\T{\hat{X}} - \mathcal A(\hat{\T{X}}, \V{r}) \|_F \right] & \\
	& \le C(N)\|\T{X} - \llbracket \T{X} \rrbracket_{\V{r}}\|_F + C'(N) \mathbb{E}\|\T{\hat{X}} - \llbracket \hat{\T{X}}\rrbracket_{\mathbf{r}} \|_F  & \\
	& \le C(N)\|\T{X} - \llbracket \T{X} \rrbracket_{\V{r}}\|_F + C'(N)\mathbb{E}\|\T{\hat{X}}- \llbracket \T{X}\rrbracket_{\mathbf{r}} \|_F &
	\\
	& \le C(N)\|\T{X} - \llbracket \T{X} \rrbracket_{\V{r}}\|_F + C'(N) \mathbb{E}\left(\|\T{X}- \T{\hat{X}} \|_F +  \|\T{X}- \llbracket \T{X}\rrbracket_{\mathbf{r}} \|_F \right)\\
	& \le (C'(N)+C(N)+C'(N)C(N)) \|\T{X}- \llbracket \T{X}\rrbracket_{\mathbf{r}} \|_F.
\end{align*}
\end{proof}


\begin{corollary}
Suppose $\V{k} \succeq 2\V{r} +1$ and $\V{s} \succeq 2 \V{k}$
and	the core approximation $\mathcal A$ in \cref{alg:fix_rank_appro}
is quasi-optimal with factor $\sqrt{N}$ (such as the ST-HOSVD).
Then the rank-$\V{r}$ approximation algorithm produced by composing
the two-pass approximation (\cref{alg:two_pass_low_rank_appro}),
resp.~the one-pass approximation (\cref{alg:one_pass_low_rank_appro})
with \cref{alg:fix_rank_appro}
is quasi-optimal with factor $(1+\sqrt{2})\sqrt{N} + \sqrt{2}N$, resp.~$2N+3\sqrt{N}$ for one pass.
\end{corollary}
\begin{proof}
	Use \cref{thm:fix_rank_err} together with \cref{cor-two-pass} or \cref{cor-one-pass}.
\end{proof}

\section{Numerical Experiments}\label{s-experiments}

\begin{figure}
	\centering
	\begin{subfigure}{0.3\textwidth}
		\includegraphics[scale = 0.24]{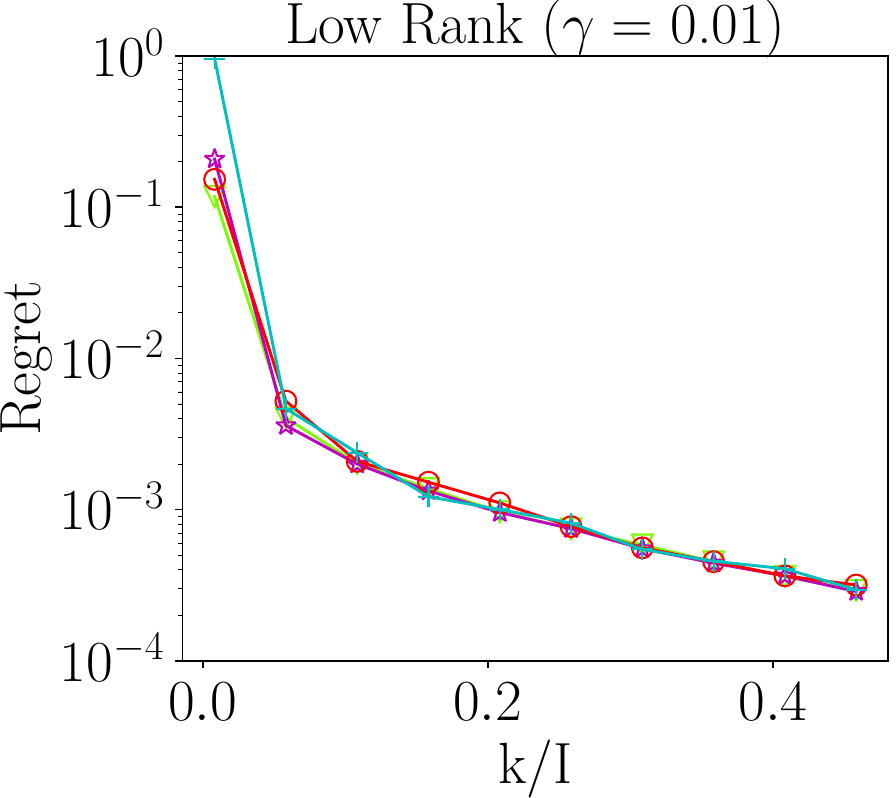}
	\end{subfigure}
	\begin{subfigure}{0.3\textwidth}
		\includegraphics[scale = 0.24]{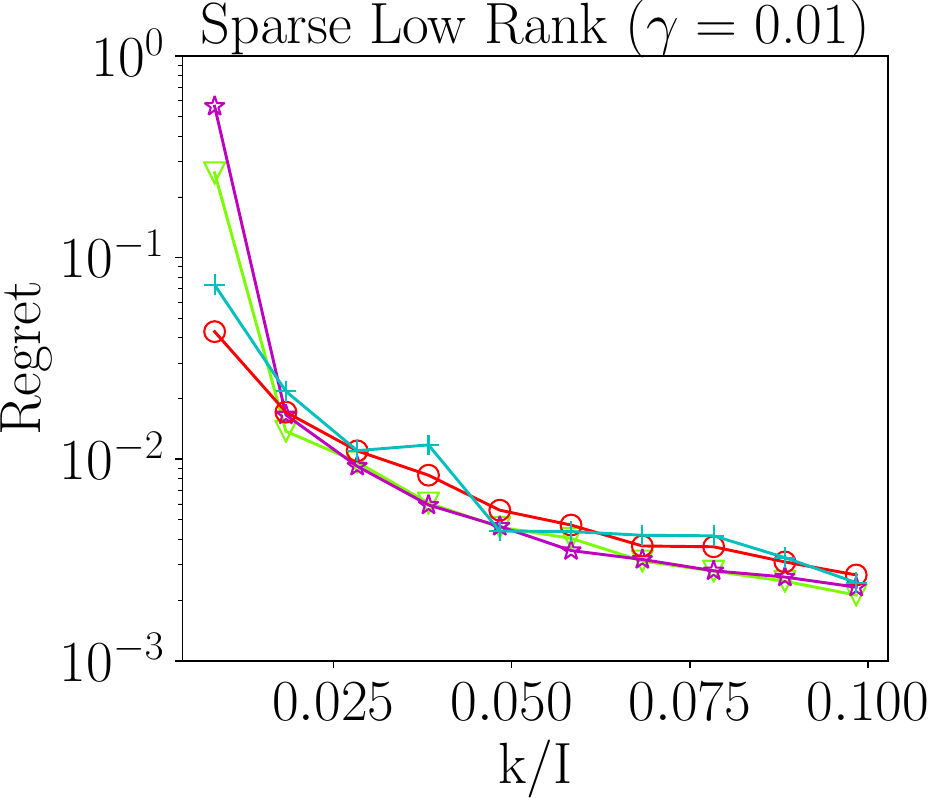}
	\end{subfigure}
	\begin{subfigure}{0.3\textwidth}
		\includegraphics[scale = 0.24]{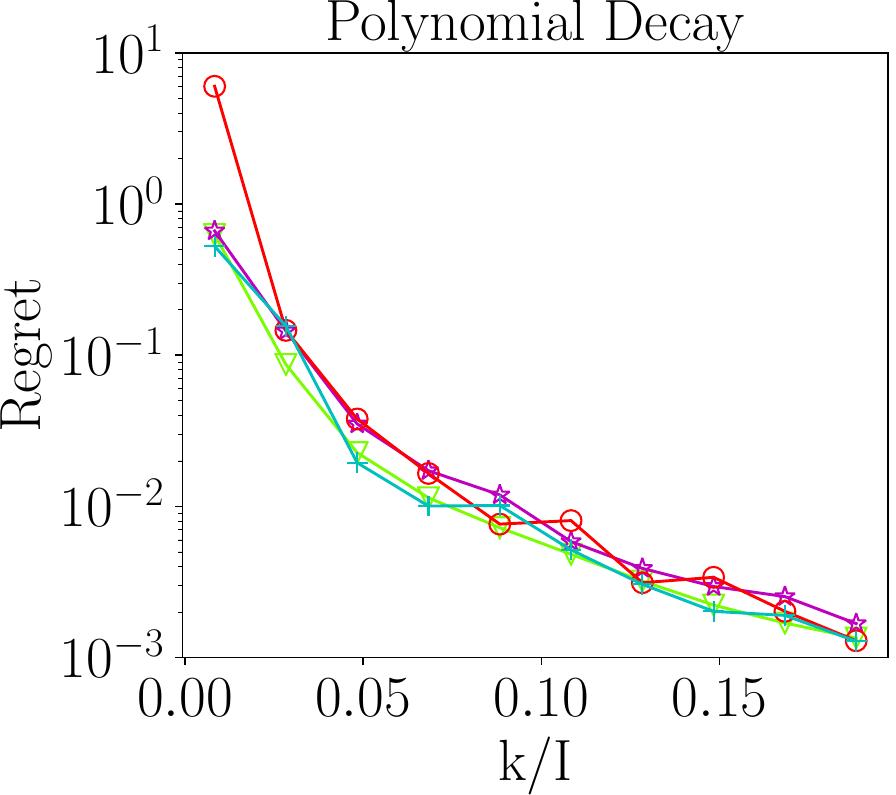}
	\end{subfigure}\\
	\begin{subfigure}{0.3\textwidth}
		\includegraphics[scale = 0.24]{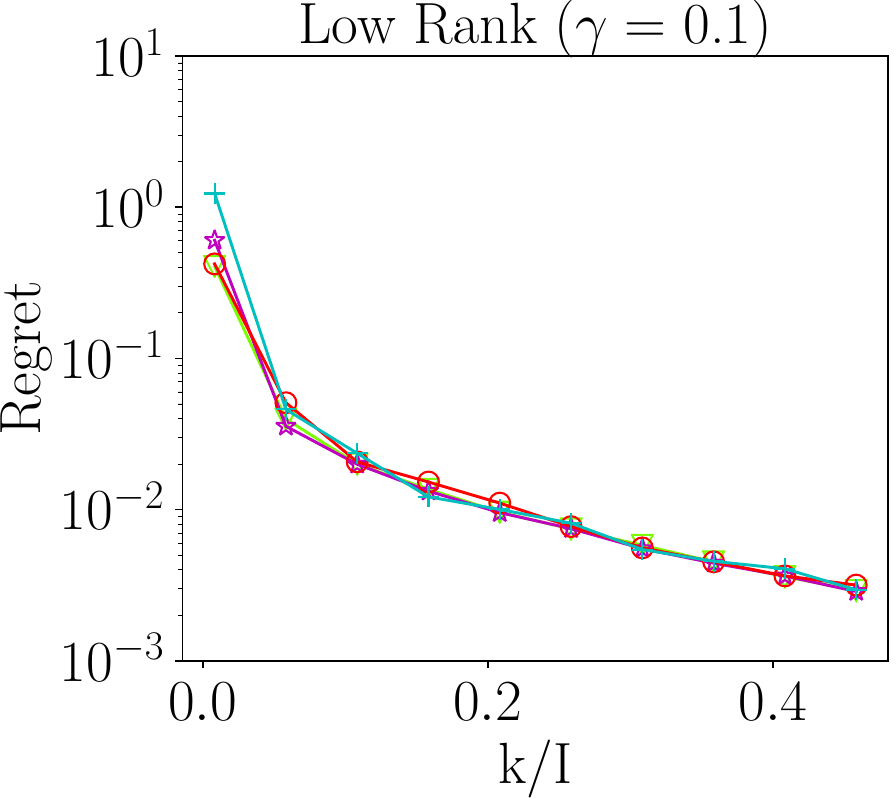}
	\end{subfigure}
	\begin{subfigure}{0.5\textwidth}
		\includegraphics[scale = 0.24]{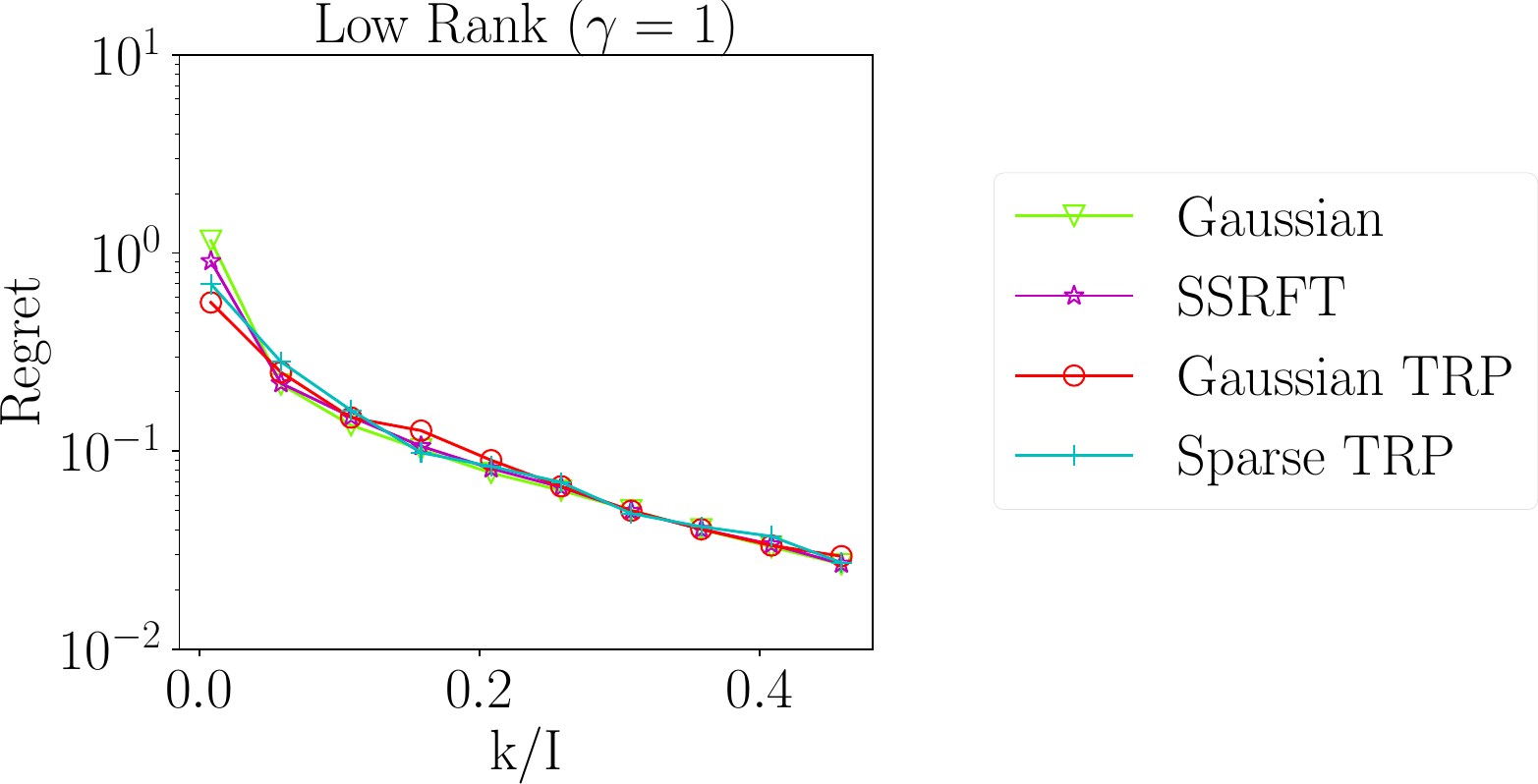}
	\end{subfigure}
	\\
	\caption{\textit{Different DRMs perform similarly.}
		We approximate three-dimensional synthetic tensors (see \cref{s-synthetic-data}) with $I = 600$,
		using our one-pass algorithm with $r = 5$ and varying $k$ ($s = 2k+1$),
		using different DRMs in the Tucker sketch.
		\label{fig:vary-k-600}
	}
\end{figure}

\begin{figure}
	\centering
	\begin{subfigure}{0.3\textwidth}
		\includegraphics[scale = 0.24]{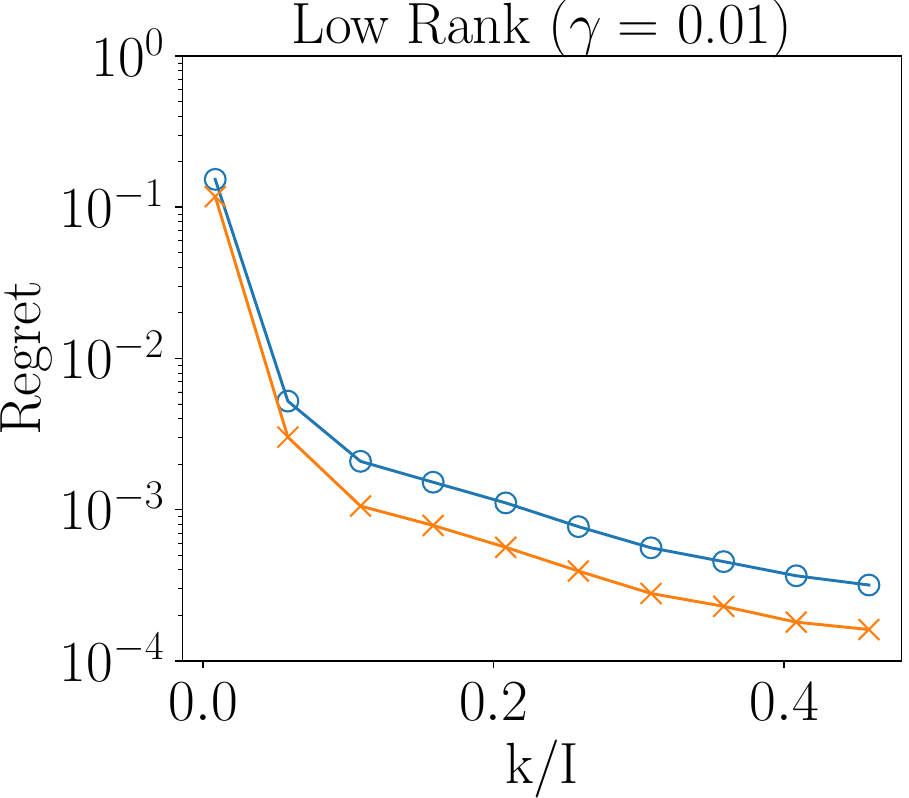}
	\end{subfigure}
	\begin{subfigure}{0.3\textwidth}
		\includegraphics[scale = 0.24]{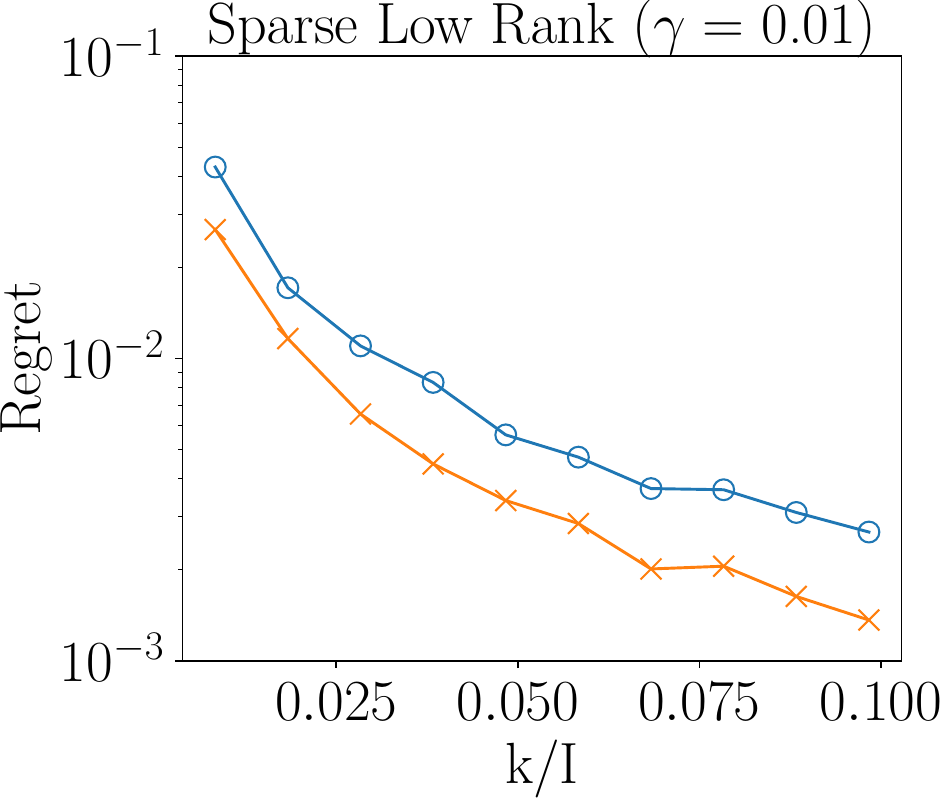}
	\end{subfigure}
	\begin{subfigure}{0.3\textwidth}
		\includegraphics[scale = 0.24]{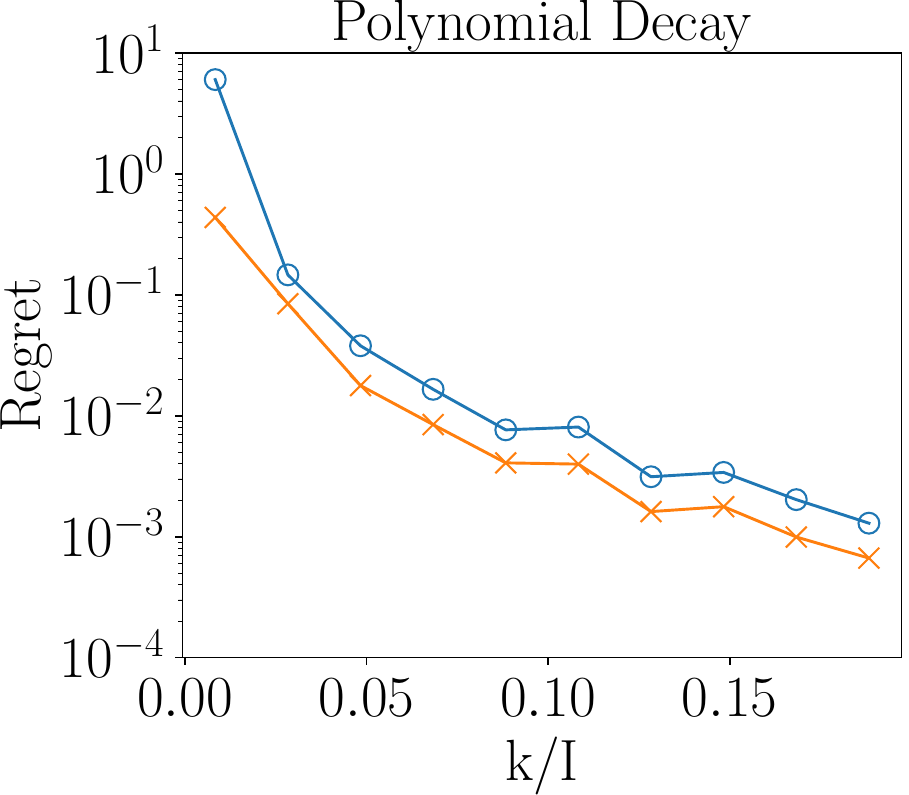}
	\end{subfigure}\\
	\begin{subfigure}{0.3\textwidth}
		\includegraphics[scale = 0.24]{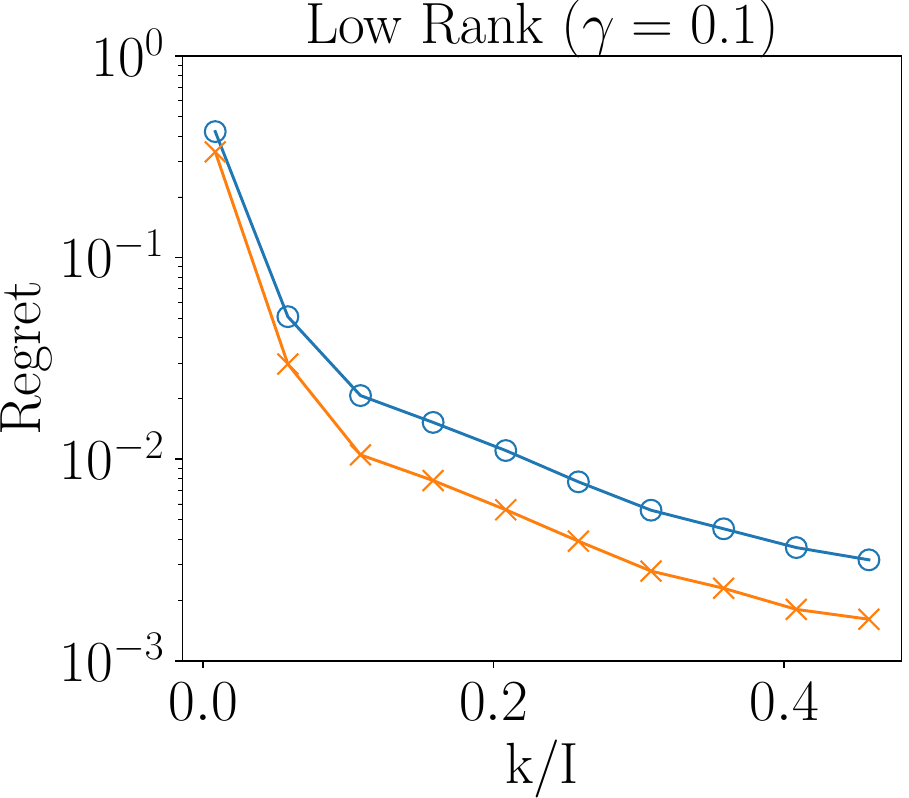}
	\end{subfigure}
	\begin{subfigure}{0.55\textwidth}
		\includegraphics[scale = 0.24]{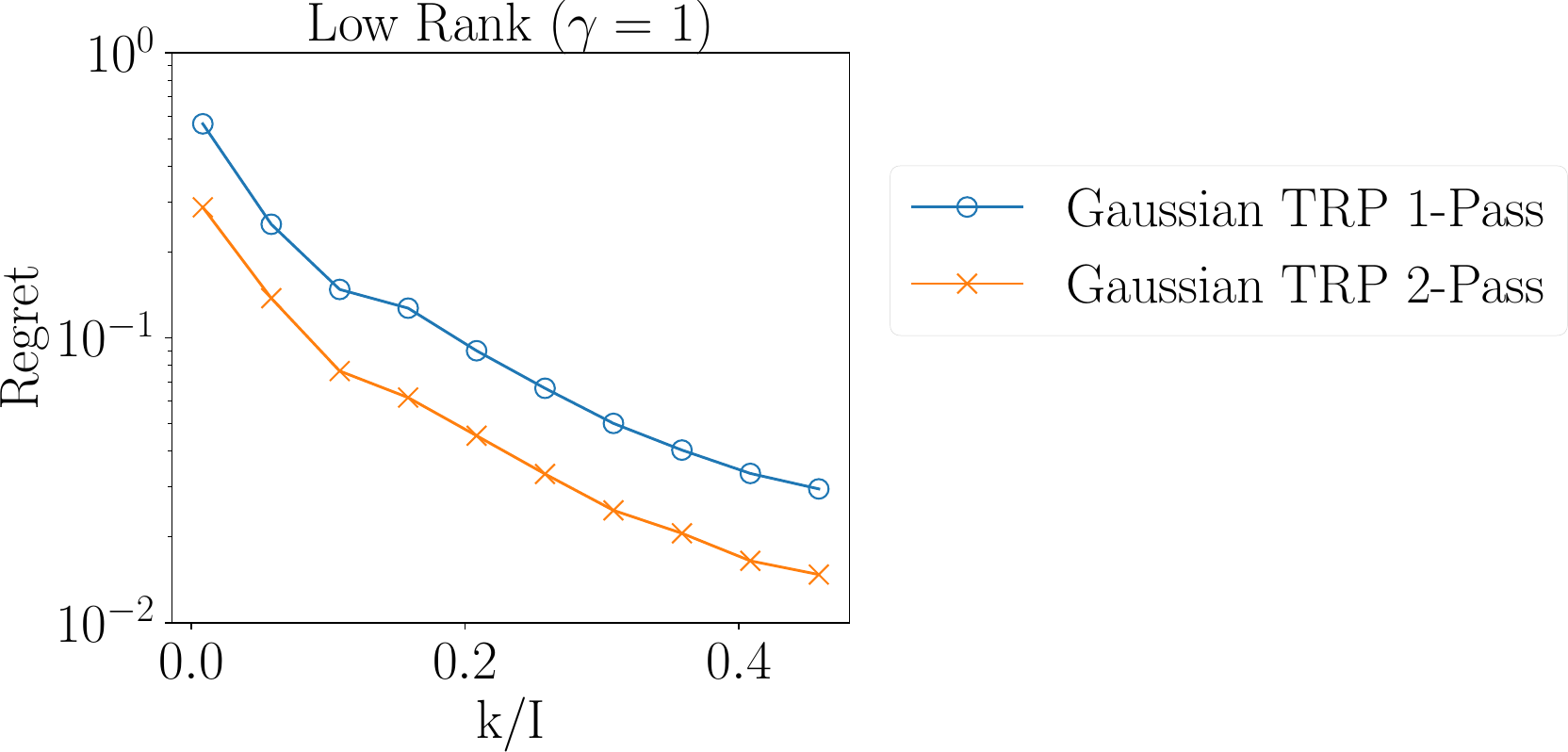}
	\end{subfigure}
	\\
	\caption{\textit{Two-pass improves on one-pass.}
	We approximate three-dimensional synthetic tensors (see \cref{s-synthetic-data}) with $I = 600$,
	using our one-pass and two-pass algorithms with $r = 5$ and varying $k$ ($s = 2k+1$),
	using the Gaussian TRP in the Tucker sketch.
	\label{fig:vary-k-600-compare}
  }
\end{figure}
In this section, we study the performance of our streaming Tucker approximation methods.
We compare the performance using various different
DRMs, including the TRP.
We also compare our method with the algorithm proposed by \cite{malik2018low}
to show that, for the same storage budget, our method produces better approximations.
Our two-pass algorithm outperforms the one-pass version, as expected.
(Contrast this to \cite{malik2018low}, where the multi-pass method
performs less well than the one-pass version.)

%


\subsection{Error metrics}

We measure the quality of an approximation $\hat{\T{X}}$
to the original tensor $\T{X}$ using two different metrics.
One is the relative error:
\[
\mbox{relative error:} \qquad \|\T{X} - \hat{\T{X}}\|_F / \|\T{X}\|_F.
\]
However, many tensors are not close to low rank, in which case
every low-rank approximation will incur high relative error.

Our methods cannot solve the problem that many tensors are not low rank; instead,
our goal is just to propose a faster, cheaper, more memory-efficient method to
compute a low-rank approximation to the tensor that is \emph{almost} as good
as one computed using a more expensive method like the HOOI, HOSVD, or ST-HOSVD.
To facilitate comparisons among approximation algorithms,
we define another metric that we call \emph{regret}. 
We found that the HOOI performs marginally better than the ST-HOSVD approximation
on the examples featured in this section.
To simplify plots and interpretations, we treat the HOOI as the gold standard,
and we define the \emph{regret} of an approximation relative to the HOOI as
\[
\left(\|\T{X} - \hat{\T{X}}\|_F - \|\T{X} - \T{X}_\text{HOOI}\|_F \right) / \|\T{X}\|_F.
\]
The regret measures the increase in error incurred by
using the approximation $\hat{\T{X}}$ rather than $\T{X}_\text{HOOI}$.
The regret of HOOI is 0.
An approximation with a regret of $.01$ is only 1\% worse
than HOOI, relative to the norm of the target tensor $\T{X}$.

\subsection{Computational platform}

We ran all experiments on a server with 128 Intel Xeon E7-4850 v4 2.10GHz CPU cores and 1056GB memory.
All experiments are implemented in Python.
We use the default implementations available in the
Python package \emph{tensorly} \cite{kossaifi2019tensorly}
for tensor algorithms such as the HOOI and ST-HOSVD.
Code for the one- and two-pass approximation algorithms
is available on Github at \url{https://github.com/udellgroup/tensorsketch},
as is the code that generates the experiments in this paper.

\subsection{Synthetic experiments}\label{s-synthetic-data}
\begin{figure}
	\centering
	\begin{subfigure}{0.3\textwidth}
		\includegraphics[scale = 0.24]{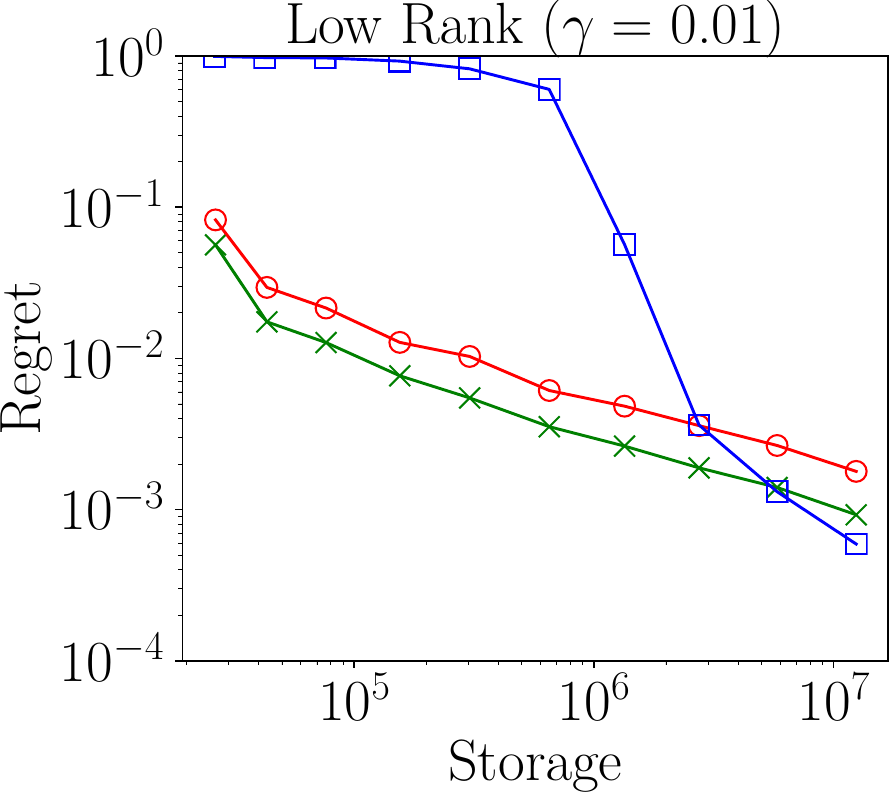}
	\end{subfigure}
	\begin{subfigure}{0.3\textwidth}
		\includegraphics[scale = 0.24]{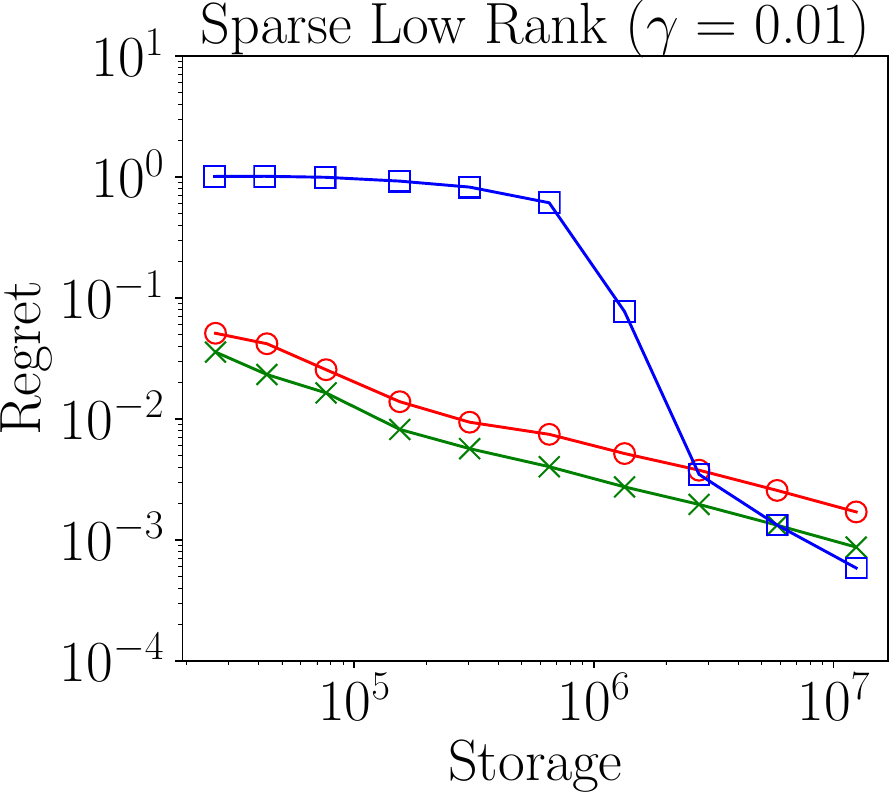}
	\end{subfigure}
	\begin{subfigure}{0.3\textwidth}
		\includegraphics[scale = 0.24]{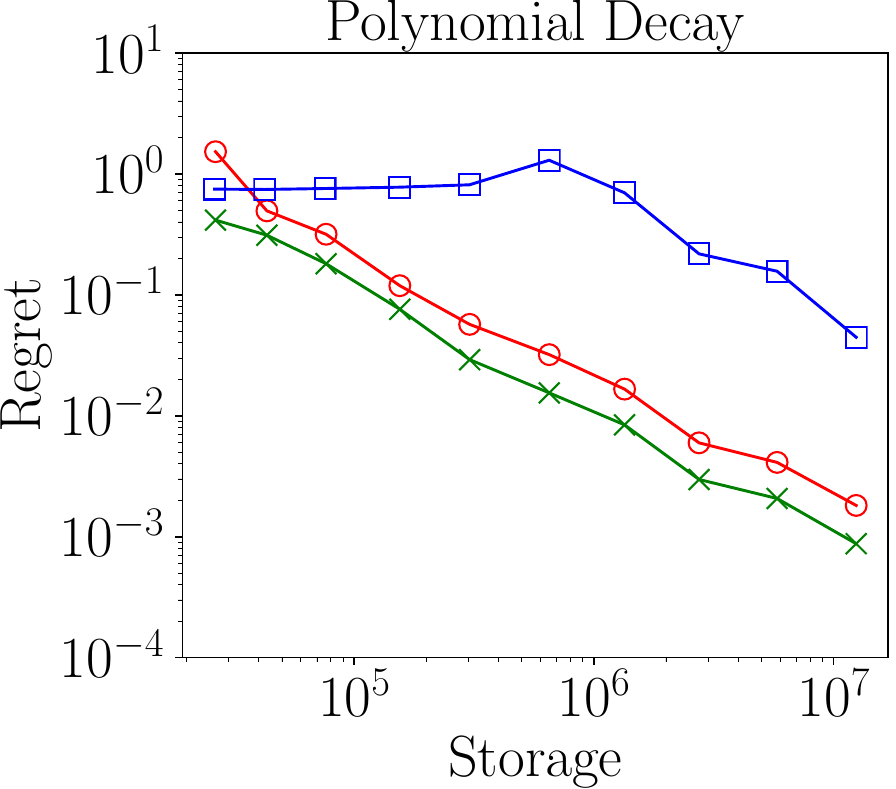}
	\end{subfigure}\\
	\begin{subfigure}{0.3\textwidth}
		\includegraphics[scale = 0.24]{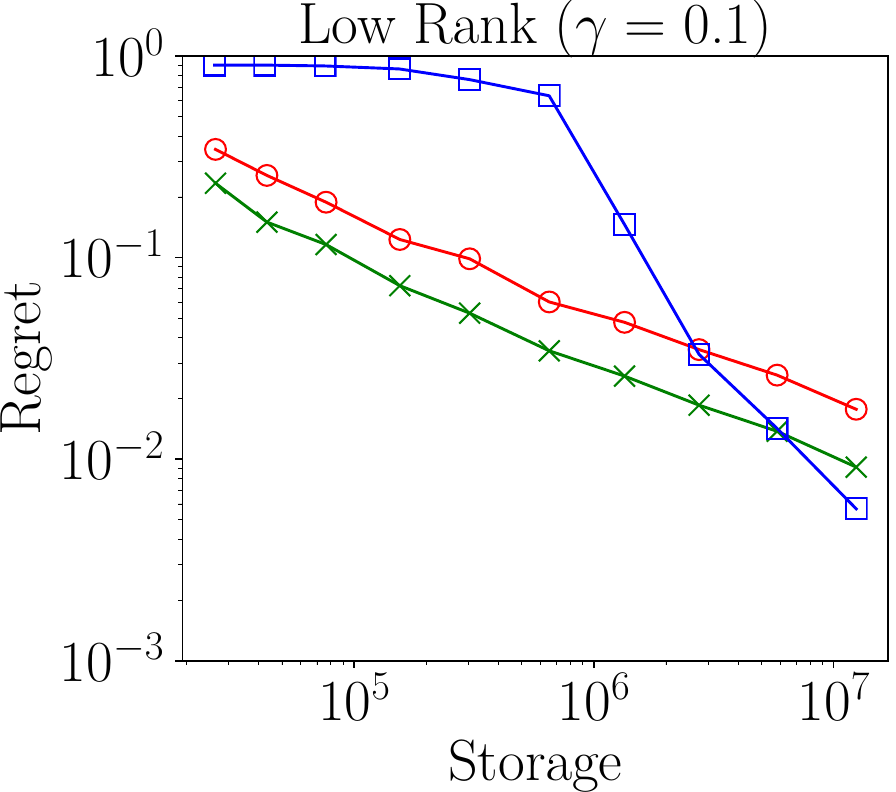}
	\end{subfigure}
	\begin{subfigure}{0.43\textwidth}
		\includegraphics[scale = 0.24]{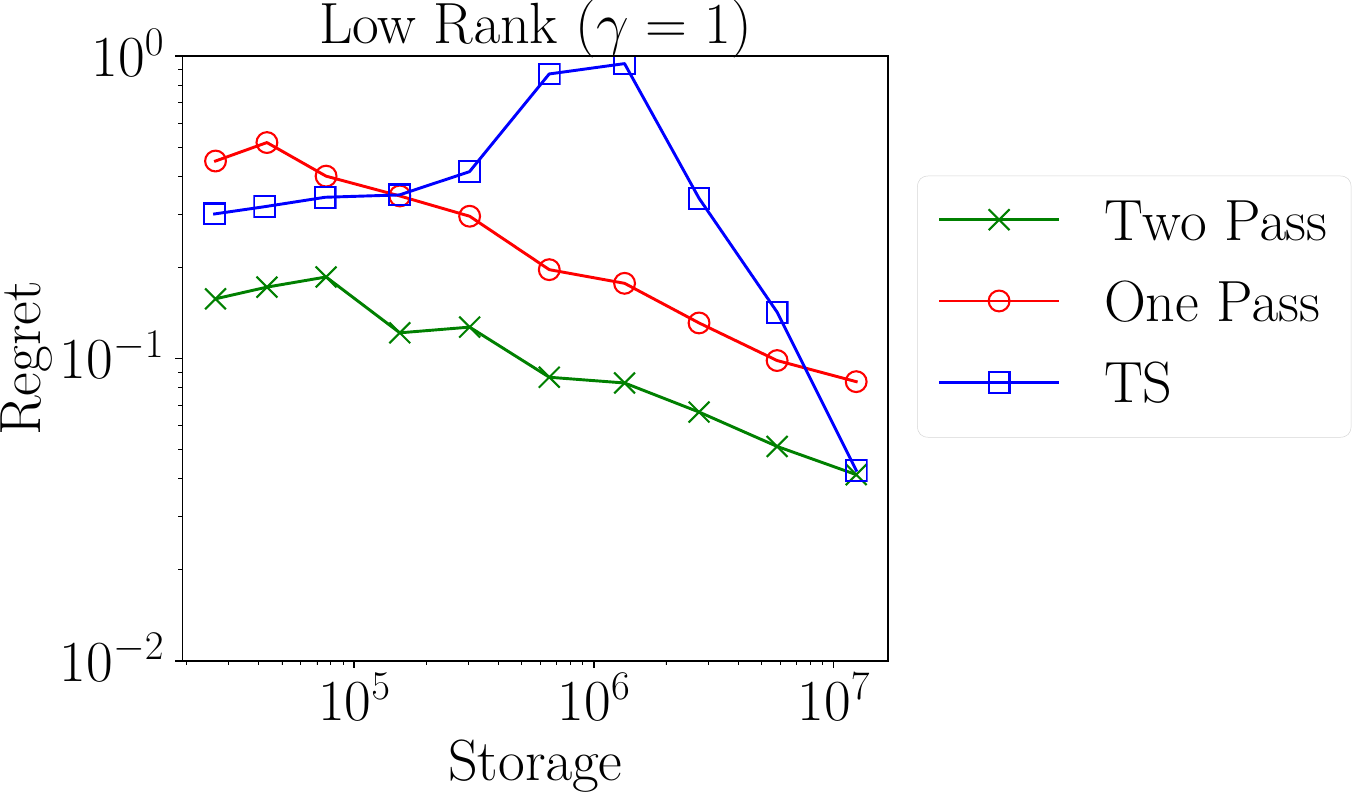}
	\end{subfigure}\\
	\caption{\textit{Approximations improve with more memory: synthetic data.}
	We approximate three-dimensional synthetic tensors (see \cref{s-synthetic-data}) with $I = 300$,
	using  T.-TS and our one-pass and two-pass algorithms
	with the Gaussian TRP to produce approximations with equal ranks $r=10$.
	Notice every marker on the plot corresponds to a 2700$\times$ compression!}\label{fig:vary-memory}
\end{figure}

All synthetic experiments use an input tensor with equal side lengths $I$.
We consider three different data generation schemes:
\begin{itemize}
\item \emph{Low-rank noise.} Generate a core tensor $\T{C} \in \mathbb{R}^{r^N}$
with entries drawn i.i.d. from the uniform distribution $U(0,1)$.
Generate $N$ random matrices $\mathbf{B}_1, \dots, \mathbf{B}_N \in \mathbb{R}^{r \times I}$ with i.i.d. $\mathcal{N}(0,1)$ entries,
and let $\mathbf{A}_1, \dots, \mathbf{A}_N \in \mathbb{R}^{r \times I}$ be orthonormal bases for their respective column spaces.
Define $\T{X}^\natural = \T{C} \times_1 \mathbf{A}_1 \cdots \times_N \mathbf{A}_N$
and the noise parameter $\gamma > 0$.
Generate an input tensor as
	$\T{X} = \T{X}^\natural + (\gamma \|\T{X}^\natural\|_F / I^{N/2})\T{\epsilon}$
	where the noise $\T{\epsilon}$ has i.i.d. $\mathcal{N}(0,1)$ entries.
\item \emph{Sparse low-rank noise.} We construct the input tensor $\T{X}$ as above (low-rank noise),
but with sparse factor matrices $\mathbf{A}_n$:
If $\delta_n$ is the sparsity (proportion of nonzero elements) of $\mathbf{A}_n$,
then the sparsity of the true signal $\T{X}^\natural$ scales as $r^N \prod_{n=1}^N \delta_n$.
We use $\delta_n = 0.2$ unless otherwise specified.
\item \emph{Polynomial decay.} We construct the input tensor $\T{X}$ as
  \begin{equation*}
      \T{X} = \mathop{\textbf{superdiag}}(1,\dots,1, 2^{-t},3^{-t},\dots, (I-r)^{-t}).
  \end{equation*}
  The first $r$ entries are 1.
	Recall that $\mathop{\textbf{superdiag}}$ converts a vector to an $N$-dimensional superdiagonal tensor.
  Our experiments use $t = 1$.
\end{itemize}
Our goal in including the polynomial and sparse setups
is to demonstrate that the method performs robustly and reliably even when the
distribution of the data is far from ideal for the method.
In the polynomial decay setup, the original tensor is not particularly low rank,
so even a rather expensive and accurate method (the HOOI) cannot achieve low error;
yet \cref{fig:vary-k-600,fig:vary-k-600-compare,fig:vary-memory}
tell us that the penalty from using our cheaper methods is
essentially the same regardless of the data distribution.

\subsubsection{Different dimension reduction maps perform similarly}
We first investigate the performance of our one-pass fixed-rank algorithm
as the sketch size (hence, the required storage) varies
for several types of dimension reduction maps. 
We generate synthetic data as described above with $\mathbf{r} = (5,5,5)$, $I = 600$.
\cref{fig:vary-k-600} shows the error of the rank-$\mathbf{r}$ approximation
as a function of the compression factor $k/I$.
\ifdefined \issupplement
(Results for other input tensors appear in the supplement.)
\else
(Results for other input tensors are presented as
\cref{fig:vary-k-400-app} and \cref{fig:vary-k-200-app} in \cref{appendix:more_result}.)
\fi
In general, the performance for different maps are similar,
although our theory only guarantees results for the Gaussian map.
We see that for all input tensors, the performance of our one-pass algorithm
converges to that of HOOI as $k$ increases.

\begin{figure}
	\centering
	\includegraphics[height=2.9cm]{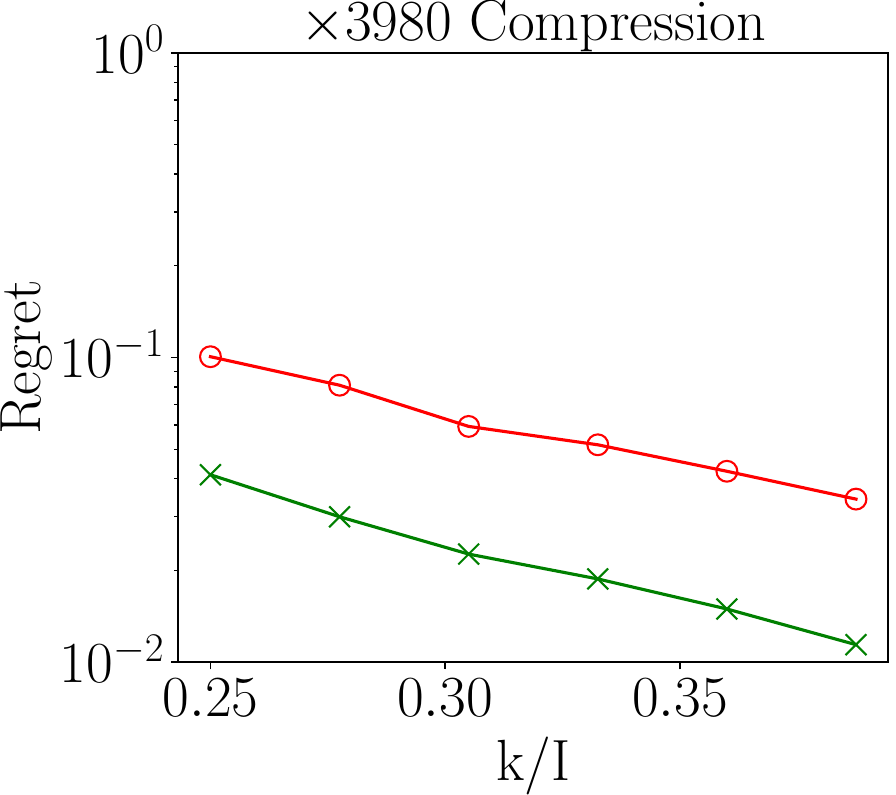}
	\includegraphics[height=2.9cm]{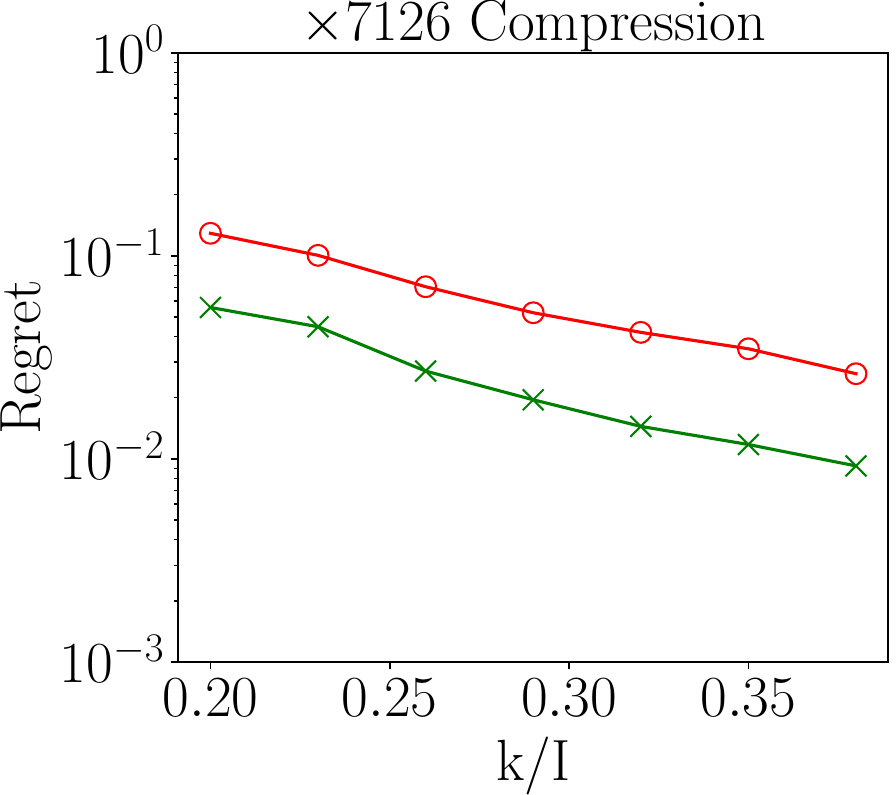}
	\includegraphics[height=2.9cm]{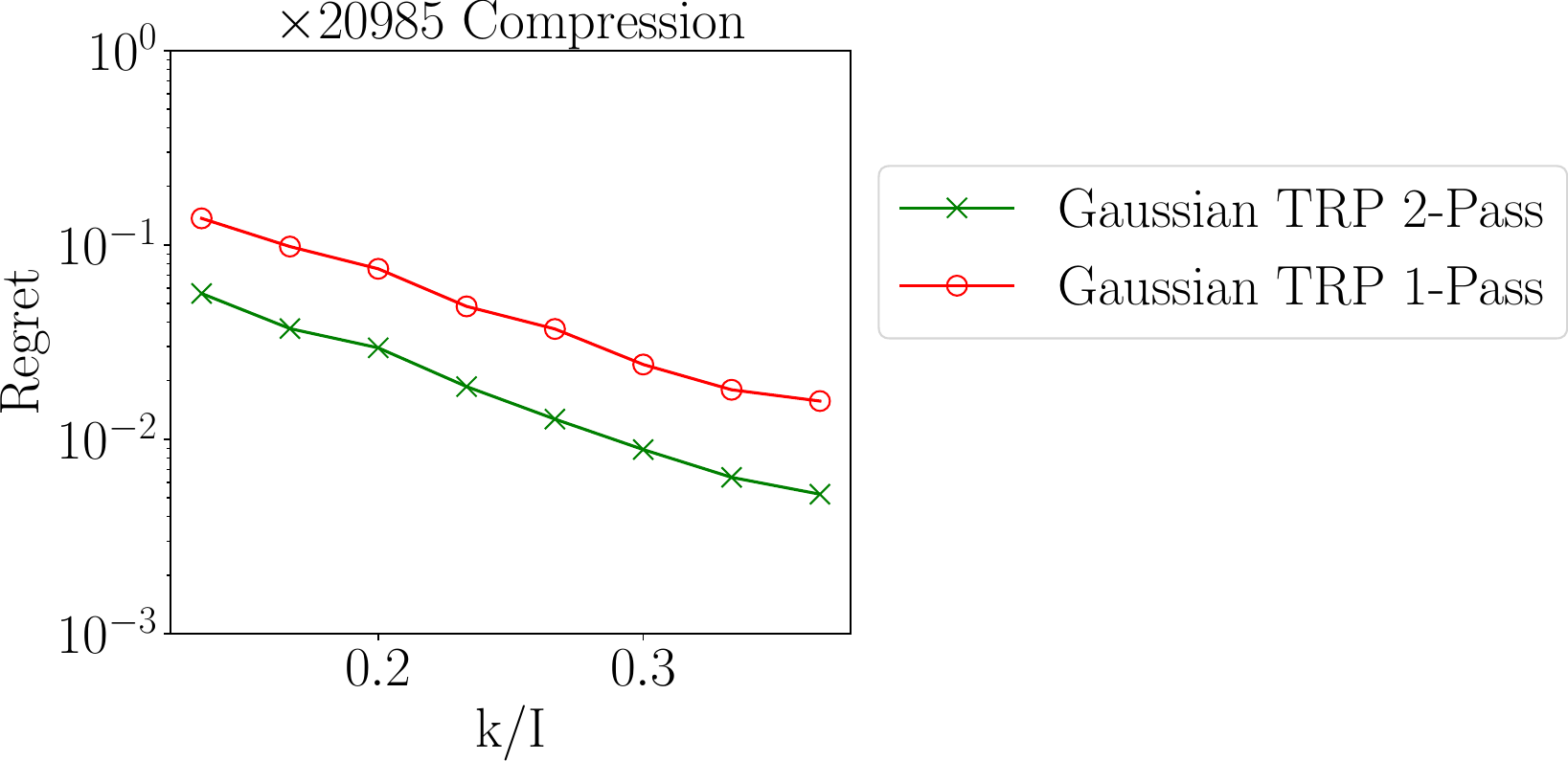}\\
	\textbf{Aerosol Absorption}\\~\\
	\includegraphics[height=2.9cm]{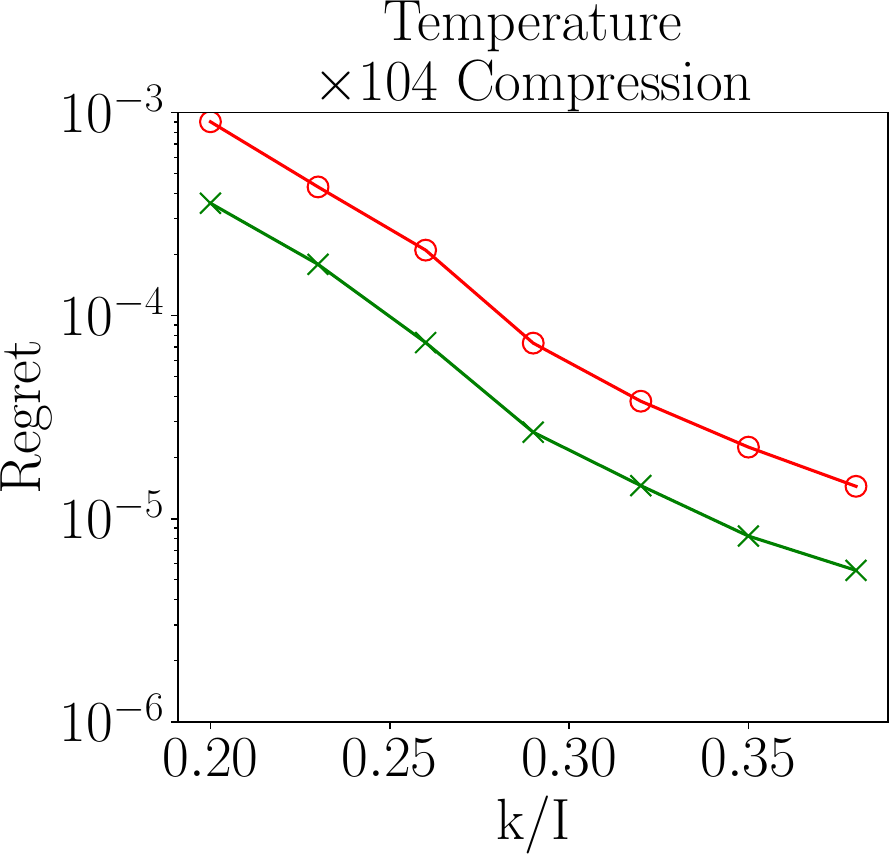}
	\includegraphics[height=2.9cm]{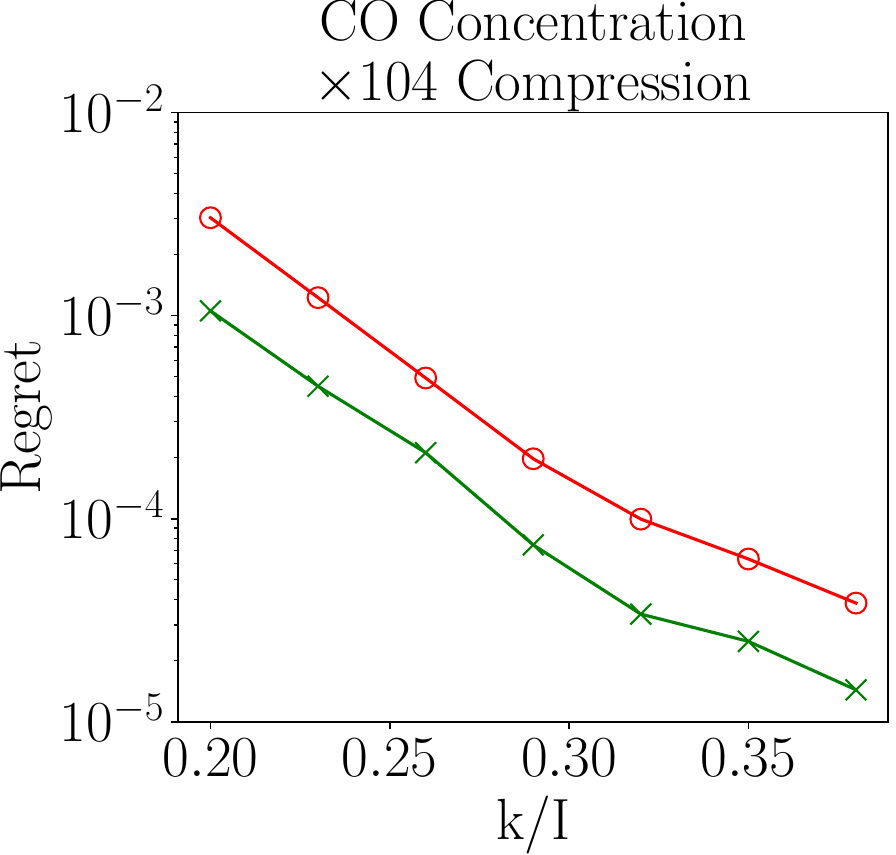}
	\includegraphics[height=2.9cm]{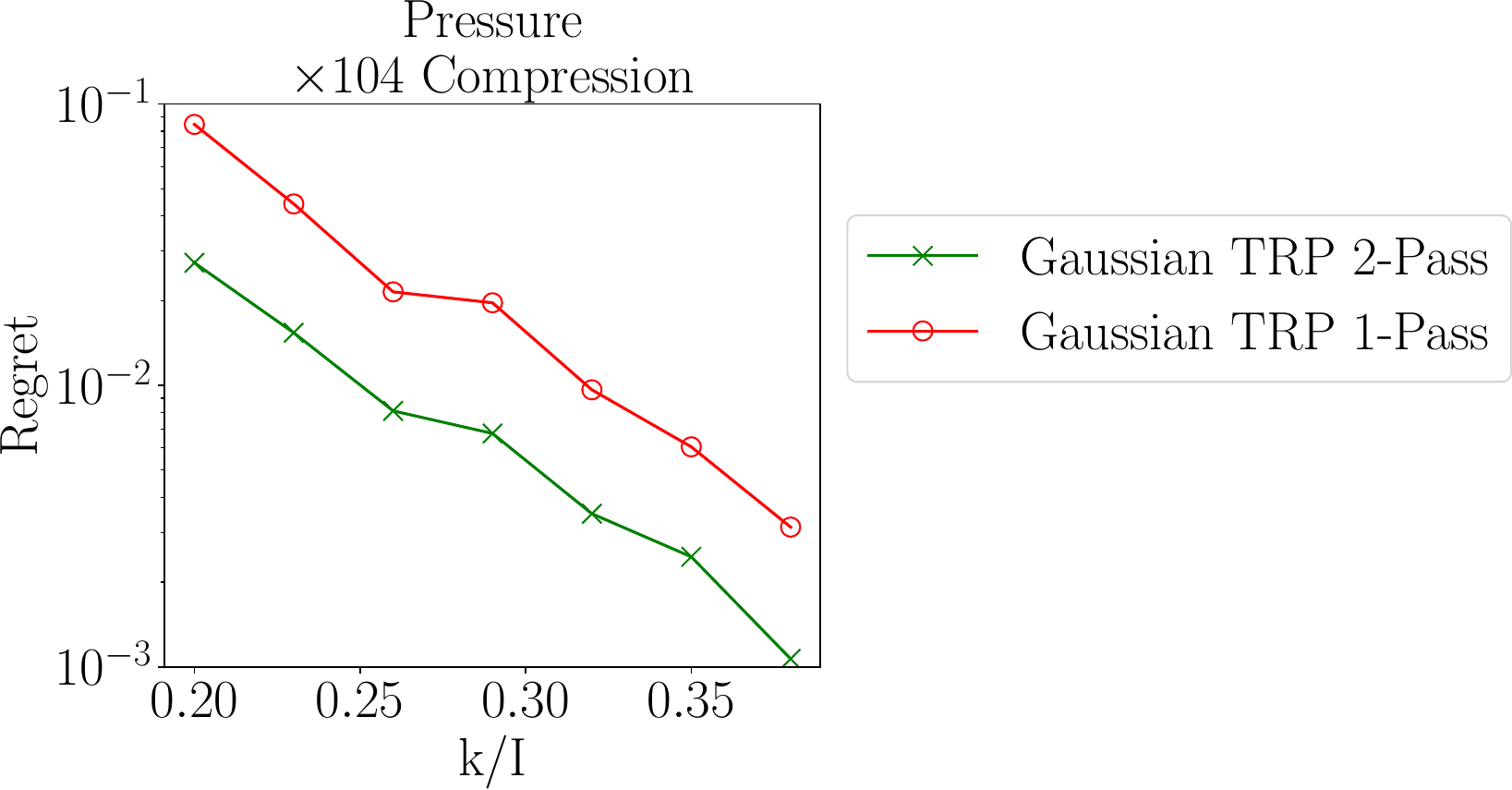}\\
	\textbf{Combustion Simulation}\\
	\caption{\textit{Approximations improve with more memory: real data.}
		We approximate aerosol absorption and combustion data
		using our one-pass and two-pass algorithms with the Gaussian TRP.
		We compare three target ranks ($r/I = 0.125,0.1,0.067$) for the former,
		and use the same target rank ($r/I = 0.1$) for each measured quantity in
		the combustion dataset.
		Notice that $r/I = 0.1$ gives a hundred-fold compression.
		For reference, on the aerosol data, the HOOI gives an approximation with relative errors
		.23, .26, and .33 for each of the three ranks, respectively;
		on the combusion data, the relative error of HOOI is
		.0063, .032, and .28 for temperature, CO, and pressure, respectively.
	}\label{fig:climate}
\end{figure}

\begin{figure}[h!]
	\includegraphics[height=3cm]{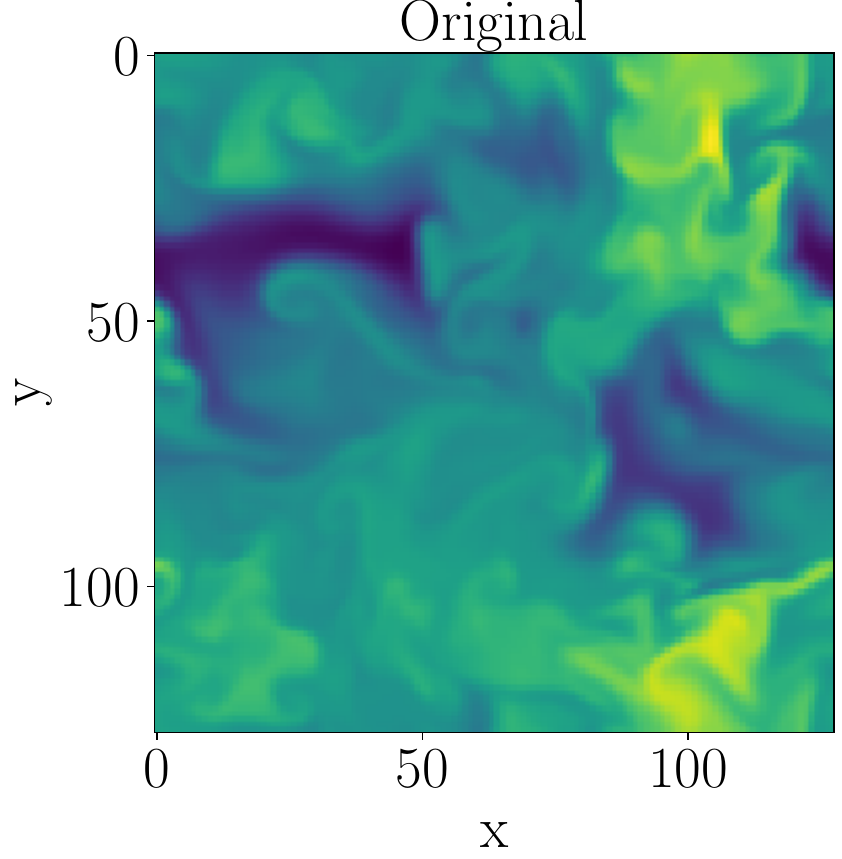}
	\includegraphics[height=3cm]{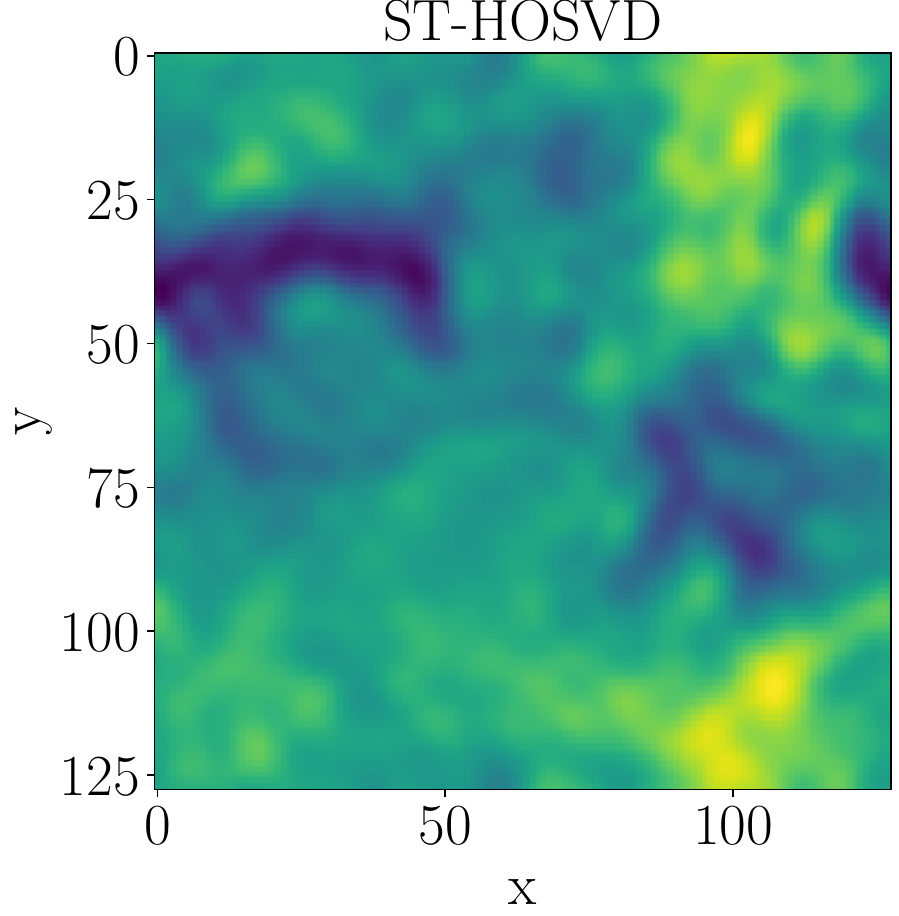}
	\includegraphics[height=3cm]{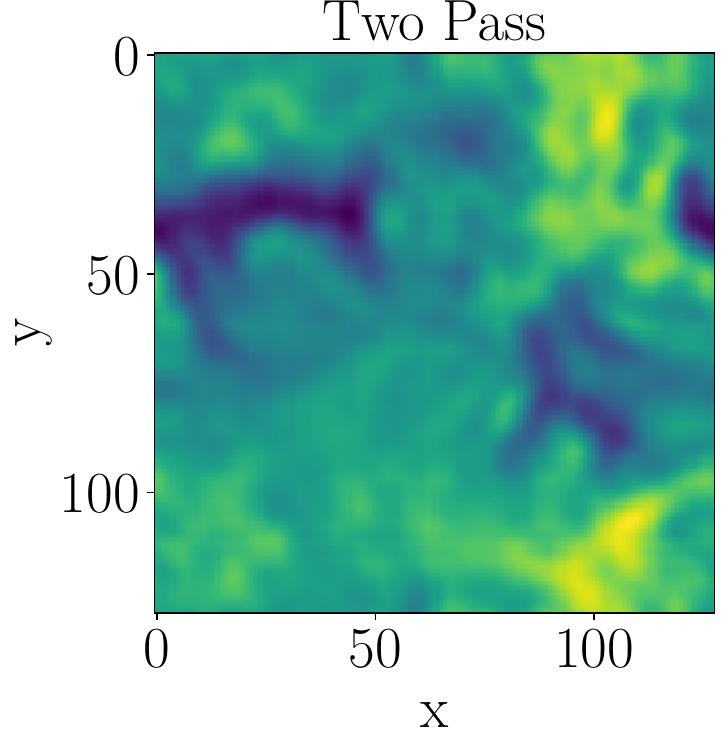}
	\includegraphics[height=3cm]{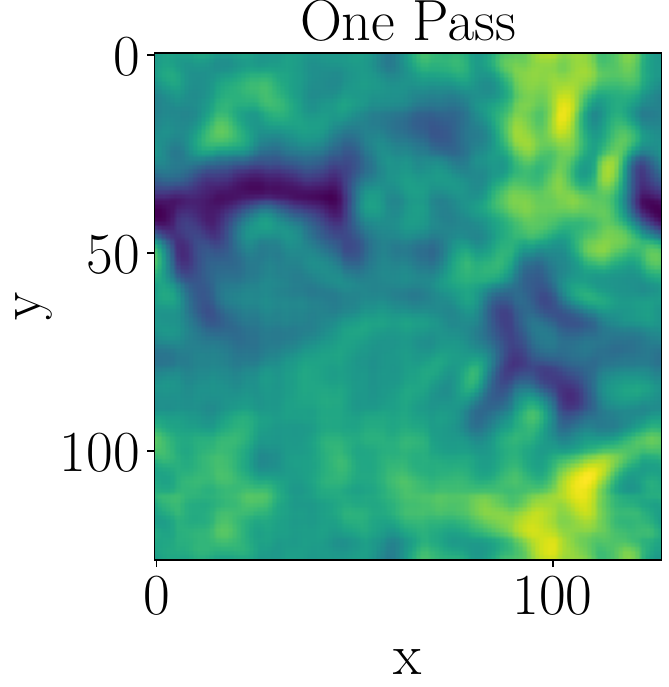}
	\centering
	\caption{\label{fig:T100}\textit{Visualizing combustion simulation:}
	All four figures show a slice of the temperature data along the first dimension.
	The approximation uses
	$\mathbf{r} = (281,25,25)$,
	$\mathbf{k} = (562,50,50)$,
	$\mathbf{s} = (1125, 101, 101)$,
	with the Gaussian TRP in the Tucker sketch.}
\end{figure}

\begin{figure}[h!]
	\includegraphics[height=7cm]{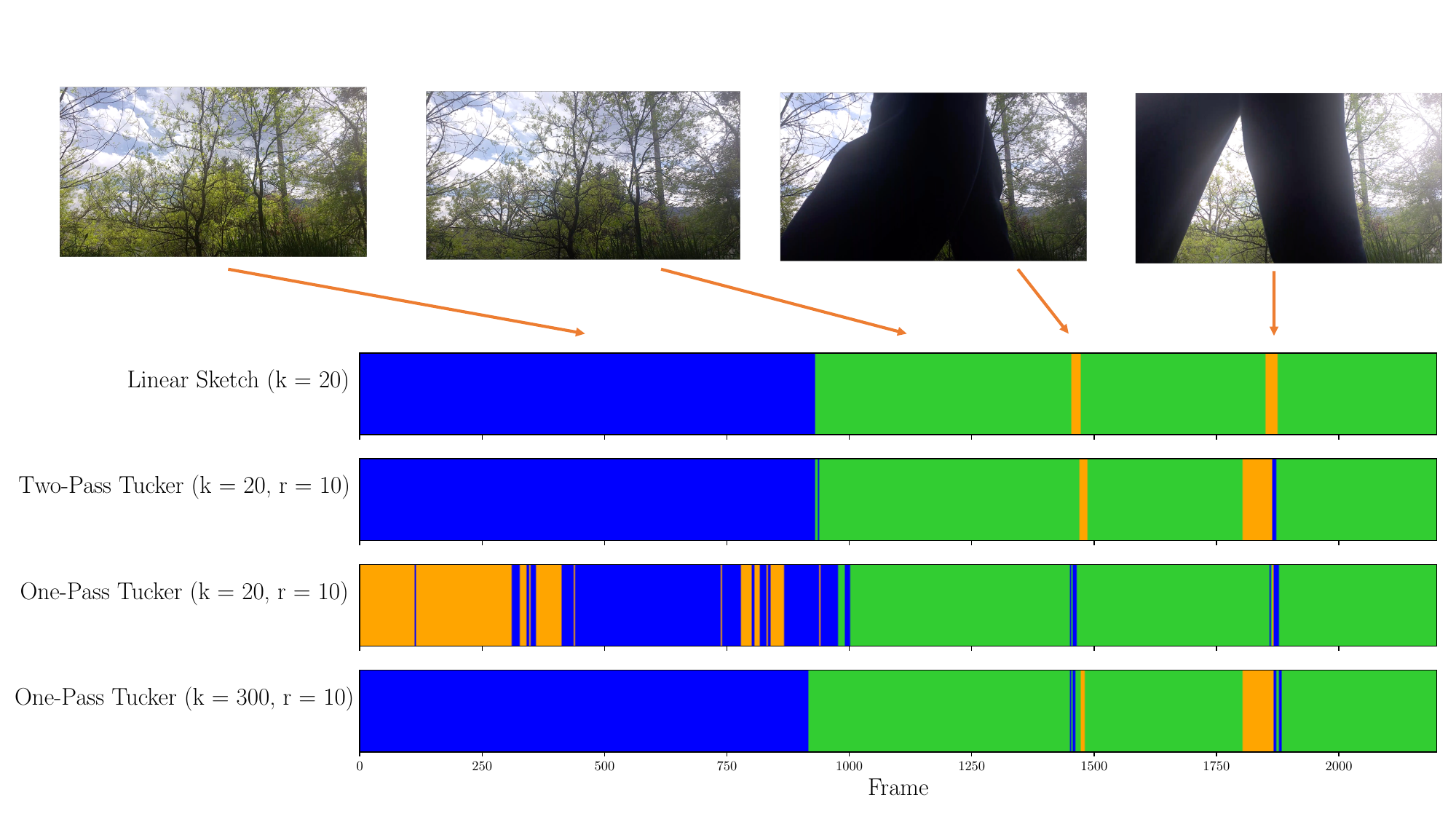} \\
	\centering
	\textbf{Video scene classification}
	\caption{\textit{Video scene classification}
		($2200 \times 1080 \times 1980$):
		We classify frames from the video data
		from \cite{malik2018low} (collected as a third order tensor with size $2200 \times 1080 \times 1980$) using $K$-means with $K$=3 on vectors computed using four different methods. $s = 2k+1$ throughout.
		(1) The linear sketch along the time dimension (row 1).
		(2-3) the Tucker factor along the time dimension,
		computed via our two-pass (row 2) and one-pass (row 3) algorithms.
		(4) The Tucker factor along the time dimension,
		computed via our one-pass (row 4) algorithm
		}\label{fig:video}
\end{figure}

\begin{figure}[h!]
	\includegraphics[height=2.4cm]{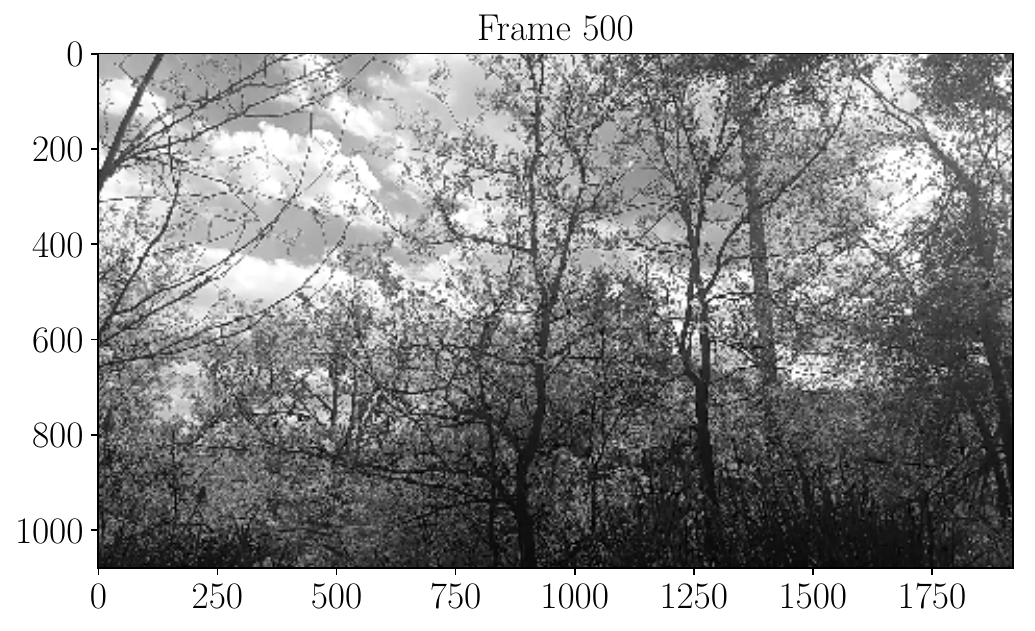}
	\includegraphics[height=2.4cm]{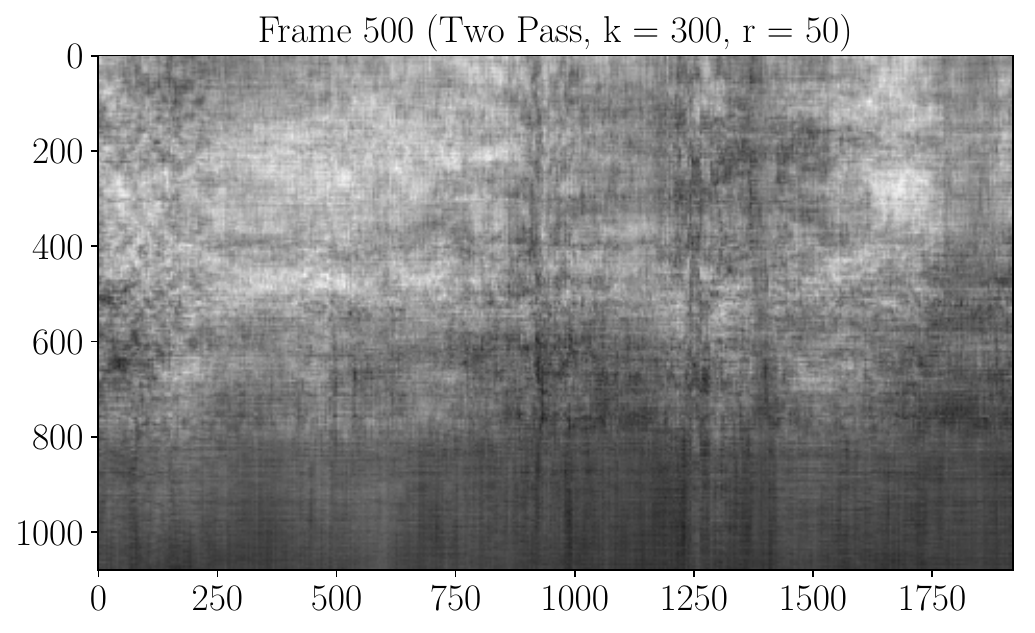}
	\includegraphics[height=2.4cm]{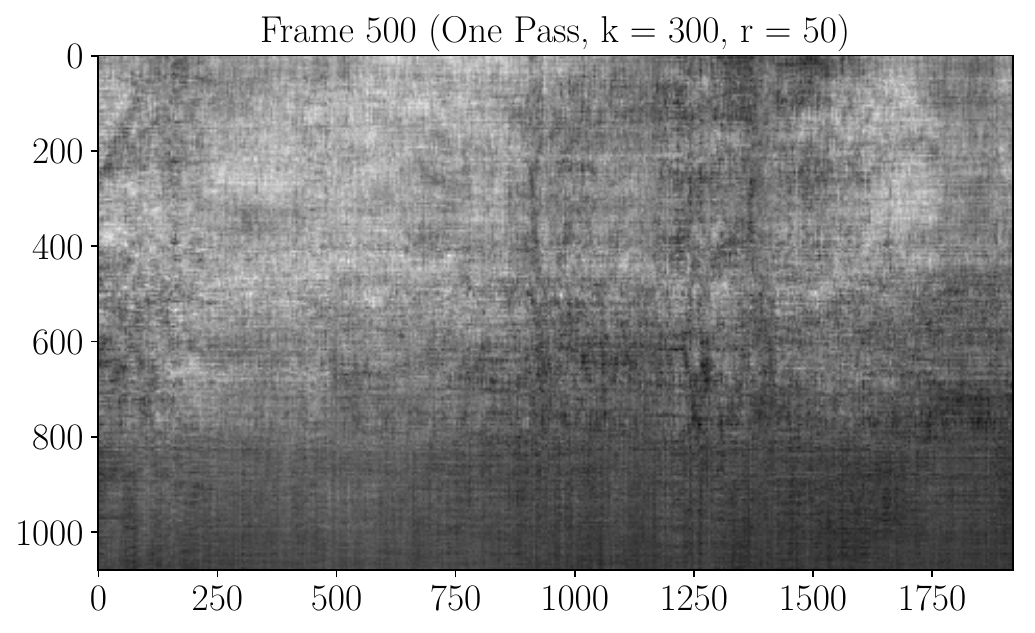}
	\centering
	\caption{\textit{Visualizing video recovery:}
	Original frame (left);
	approximation by two-pass sketch (middle);
	approximation by one-pass sketch (right).
	}\label{fig:Frame500}
\end{figure}

\subsubsection{A second pass reduces error}
The second experiment compares our two-pass and one-pass algorithms.
The design is similar to the first experiment.
\cref{fig:vary-k-600-compare} shows that the two-pass algorithm
typically outperforms the one-pass algorithm,
especially in the high-noise, sparse, or rank-decay case.
Both converge at the same asymptotic rate.
(Results for other input tensors are available
\ifdefined \issupplement
in the supplement.)
\else
in \cref{appendix:more_result}.)
\fi

\subsubsection{Improvement on state-of-the-art}
The third experiment compares the performance of our two-pass and one-pass algorithms
and Tucker TensorSketch (T.--TS), as described in \cite{malik2018low},
the only extant one-pass algorithm.
For a fair comparison, we allocate the same storage budget to each algorithm
and compare the relative error of the resulting fixed-rank approximations.
We approximate synthetic three-dimensional tensors with equal side lengths $I_1 = I_2 = I_3 = I = 300$
and of equal multilinear rank $\mathbf{r} = (r,r,r)$ with $r = 10$.
We use the suggested parameter settings for each algorithm:
$\mathbf{k} = 2\mathbf{r}+1$ and $\mathbf{s} = 2\mathbf{k}+1$ for our methods;
$K = 10$ for T.--TS.
\mnote{should this be $K=r$?}
Our one-pass algorithm
(with the Gaussian TRP)
uses $((4k+3)^N + (2r+1)IN)$ storage,
whereas T.-TS uses $(Kr^{2N}+Kr^{2N-2})$ storage
\ifdefined \issupplement
(see supplement).
\else
(see \cref{tab:storage-comparison} in \cref{appendix: time-complexity}).
\fi

\cref{fig:vary-memory} shows that our algorithms generally perform as well as T.--TS
and dramatically outperform for small storage budgets.
One nice property of our method is that the regret consistently decreases with increasing storage.
In contrast, the tensor sketch method behaves unpredictably as storage increases: there are wide plateaus
where increasing storage hardly helps at all, and occasionally, increasing storage hurts performance.
The performance of T.-TS is comparable with that of
the algorithms presented in this paper only when the storage budget is large.

\begin{remark}
	The paper \cite{malik2018low} proposes a multi-pass method, Tucker Tensor-Times-Matrix-TensorSketch (TTMTS) that is dominated by the one-pass method Tucker TensorSketch(TS) in all numerical experiments;
  hence we compare only with T.-TS.
\end{remark}

\subsection{Applications}\label{s-real-data}

We also apply our method to datasets drawn from three application domains:
climate, combustion, and video.
\begin{itemize}
\item \emph{Climate data.}
We consider global climate simulation datasets from
the Community Earth System Model (CESM) Community Atmosphere Model (CAM) 5.0 \cite{hurrell2013community,kay2015community}.
The dataset on aerosol absorption has four dimensions:
times, altitudes, longitudes, and latitudes  ($240 \times 30 \times 192 \times 288$).
The data on net radiative flux at surface and dust aerosol burden have three dimensions:
times, longitudes, and latitudes ($1200 \times 192 \times 288$).
Each of these quantitives has a strong impact on the absorption of solar radiation and on cloud formation.

\item \emph{Combustion data.}
We consider combustion simulation data from \cite{lapointe2015differential}.
The data consists of three measured quantities
(pressure, CO concentration, and temperature) 
each observed on a $1408 \times 128 \times 128$ spatial grid.

\item \emph{Video data.}
Consider the three-dimensional tensor from \cite{malik2018low}: each slice of the tensor is a video frame.
A low frame rate camera is mounted in a fixed position as people walk by to form the video,
which consists of 2493 frames, each of size 1080 by 1980.
Stored as a \texttt{numpy.array}, the video data is 41.4 GB in total.
\end{itemize}

\subsubsection{Data compression}
We show that our proposed algorithms are able to
successfully compress climate and combustion data
even when the full data does not fit in memory.
Since the multilinear rank of the original tensor is unknown, we perform experiments for
three different target ranks. In this experiment, we hope to understand the effect of different choices of storage budget $k$ to
 achieve the same compression ratio. We define the compression ratio
 as the ratio in size between the original input tensor and the output Tucker factors, i.e. $\frac{\prod_{i = 1}^N I_i}{\sum_{i=1}^Nr_iI_i+ \prod_{i = 1}^N r_i}$.
As in our experiments on simulated data, \cref{fig:climate} shows
that the two-pass algorithm outperforms the one-pass algorithm as expected.
However, as the storage budget $k$ increases, both methods converge to the performance of HOOI.
The rate of convergence is faster for smaller target ranks.
Performance of our algorithms on the combustion simulation is qualitatively similar
but converges faster to the performance of HOOI. \cref{fig:T100} visualizes the recovery of
the temperature data in combustion simulation for a slice along the first dimension.
We observe that the recovery for both two-pass and one-pass algorithms approximate the recovery from HOOI.
\ifdefined \issupplement
Similar results on other datasets appear in the supplement.)
\else
\cref{fig:srfrad_burden_dust} in \cref{appendix:more_real_data_result}
shows similar results on another dataset.
\fi

\subsubsection{Video scene classification}
We show how to use our single-pass method to classify scenes in the video data described above.
The goal is to identify frames in which people appear.
We remove the first 100 frames and last 193 frames where the camera setup happened,
as in \cite{malik2018low}.
We stream over the tensor and sketch it using parameters $k = 300, s = 601$.
Finally, we compute a fixed-rank approximation with $\mathbf{r} = (10,10,10)$ and $(20,20,20)$.
We apply K-means clustering to the resulting 10- or 20-dimensional vectors
corresponding to each of the remaining 2200 frames.

We experimented with clustering vectors found in three ways:
from the unfolding along the time dimension after two-pass or one-pass Tucker approximations,
or directly from the factor sketch along the time dimension, which we call the linear sketch.
In \cref{fig:video}, comparing the video frames with the classification results,
we can see that the background lighting is relatively dark at the beginning,
and initial frames are classified into \texttt{Class} $0$.
After a change in the background lighting,
most other frames of the video are classified into \texttt{Class} $1$.
When a person passes by the camera, the frames are classified into \texttt{Class} $2$.
Right after the person passes by, the frames are classified into \texttt{Class} $0$,
the brighter background scene, due to the light adjustment.

Our classification results (using the linear sketch or approximation)
are similar to those in \cite{malik2018low}
while using only $1/500$ as much storage; the one-pass approximation
requires more storage (but still less than \cite{malik2018low}) 
to achieve similar performance.
In particular, using the sketch itself, rather than the Tucker approximation,
to summarize the data enables very efficient video scene classification.
Interestingly, classification works well even though
the video is not very low rank along the spatial dimensions.
\cref{fig:Frame500} shows that the scene is poorly approximated
even with $\mathbf{s} = {601, 601, 601}$, $\mathbf{k} = (300, 300, 300)$, and $\mathbf{r} = (50, 50, 50)$.

\section*{Acknowledgments}
MU, YS, and YG were supported in part by DARPA Award FA8750-17-2-0101,
NSF Awards IIS-1943131 
and CCF-1740822,
the ONR Young Investigator Program, %
and the Simons Institute.
JAT gratefully acknowledges support from ONR Awards N00014-11-10025, N00014-17-12146, and N00014-18-12363.
The authors wish to thank Osman Asif Malik and Stephen Becker for their help
in understanding and implementing Tucker TensorSketch,
and Tamara Kolda and two anonymous reviewers for insightful comments and suggestions
that helped to improve this manuscript.

\clearpage
\bibliographystyle{siamplain}
\bibliography{bibtex}

\begin{thebibliography}{10}

\bibitem{achlioptas2003database}
{\sc D.~Achlioptas}, {\em Database-friendly random projections:
  {Johnson-Lindenstrauss} with binary coins}, Journal of computer and System
  Sciences, 66 (2003), pp.~671--687.

\bibitem{ailon2009fast}
{\sc N.~Ailon and B.~Chazelle}, {\em The fast {Johnson--Lindenstrauss}
  transform and approximate nearest neighbors}, SIAM Journal on computing, 39
  (2009), pp.~302--322.

\bibitem{arora2009computational}
{\sc S.~Arora and B.~Barak}, {\em Computational complexity: a modern approach},
  Cambridge University Press, 2009.

\bibitem{austin2016parallel}
{\sc W.~Austin, G.~Ballard, and T.~G. Kolda}, {\em Parallel tensor compression
  for large-scale scientific data}, in Parallel and Distributed Processing
  Symposium, 2016 IEEE International, IEEE, 2016, pp.~912--922.

\bibitem{ballester2019tthresh}
{\sc R.~Ballester-Ripoll, P.~Lindstrom, and R.~Pajarola}, {\em {TTHRESH}:
  Tensor compression for multidimensional visual data}, IEEE transactions on
  visualization and computer graphics,  (2019).

\bibitem{baskaran2012efficient}
{\sc M.~Baskaran, B.~Meister, N.~Vasilache, and R.~Lethin}, {\em Efficient and
  scalable computations with sparse tensors}, in High Performance Extreme
  Computing (HPEC), 2012 IEEE Conference on, IEEE, 2012, pp.~1--6.

\bibitem{battaglino2018practical}
{\sc C.~Battaglino, G.~Ballard, and T.~G. Kolda}, {\em A practical randomized
  cp tensor decomposition}, SIAM Journal on Matrix Analysis and Applications,
  39 (2018), pp.~876--901.

\bibitem{battaglino2019faster}
{\sc C.~Battaglino, G.~Ballard, and T.~G. Kolda}, {\em Faster parallel tucker
  tensor decomposition using randomization},  (2019).

\bibitem{boutsidis2013improved}
{\sc C.~Boutsidis and A.~Gittens}, {\em Improved matrix algorithms via the
  subsampled randomized hadamard transform}, SIAM Journal on Matrix Analysis
  and Applications, 34 (2013), pp.~1301--1340.

\bibitem{breiding2018riemannian}
{\sc P.~Breiding and N.~Vannieuwenhoven}, {\em A riemannian trust region method
  for the canonical tensor rank approximation problem}, SIAM Journal on
  Optimization, 28 (2018), pp.~2435--2465.

\bibitem{cichocki2013tensor}
{\sc A.~Cichocki}, {\em Tensor decompositions: a new concept in brain data
  analysis?}, arXiv preprint arXiv:1305.0395,  (2013).

\bibitem{clarkson2017low}
{\sc K.~L. Clarkson and D.~P. Woodruff}, {\em Low-rank approximation and
  regression in input sparsity time}, Journal of the ACM (JACM), 63 (2017),
  p.~54.

\bibitem{cormode2008finding}
{\sc G.~Cormode and M.~Hadjieleftheriou}, {\em Finding frequent items in data
  streams}, Proceedings of the VLDB Endowment, 1 (2008), pp.~1530--1541.

\bibitem{de2000multilinear}
{\sc L.~De~Lathauwer, B.~De~Moor, and J.~Vandewalle}, {\em A multilinear
  singular value decomposition}, SIAM journal on Matrix Analysis and
  Applications, 21 (2000), pp.~1253--1278.

\bibitem{de2000best}
{\sc L.~De~Lathauwer, B.~De~Moor, and J.~Vandewalle}, {\em On the best rank-1
  and rank-(r1, r2, ..., rn) approximation of higher-order tensors}, SIAM
  Journal on Matrix Analysis and Applications, 21 (2000), pp.~1324--1342,
  \url{https://doi.org/10.1137/S0895479898346995},
  \url{https://doi.org/10.1137/S0895479898346995},
  \url{https://arxiv.org/abs/https://doi.org/10.1137/S0895479898346995}.

\bibitem{de2008tensor}
{\sc V.~De~Silva and L.-H. Lim}, {\em Tensor rank and the ill-posedness of the
  best low-rank approximation problem}, SIAM Journal on Matrix Analysis and
  Applications, 30 (2008), pp.~1084--1127.

\bibitem{2017arXiv171209473D}
{\sc H.~{Diao}, Z.~{Song}, W.~{Sun}, and D.~P. {Woodruff}}, {\em {Sketching for
  Kronecker Product Regression and P-splines}}, arXiv e-prints,  (2017),
  arXiv:1712.09473, p.~arXiv:1712.09473,
  \url{https://arxiv.org/abs/1712.09473}.

\bibitem{grasedyck2010hierarchical}
{\sc L.~Grasedyck}, {\em Hierarchical singular value decomposition of tensors},
  SIAM Journal on Matrix Analysis and Applications, 31 (2010), pp.~2029--2054.

\bibitem{gu1996efficient}
{\sc M.~Gu and S.~C. Eisenstat}, {\em Efficient algorithms for computing a
  strong rank-revealing qr factorization}, SIAM Journal on Scientific
  Computing, 17 (1996), pp.~848--869.

\bibitem{hackbusch2012tensor}
{\sc W.~Hackbusch}, {\em Tensor spaces and numerical tensor calculus}, vol.~42,
  Springer Science \& Business Media, 2012.

\bibitem{halko2011finding}
{\sc N.~Halko, P.-G. Martinsson, and J.~A. Tropp}, {\em Finding structure with
  randomness: Probabilistic algorithms for constructing approximate matrix
  decompositions}, SIAM review, 53 (2011), pp.~217--288.

\bibitem{hurrell2013community}
{\sc J.~W. Hurrell, M.~M. Holland, P.~R. Gent, S.~Ghan, J.~E. Kay, P.~J.
  Kushner, J.-F. Lamarque, W.~G. Large, D.~Lawrence, K.~Lindsay, et~al.}, {\em
  The community earth system model: a framework for collaborative research},
  Bulletin of the American Meteorological Society, 94 (2013), pp.~1339--1360.

\bibitem{kay2015community}
{\sc J.~Kay, C.~Deser, A.~Phillips, A.~Mai, C.~Hannay, G.~Strand, J.~Arblaster,
  S.~Bates, G.~Danabasoglu, J.~Edwards, et~al.}, {\em The community earth
  system model (cesm) large ensemble project: A community resource for studying
  climate change in the presence of internal climate variability}, Bulletin of
  the American Meteorological Society, 96 (2015), pp.~1333--1349.

\bibitem{kaya2016high}
{\sc O.~Kaya and B.~U{\c{c}}ar}, {\em High performance parallel algorithms for
  the tucker decomposition of sparse tensors}, in Parallel Processing (ICPP),
  2016 45th International Conference on, IEEE, 2016, pp.~103--112.

\bibitem{kolda2009tensor}
{\sc T.~G. Kolda and B.~W. Bader}, {\em Tensor decompositions and
  applications}, SIAM review, 51 (2009), pp.~455--500.

\bibitem{kolda2008scalable}
{\sc T.~G. Kolda and J.~Sun}, {\em Scalable tensor decompositions for
  multi-aspect data mining}, in 2008 Eighth IEEE International Conference on
  Data Mining, IEEE, 2008, pp.~363--372.

\bibitem{kossaifi2019tensorly}
{\sc J.~Kossaifi, Y.~Panagakis, A.~Anandkumar, and M.~Pantic}, {\em Tensorly:
  Tensor learning in python}, The Journal of Machine Learning Research, 20
  (2019), pp.~925--930.

\bibitem{lapointe2015differential}
{\sc S.~Lapointe, B.~Savard, and G.~Blanquart}, {\em Differential diffusion
  effects, distributed burning, and local extinctions in high karlovitz
  premixed flames}, Combustion and flame, 162 (2015), pp.~3341--3355.

\bibitem{li2015input}
{\sc J.~Li, C.~Battaglino, I.~Perros, J.~Sun, and R.~Vuduc}, {\em An
  input-adaptive and in-place approach to dense tensor-times-matrix multiply},
  in High Performance Computing, Networking, Storage and Analysis, 2015
  SC-International Conference for, IEEE, 2015, pp.~1--12.

\bibitem{li2006very}
{\sc P.~Li, T.~J. Hastie, and K.~W. Church}, {\em Very sparse random
  projections}, in Proceedings of the 12th ACM SIGKDD international conference
  on Knowledge discovery and data mining, ACM, 2006, pp.~287--296.

\bibitem{malik2018low}
{\sc O.~A. Malik and S.~Becker}, {\em Low-rank tucker decomposition of large
  tensors using tensorsketch}, in Advances in Neural Information Processing
  Systems, 2018, pp.~10116--10126.

\bibitem{minster2019randomized}
{\sc R.~Minster, A.~K. Saibaba, and M.~E. Kilmer}, {\em Randomized algorithms
  for low-rank tensor decompositions in the tucker format}, arXiv preprint
  arXiv:1905.07311,  (2019).

\bibitem{muthukrishnan2005data}
{\sc S.~Muthukrishnan et~al.}, {\em Data streams: Algorithms and applications},
  Foundations and Trends{\textregistered} in Theoretical Computer Science, 1
  (2005), pp.~117--236.

\bibitem{oymak2015universality}
{\sc S.~Oymak and J.~A. Tropp}, {\em Universality laws for randomized dimension
  reduction, with applications}, Information and Inference: A Journal of the
  IMA,  (2015).

\bibitem{rudelson2012row}
{\sc M.~Rudelson}, {\em Row products of random matrices}, Advances in
  Mathematics, 231 (2012), pp.~3199--3231.

\bibitem{sun2008incremental}
{\sc J.~Sun, D.~Tao, S.~Papadimitriou, P.~S. Yu, and C.~Faloutsos}, {\em
  Incremental tensor analysis: Theory and applications}, ACM Transactions on
  Knowledge Discovery from Data (TKDD), 2 (2008), p.~11.

\bibitem{sun2018tensor}
{\sc Y.~Sun, Y.~Guo, J.~A. Tropp, and M.~Udell}, {\em Tensor random projection
  for low memory dimension reduction}, in NeurIPS Workshop on Relational
  Representation Learning, 2018,
  \url{https://r2learning.github.io/assets/papers/CameraReadySubmission%2041.pdf}.

\bibitem{tropp2011improved}
{\sc J.~A. Tropp}, {\em Improved analysis of the subsampled randomized hadamard
  transform}, Advances in Adaptive Data Analysis, 3 (2011), pp.~115--126.

\bibitem{tropp2017practical}
{\sc J.~A. Tropp, A.~Yurtsever, M.~Udell, and V.~Cevher}, {\em Practical
  sketching algorithms for low-rank matrix approximation}, SIAM Journal on
  Matrix Analysis and Applications, 38 (2017), pp.~1454--1485.

\bibitem{tropp2018more}
{\sc J.~A. Tropp, A.~Yurtsever, M.~Udell, and V.~Cevher}, {\em More practical
  sketching algorithms for low-rank matrix approximation}, Tech. Report
  2018-01, California Institute of Technology, Pasadena, California, 2018.

\bibitem{tropp2019streaming}
{\sc J.~A. Tropp, A.~Yurtsever, M.~Udell, and V.~Cevher}, {\em Streaming
  low-rank matrix approximation with an application to scientific simulation},
  {SIAM} Journal on Scientific Computing {(SISC)},  (2019),
  \url{https://arxiv.org/abs/1902.08651}.

\bibitem{tsourakakis2010mach}
{\sc C.~E. Tsourakakis}, {\em Mach: Fast randomized tensor decompositions}, in
  Proceedings of the 2010 SIAM International Conference on Data Mining, SIAM,
  2010, pp.~689--700.

\bibitem{tucker1966some}
{\sc L.~R. Tucker}, {\em Some mathematical notes on three-mode factor
  analysis}, Psychometrika, 31 (1966), pp.~279--311.

\bibitem{vannieuwenhoven2012new}
{\sc N.~Vannieuwenhoven, R.~Vandebril, and K.~Meerbergen}, {\em A new
  truncation strategy for the higher-order singular value decomposition}, SIAM
  Journal on Scientific Computing, 34 (2012), pp.~A1027--A1052.

\bibitem{vasilescu2002multilinear}
{\sc M.~A.~O. Vasilescu and D.~Terzopoulos}, {\em Multilinear analysis of image
  ensembles: Tensorfaces}, in European Conference on Computer Vision, Springer,
  2002, pp.~447--460.

\bibitem{wang2015fast}
{\sc Y.~Wang, H.-Y. Tung, A.~J. Smola, and A.~Anandkumar}, {\em Fast and
  guaranteed tensor decomposition via sketching}, in Advances in Neural
  Information Processing Systems, 2015, pp.~991--999.

\bibitem{woodruff2014sketching}
{\sc D.~P. Woodruff et~al.}, {\em Sketching as a tool for numerical linear
  algebra}, Foundations and Trends{\textregistered} in Theoretical Computer
  Science, 10 (2014), pp.~1--157.

\bibitem{woolfe2008fast}
{\sc F.~Woolfe, E.~Liberty, V.~Rokhlin, and M.~Tygert}, {\em A fast randomized
  algorithm for the approximation of matrices}, Applied and Computational
  Harmonic Analysis, 25 (2008), pp.~335--366.

\bibitem{zhou2014decomposition}
{\sc G.~Zhou, A.~Cichocki, and S.~Xie}, {\em Decomposition of big tensors with
  low multilinear rank}, arXiv preprint arXiv:1412.1885,  (2014).

\end{thebibliography}

\appendix
\begin{appendices}

\section{Probabilistic Analysis of Core Sketch Error}
This section contains the most technical part of our proof.
We provide a probabilistic error bound for the difference between the
two-pass core approximation $\T{W}_2$
from \cref{alg:two_pass_low_rank_appro}
and the one-pass core approximation $\T{W}_1$
from \cref{alg:one_pass_low_rank_appro}.

Introduce for each $n\in[N]$ the orthonormal matrix $\M{Q}_n^\bot$ that forms
a basis for the subspace orthogonal to $\M{Q}_n$, so that
$\M{Q}_n^\bot (\M{Q}_n^\bot)^\top = \M{I} - \M{Q}_n\M{Q}_n^\top$.
Next, define
\begin{equation}
\label{eq: def-proj-Q}
\begin{aligned}
\M{\Phi}_n^Q = \M{\Phi}^\top_n \M{Q}_n  ,~~~~\M{\Phi}_n^{Q^\bot} = \M{\Phi}^\top_n \M{Q}_n^\bot.
\end{aligned}
\end{equation}
Recall that the DRMs $\M{\Phi}_n$ are i.i.d. Gaussian. Thus, conditional on $\M{Q}_n$,
the random matrices $\M{\Phi}_n^Q$ and $\M{\Phi}_n^{Q^\bot}$ are statistically independent.

\subsection{Decomposition of Core Approximation Error}
In this section, we characterize the
difference between the one- and two-pass core approximations
$\T{W}_1-\T{W}_2 = \T{W}_1 - \T{X}\times_1 \M{Q}_1^\top \dots \times_N \M{Q}_N^\top$.
\begin{lem}
\label{lemma:core_error_decomposition}
Suppose that $\M{\Phi}_n$ has full column rank for each $n \in [N]$.
We define $\mathbbm{1}_{a = b} = 1$ if $a=b$ and 0 otherwise.
Then
\begin{equation}
\T{W}_1-\T{W}_2 = \T{W}_1 - \T{X}\times_1 \M{Q}_1^\top \dots \times_N \M{Q}_N^\top =
\sum_{(i_1,\dots, i_N) \in \{0,1\}^N, \sum_{j=1}^N i_j \geq 1} \T{Y}_{i_1\dots i_N}, \nonumber
\end{equation}
where
\begin{equation}
\label{eq:def_each_part}
\begin{aligned}
\T{Y}_{i_1\dots i_N} &= \T{X}\times_1 \left(\mathbbm{1}_{i_1=0}\M{Q}_1^\top + \mathbbm{1}_{i_1=1}(\M{\Phi}_1^{Q_1})^\dag  \M{\Phi}_1^{Q_1^\bot}(\M{Q}_1^\bot)^\top \right)\\
&\times_2 \cdots \times_N \left(\mathbbm{1}_{i_N=0}\M{Q}_N^\top + \mathbbm{1}_{i_1=1}(\M{\Phi}_N^{Q_N})^\dag  \M{\Phi}_N^{Q_N^\bot}(\M{Q}_N^\bot)^\top \right).
\end{aligned}
\end{equation}
\end{lem}
\begin{proof}
Let $\T{H}$ be the core sketch from \cref{alg:tensor_sketch}.
Write 
$\T{W}_1$ as
\begin{equation}
\begin{aligned}
\T{W}_1 &= \T{H}\times_1 (\M{\Phi}_1^\top \M{Q}_1)^\dag \times_2 \cdots \times_N (\M{\Phi}^\top_N \M{Q}_N)^\dag \\
&= (\T{X} -  \hat{\T{X}}_2)\times_1 \M{\Phi}^\top_1 \times_2 \cdots \times_N \M{\Phi}^\top_N  \times_1 (\M{\Phi}^\top_1 \M{Q}_1)^\dag \times_2 \cdots \times_N (\M{\Phi}^\top_N \M{Q}_N)^\dag
\\
&+ \hat{\T{X}}_2\times_1 \M{\Phi}^\top_1 \times_2 \cdots \times_N \M{\Phi}^\top_N \times_1 (\M{\Phi}^\top_1 \M{Q}_1)^\dag \times_2 \cdots \times_N (\M{\Phi}^\top_N \M{Q}_N)^\dag.\nonumber
\end{aligned}
\end{equation}
Using the fact that $(\M{\Phi}^\top_n \M{Q}_n)^\dag (\M{\Phi}^\top_n \M{Q}_n) =\M{I}$, we can simplify the second term as
\begin{equation}
\begin{aligned}
&\hat{\T{X}}_2\times_1 \M{\Phi}^\top_1 \times_2 \cdots \times_N \M{\Phi}^\top_N \times_1 (\M{\Phi}^\top_1 \M{Q}_1)^\dag \times_2 \cdots \times_N (\M{\Phi}^\top_N \M{Q}_N)^\dag   \\
& = \T{X}\times_1 (\M{\Phi}^\top_1 \M{Q}_1)^\dag \M{\Phi}^\top_1\M{Q}_1\M{Q}_1^\top \times_2 \cdots \times_N (\M{\Phi}_N^\top \M{Q}_N)^\dag \M{\Phi}_N^\top\M{Q}_N\M{Q}_N^\top\\
& = \T{X}\times_1 \M{Q}_1^\top \times_2 \cdots \times_N \M{Q}_N^\top, \nonumber
\end{aligned}
\end{equation}
which is exactly the two-pass core approximation $\T{W}_2$.
Therefore
\begin{equation}
\begin{aligned}
&\T{W}_1 -\T{W}_2= (\T{X} -  \hat{\T{X}}_2)\times_1 \M{\Phi}^\top_1 \times_2 \cdots \times_N \M{\Phi}^\top_N  \times_1 (\M{\Phi}^\top_1 \M{Q}_1)^\dag \times_2 \cdots \times_N (\M{\Phi}^\top_N \M{Q}_N)^\dag. \nonumber
\end{aligned}
\end{equation}
We continue to simplify this difference:
\begin{equation}
\begin{aligned}
(\T{X} -  \tilde{\T{X}})&\times_1 \M{\Phi}^\top_1 \times_2 \cdots \times_N \M{\Phi}^\top_N  \times_1 (\M{\Phi}^\top_1 \M{Q}_1)^\dag \times_2\cdots \times_N (\M{\Phi}^\top_N \M{Q}_N)^\dag \label{eq:core_err_decom} \\
& =(\T{X} -  \tilde{\T{X}})\times_1 (\M{\Phi}^\top_1\M{Q}_1)^\dag \M{\Phi}^\top_1 \times_2\cdots \times_N (\M{\Phi}^\top_N\M{Q}_N)^\dag \M{\Phi}_N^\top \\
& =  (\T{X} -  \tilde{\T{X}}) \times_1 (\M{\Phi}^\top_1\M{Q}_1)^\dag \M{\Phi}^\top_1(\M{Q}_1\M{Q}_1^\top + \M{Q}_1^\bot (\M{Q}_1^\bot)^\top)\dots  \\
& \times_N (\M{\Phi}^\top_N\M{Q}_N)^\dag \M{\Phi}^\top_N(\M{Q}_N\M{Q}_N^\top + \M{Q}_N^\bot (\M{Q}_N^\bot)^\top)\\
& = (\T{X} -  \tilde{\T{X}}) \times_1 (\M{Q}_1^\top + (\M{\Phi}_1^Q)^\dag  \M{\Phi}_1^{Q^\bot}(\M{Q}_1^\bot)^\top)\times_2\dots \\
 &\times_N (\M{Q}_N^\top + (\M{\Phi}_N^{Q_N})^\dag  \M{\Phi}_N^{Q_N^\bot}(\M{Q}_N^\bot)^\top).
\end{aligned}
\end{equation}
Many terms in this sum are zero. We use the following two facts:
\begin{enumerate}
\item $(\T{X} - \tilde{\T{X}})\times_1 \M{Q}_1^\top\dots \times_N \M{Q}_N^\top = 0$.
\item For each $n \in [N]$, $\tilde{\T{X}}\times_n (\M{\Phi}_n^{Q_n})^\dag  \M{\Phi}_n^{Q_n^\bot}(\M{Q}_n^\bot)^\top =  0$.
\end{enumerate}
Here, $0$ denotes a tensor with all zero elements.
These facts can be obtained from the exchange rule of the mode product and the orthogonality between $\M{Q}_n^\bot$ and $\M{Q}_n$.
Using these two facts, we find that only the terms $\T{Y}_{i_1\dots i_N}$ (defined in \eqref{eq:def_each_part})
remain in the expression.
Therefore, to complete the proof, we write \eqref{eq:core_err_decom} as
\begin{equation}
\sum_{(i_1,\dots, i_N) \in \{0,1\}^N, \sum_{n=1}^N i_n\neq 0} \T{Y}_{i_1\dots i_N}.\nonumber
\end{equation}
\end{proof}

\subsection{Probabilistic Core Error Bound}
In this section, we derive a probabilistic error bound
based on the core error decomposition from  \cref{lemma:core_error_decomposition}.
\begin{lem}
\label{lemma:err_core_sketch}
Sketch the tensor $\T{X}$ using a Tucker sketch with parameters $\V{k}$ and $\V{s} > 2 \V{k}$
with i.i.d. Gaussian $\mathcal N(0,1)$ DRMs.
Define $\Delta = \max_{n=1}^N \frac{k_n}{s_n-k_n-1}$. Let $\hat{\T{X}_2}$ be the output from the two-pass low-rank approximation method (\cref{alg:two_pass_low_rank_appro}).  Then
\begin{equation}
\mathbb{E} \|\T{W}_1 - \T{X}\times_1 \M{Q}_1^\top \dots \times_N \M{Q}_N^\top\|_F^2 \le \Delta \|\T{X} - \hat{\T{X}}_2\|
\end{equation}
\end{lem}
\begin{proof}
We use the fact that the core DRMs $\{\M{\Omega}_n\}_{n \in [N]}$
are independent of the factor matrix DRMs $\{\M{\Phi}_n\}_{n \in [N]}$,
and that the randomness in
each factor matrix approximation $\M{Q}_n$
comes solely from $\M{\Omega}_n$.

For $i\in \{0,1\}^N$, define $\T{B}_{i_1\dots i_N} =$
\[
\T{X}\times_1 (\mathbbm{1}_{i_1=0}\M{Q}_1\M{Q}_1^\top + \mathbbm{1}_{i_1=1}\M{Q}_1^\bot(\M{Q}_1^\bot)^\top)\cdots\times_N(\mathbbm{1}_{i_N=0}\M{Q}_N\M{Q}_N^\top + \mathbbm{1}_{i_N=1}\M{Q}_N^\bot(\M{Q}_N^\bot)^\top).
\]
\cref{lemma:core_error_decomposition} decomposes the core error as the sum of
$\T{Y}_{i_1\cdots i_n}$ where $\sum_{n=1}^N i_n \geq 1$.
Applying \cref{lemma:expectation_inverse_gaussian} and using
the orthogonal invariance of the Frobenius norm,  we observe
\begin{equation}
\mathbb{E} \left[ \|\T{Y}_{i_1\dots i_N}\|_F^2 \mid \M{\Omega}_1 \cdots \M{\Omega}_N \right] =\left(\prod_{n=1}^N \Delta_n^{i_n}\right)
\|\T{B}_{i_1\dots i_N}\|_F^2 \le \Delta \|\T{B}_{i_1\dots i_N}\|_F^2\nonumber
\end{equation}
when $\sum_{n=1}^N i_n \geq 1$,
where $\Delta_n = \frac{k_n}{s_n-k_n-1}<1$ and $\Delta = \max_{n=1}^N \Delta_n$.

Suppose $\mathbf{q}_1, \mathbf{q}_2 \in \{0,1\}^N$ are  index (binary) vectors of length $N$.
For different indices $\mathbf{q}_1$ and $\mathbf{q}_2$, there exists some $1\le r\le N$
such that their $r$th element is different.
Without loss of generality, assume $\mathbf{q}_1(r) = 0$ and $\mathbf{q}_2(r)=1$ to see
\begin{equation}\label{eq:inner_prod2}
\langle \T{B}_{q_1}, \T{B}_{q_2}\rangle = \langle \dots \M{Q}_r^\top \M{Q}_r^\bot \dots\rangle  = 0.
\end{equation}
Similarly we can show that the inner product between $\T{Y}_{q_1}$ and $\T{Y}_{q_2}$ is zero with different $\mathbf{q}_1, \mathbf{q}_2$.  Noticing that  $\T{B}_{0,\ldots, 0} = \hat{\T{X}}_2$, we have
\begin{align*}
\|\T{X} - \hat{\T{X}}_2\|_F^2
&= \left\|\sum_{(i_1,\dots, i_N) \in \{0,1\}^N, \sum_{n=1}^N i_n \geq 1}  \T{B}_{i_1\dots i_N}\right \|_F^2
&= \sum_{\substack{(i_1,\dots, i_N) \in \{0,1\}^N, \\ \sum_{n=1}^N i_n \geq 1}} \|\T{B}_{i_1\dots i_N}\|_F^2.
\end{align*}
Put these together and use the Pythagorean theorem
to finish the proof: 
\begin{align*}
&\mathbb{E}\left[ \|\T{W} - \T{X}\times_1 \M{Q}_1^\top \dots \times_N \M{Q}_N^\top\|_F^2 \mid \M{\Omega}_1, \cdots, \M{\Omega}_N \right] \\
& = \sum_{(i_1,\dots, i_N) \in \{0,1\}^N, \sum_{n=1}^N i_n \geq 1} \mathbb{E} \left[\|\T{Y}_{i_1\dots i_N}\|_F^2 \mid \M{\Omega_1}, \dots, \M{\Omega}_N\right]\\
&\le \Delta \left(\sum_{(i_1,\dots, i_N) \in \{0,1\}^N, \sum_{n=1}^N i_n \geq 1} \|\T{B}_{i_1\dots i_N}\|_F^2 \right)
= \Delta \|\T{X} - \hat{\T{X}}_2\|_F^2. \nonumber
\end{align*}
\end{proof}

\section{Random matrix projections} \label{s-matrix-projections}

Proofs for the lemmas in this section can be found in \cite[sections 9 and 10]{halko2011finding}.
\begin{lem}
\label{lemma:expectation_inverse_gaussian}
Assume that $t>q$.
Suppose $\mathbf{G}_1\in \mathbb{R}^{t\times q}$ and $\mathbf{G}_2\in \mathbb{R}^{t\times p}$
have i.i.d. standard normal entries.
For any matrix $\mathbf{B}$ with conforming dimensions,
\begin{equation}
\mathbb{E} \|\mathbf{G}_1^\dag \mathbf{G}_2 \mathbf{B}\|_F^2 = \frac{q}{t-q-1} \|\mathbf{B}\|_F^2. \nonumber
\end{equation}
\end{lem}

\begin{lem}
\label{lemma:sketchy_column_space_err}
Given a fixed $\mathbf{A} \in \reals^{m \times n}$
and random $\M{\Omega} \in \reals^{n\times k}$ with
i.i.d. standard normal entries,
let $\M{Q}_\rho = \textup{TruncatedQR}(\M{A\Omega}, \rho) \in \in \reals^{n\times \rho}$
for $\rho < k-1$ \cite{gu1996efficient}.
Then
\begin{equation}
\label{eq:fix_rank_tail_bound}
\mathbb{E}\|(\mathbf{I} - \M{Q_\rho Q_\rho^\top})\mathbf{A}\|_F^2\le
 \frac{\rho}{k-\rho-1} \tau_\rho .
\end{equation}
\end{lem}

\begin{corollary}\label{cor:err-decreasing}
Under the same conditions as in \cref{lemma:sketchy_column_space_err},
suppose $\M{Q} = \textup{QR}(\M{A\Omega})$ is an orthogonal matrix
spanning the column space of $\M{A\Omega}$. Then
\begin{equation}
\label{eq:low_rank_tail_bound}
\mathbb{E}\|(\mathbf{I} - \M{QQ^\top})\mathbf{A}\|_F^2\le
\min_{1\le \rho<k-1} \frac{\rho}{k-\rho-1} \tau_\rho.
\end{equation}
\end{corollary}
\begin{proof}
For each $\rho < k-1$,
\[
\|(\mathbf{I} - \M{Q_\rho Q_\rho^\top})\mathbf{A}\|_F^2
\leq \|(\mathbf{I} - \M{Q_\rho Q_\rho^\top})\mathbf{A}\|_F^2
\leq \frac{r}{k-\rho-1} \tau_\rho,
\]
using \cref{eq:fix_rank_tail_bound} for the second inequality.
Minimize over $\rho < k-1$ to reach the result.
\end{proof}

\section{Time and Storage Complexity}\label{appendix: time-complexity}
\subsection{Comparison Between \cref{alg:fix_rank_appro} and T.-TS \cite{malik2018low}}
Here we compare the time and storage complexity of the two extant methods for
streaming Tucker approximation: our one-pass method, and T.-TS \cite{malik2018low}.

To compare the storage and time costs of both T.-TS and the one-pass algorithm,
we separate the cost into two parts:
one for forming the sketch,
the other for each iteration of ALS.
Assume the tensor to approximate has equal side lengths $I_1=\cdots = I_N = I$
and that the target rank for each mode is $R$.

The suggested default parameters for the sketch in \cite{malik2018low}
are $J_1 = 10R^{N-1}$ and $J_2 = 10R^{N}$.
Our suggested default parameters are $k=2r, s=2k+1$. Under the choice of the default parameter, we compare the the cost of storage and time in \cref{tab:storage-comparison} and \cref{tab:time-comparison}.
In most problems with data that is not exactly low rank, i.e. $R > 4$, the suggested default setting of T.-TS typically leads to a higher storage cost. Moreover, our algorithm uses less storage and is faster to compute, particularly for tensors with many modes $N$.

However, the evaluation of the two algorithms should not be solely based on their default setups. If the memory constraint is set to be the same, our one-pass algorithm performs much better in the low-memory case,
but slightly worse in the high-memory case (see \cref{fig:vary-memory}).
The memory required by our default parameters is typically much smaller than
that required with the default parameters of \cite{malik2018low}.

\subsection{Computational Complexity of \cref{alg:fix_rank_appro}}
Here, we will calculate the computational complexity for our one-pass fixed-rank approximation algorithm.

In the sketching stage of the streaming algorithm,
we first need to compute the factor sketches, $\mathbf{G}_n = \mathbf{X}\mathbf{\Omega}_n, n \in [N]$ with $kN\hat{I}$ flops in total.
Then we need to compute the core tensor sketch $\mathscr{Z}$ by recursively multiplying $\mathscr{X}$ by $\mathbf{\Phi}_n, n \in [N]$.
We can upper bound the number of flops by $\frac{s(1-\delta_1^N)}{1-\delta_1}\bar{I}$.
Then in the approximation stage, we first perform ``economy size'' QR factorizations
of $\mathbf{G}_1, \dots, \mathbf{G}_N$ with $\mathscr{O}(k^2(\sum_{n =1}^N I_n))$
to find the orthonormal bases $\mathbf{Q}_1, \dots, \mathbf{Q}_N$.
To find the linkage tensor $\mathscr{W}$, we need to recursively solve linear square problems
with $\frac{k^2s^N(1-(k/s)^N)}{1-k/s}$ flops.
Overall, the sketch computation dominates the total time complexity.

The HOSVD directly acts on $\mathscr{X}$ by first computing the SVD for each unfolding ($\mathscr{O}(kN\bar{I})$)
and then multiplying $\mathscr{X}$ by $\mathbf{U}_1^\top, \dots, \mathbf{U}_N^\top$ ($\mathcal{O}(\frac{k(1-\delta_1^N)\bar{I}}{1-\delta_1})$).
The total time cost is less than the streaming algorithm with a constant factor.
Note: we can use the randomized SVD in the first step of the HOSVD
to improve the computational cost to $\bar{I}N\log k + \sum_{n = 1}^N(I_{n}+I_{(-n)})k^2$ \cite{halko2011finding,}.

\begin{table}[h!]
	\centering
	\begin{tabular}{c c c c }
		Algorithm  & & Storage Cost ($I=o(r^{2N})$) \\
		\hline

		\multirow{2}{*}{T.-TS} & Sketching & $\mathcal{O}(r^{2N})$ \\
		& Recovery &  $\mathcal{O}(r^{2N})$ & \\
		\hline\hline
		\multirow{2}{*}{\cref{alg:one_pass_low_rank_appro} (One Pass)} & Sketching &  $\mathcal{O}(4^Nr^N)$  \\
		& Recovery  & $\mathcal{O}(4^Nr^N)$  \\
		\hline
	\end{tabular}
	\caption{Storage complexity of \cref{alg:one_pass_low_rank_appro} and T.-TS
	on tensor $\T{X} \in \mathbb{R}^{I \times \dots \times I}$.
	\cref{alg:one_pass_low_rank_appro} uses parameters $(k,s) = (2r, 4r+1)$
	and uses a TRP composed of Gaussian DRMs inside the Tucker sketch.
	T.-TS uses default values for hyper-parameters: $J_1=10r^{N-1}, J_2=10r^{N}$.}
	\label{tab:storage-comparison}
\end{table}

\begin{table}[h!]
	\centering
	\begin{tabular}{c c c c }
		Algorithm  & & Time Cost ($I = o(r^{2N})$) \\
		\hline

		\multirow{2}{*}{T.-TS} & Sketching & $\mathcal{O}(N\rm{nnz}(\T{X}))$ \\
		& Recovery &  $\mathcal{O}(NIr^N+Nr^{2N-1}+r^{2N})$ & \\
		\hline\hline
		\multirow{2}{*}{\cref{alg:one_pass_low_rank_appro} (One Pass)} & Sketching &  $\mathcal{O}(Nr~ \rm{nnz}(\T{X})))$  \\
		& Recovery  & $\mathcal{O}(Nr^{N+1})$  \\
		\hline
	\end{tabular}
	\caption{Time complexity of \cref{alg:one_pass_low_rank_appro} and T.-TS
	on tensor $\T{X} \in \mathbb{R}^{I \times \dots \times I}$.
	\cref{alg:one_pass_low_rank_appro} uses parameters $(k,s) = (2r, 4r+1)$
	and uses a TRP composed of Gaussian DRMs inside the Tucker sketch.
	T.-TS uses default values for hyper-parameters: $J_1=10r^{N-1}, J_2=10r^{N}$.
	}
	\label{tab:time-comparison}
\end{table}

\mnote{Instead of "Proposed", we should be clear about which of our algorithms we're referring to.}

\section{More Numerics}\label{appendix:more_result}

\begin{figure}
	\centering
	\begin{subfigure}{0.3\textwidth}
		\includegraphics[scale = 0.24]{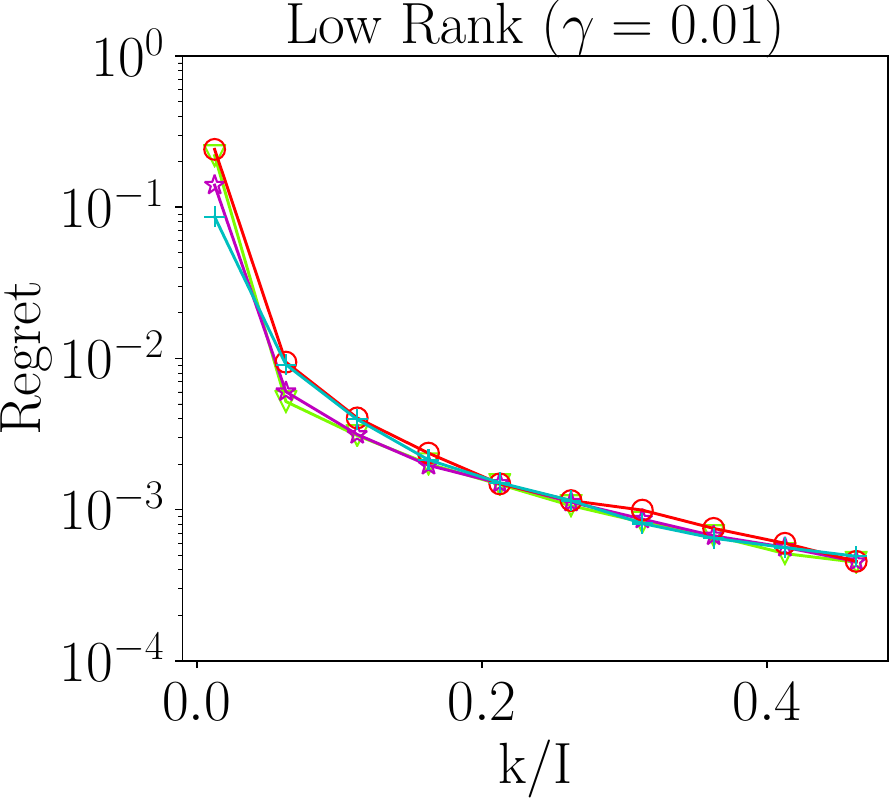}
	\end{subfigure}
	\begin{subfigure}{0.3\textwidth}
		\includegraphics[scale = 0.24]{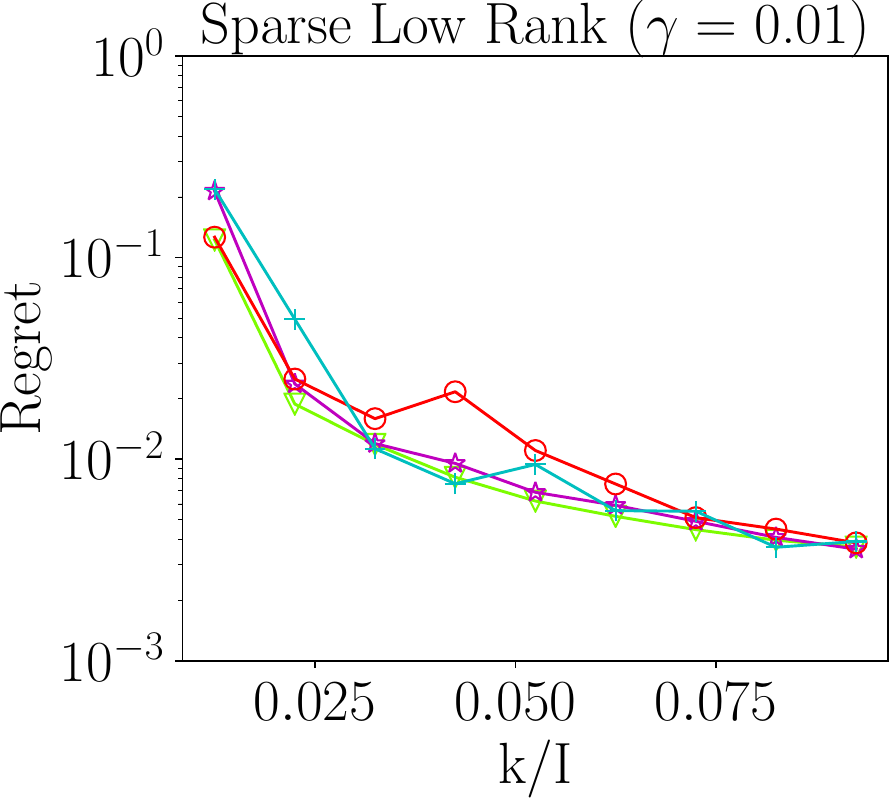}
	\end{subfigure}
	\begin{subfigure}{0.3\textwidth}
		\includegraphics[scale = 0.24]{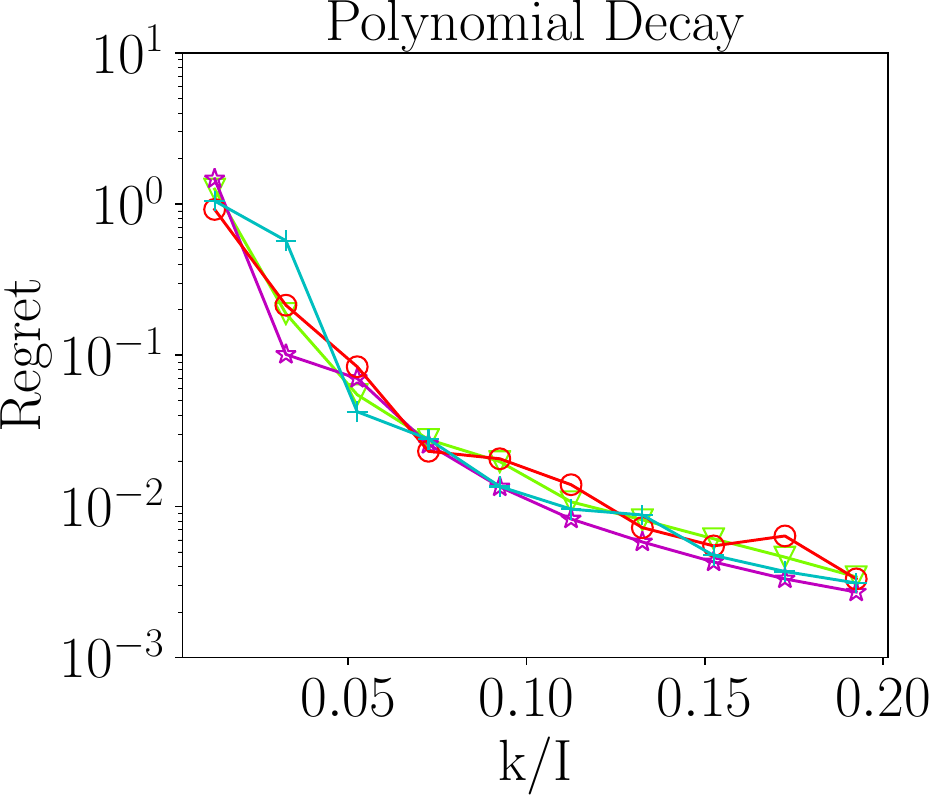}
	\end{subfigure}\\
	\begin{subfigure}{0.3\textwidth}
		\includegraphics[scale = 0.24]{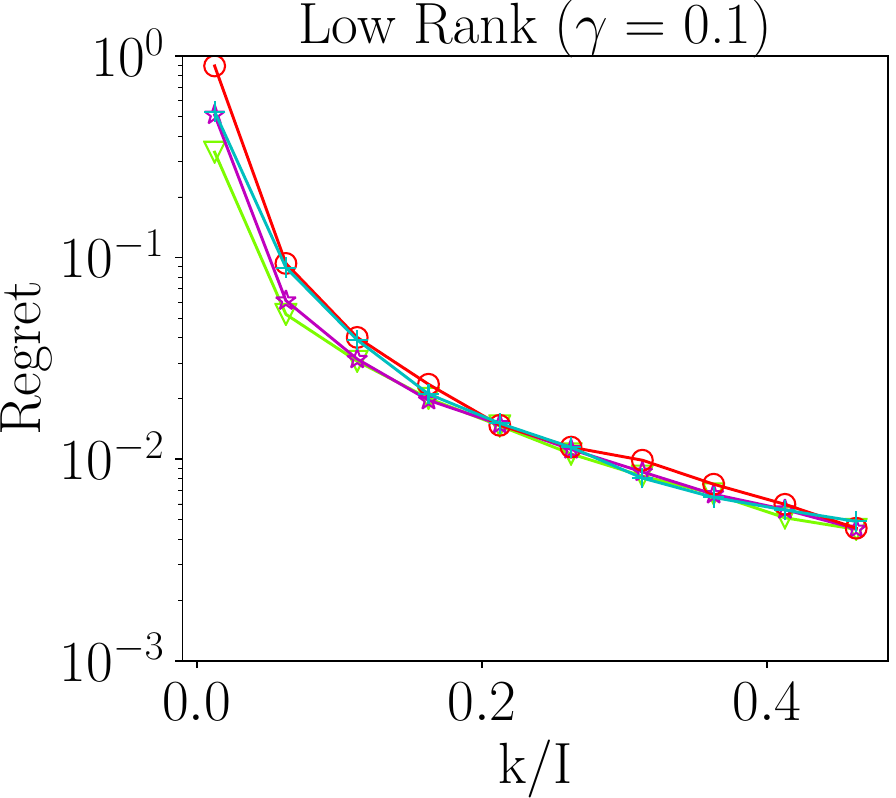}
	\end{subfigure}
	\begin{subfigure}{0.5\textwidth}
		\includegraphics[scale = 0.24]{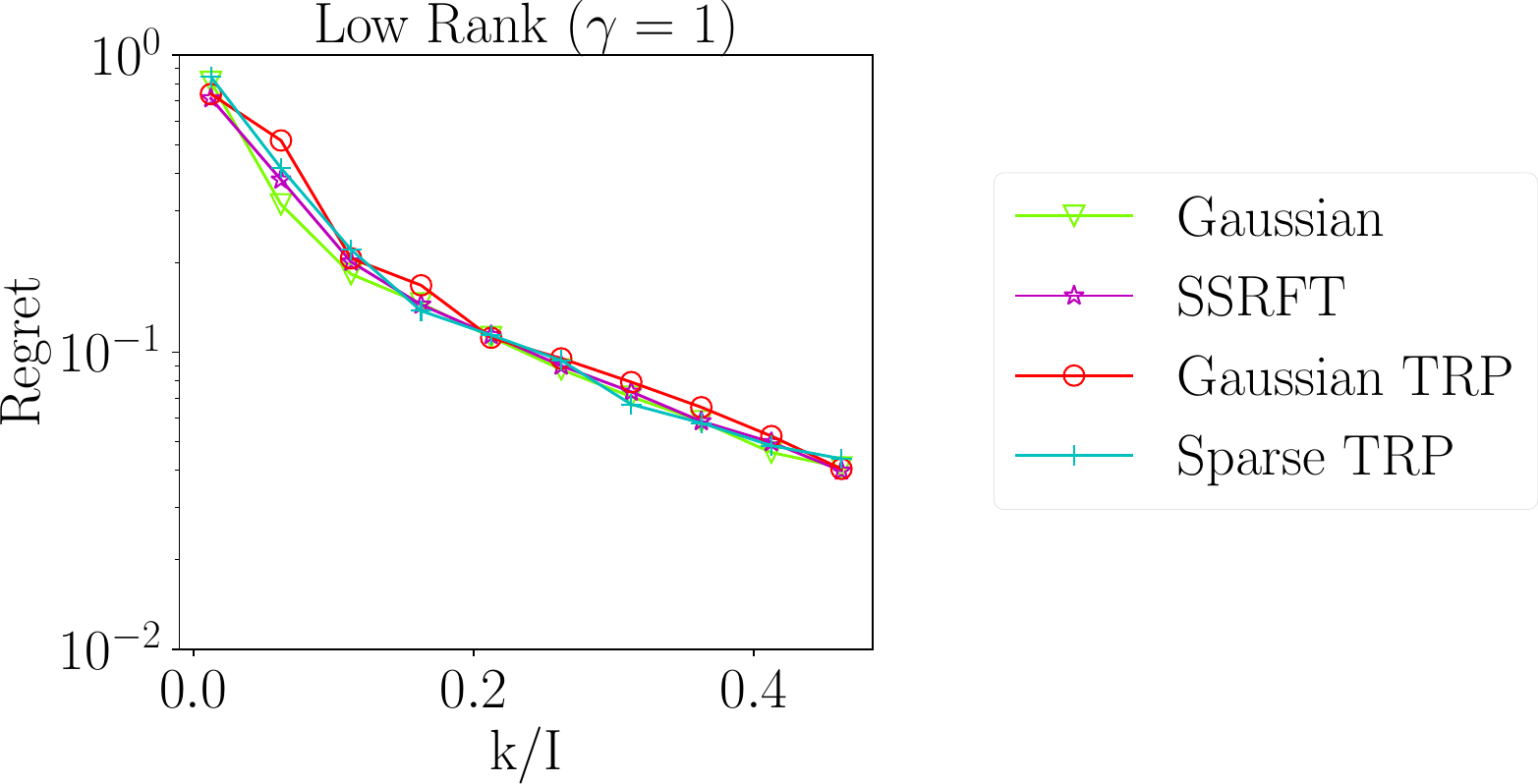}
	\end{subfigure}
	\caption{We approximate 3D synthetic tensors (see \cref{s-synthetic-data}) with $I = 400$,
		using our one-pass algorithm with $r = 5$ and varying $k$ ($s = 2k+1$),
		using a variety of DRMs in the Tucker sketch:
		Gaussian, SSRFT, Gaussian TRP, or Sparse TRP.}
		\label{fig:vary-k-400-app}
\end{figure}

\begin{figure}
	\centering
	\begin{subfigure}{0.3\textwidth}
		\includegraphics[scale = 0.25]{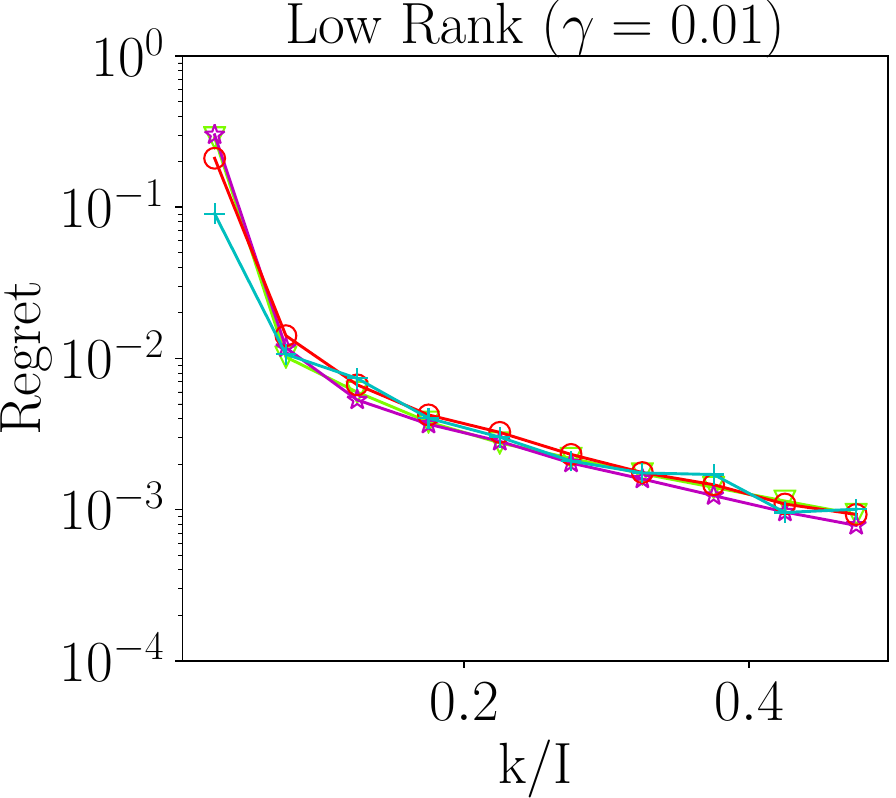}
	\end{subfigure}
	\begin{subfigure}{0.3\textwidth}
		\includegraphics[scale = 0.25]{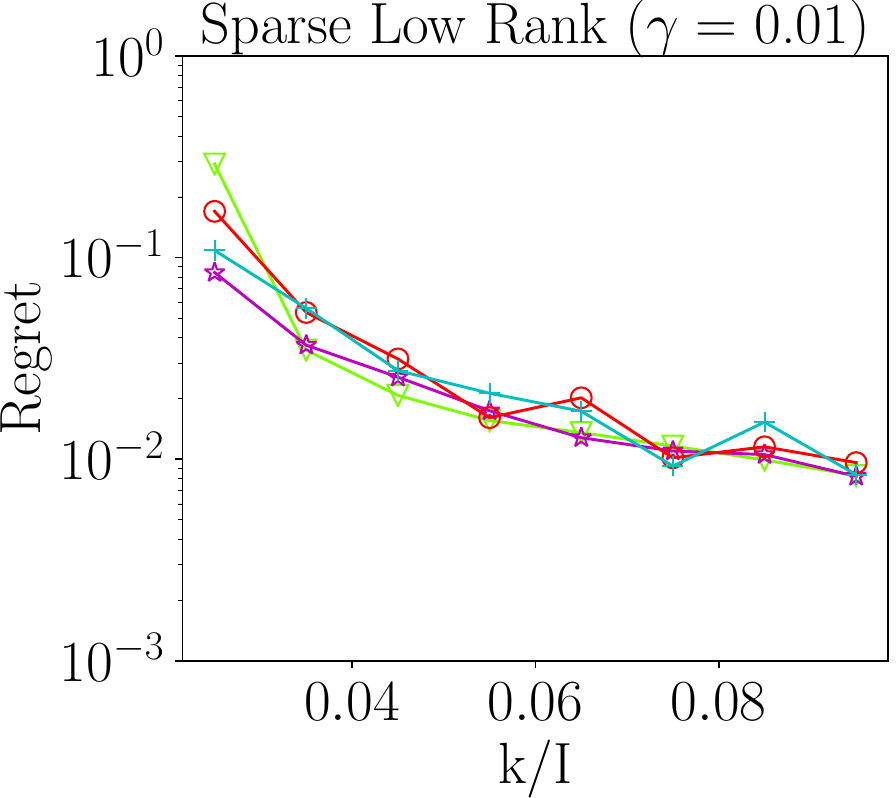}
	\end{subfigure}
	\begin{subfigure}{0.3\textwidth}
		\includegraphics[scale = 0.25]{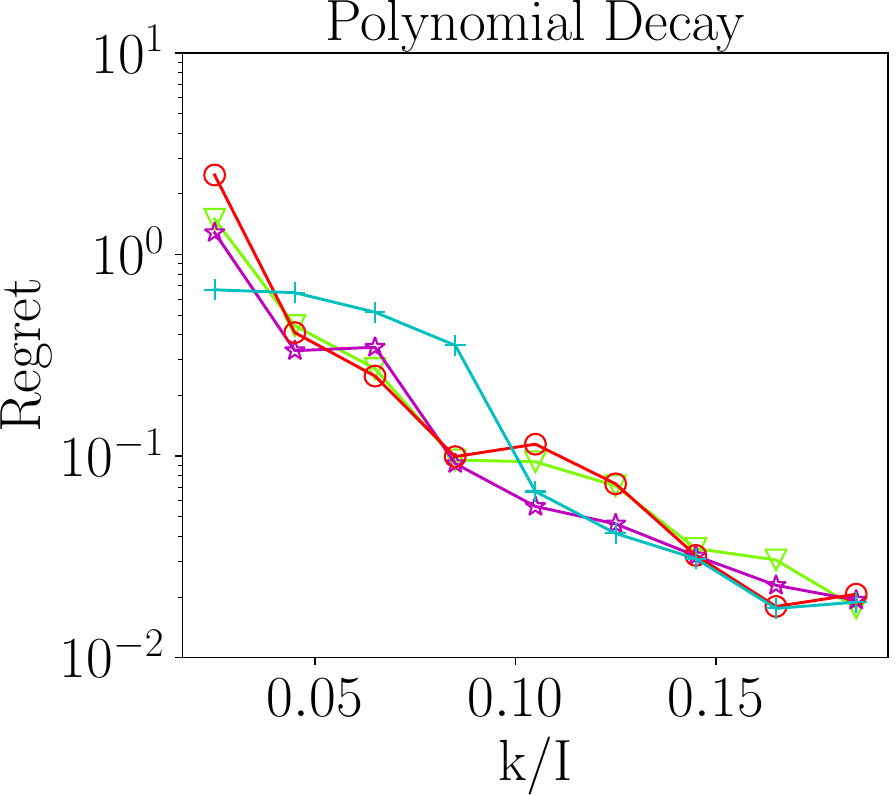}
	\end{subfigure}\\
	\begin{subfigure}{0.3\textwidth}
	\includegraphics[scale = 0.25]{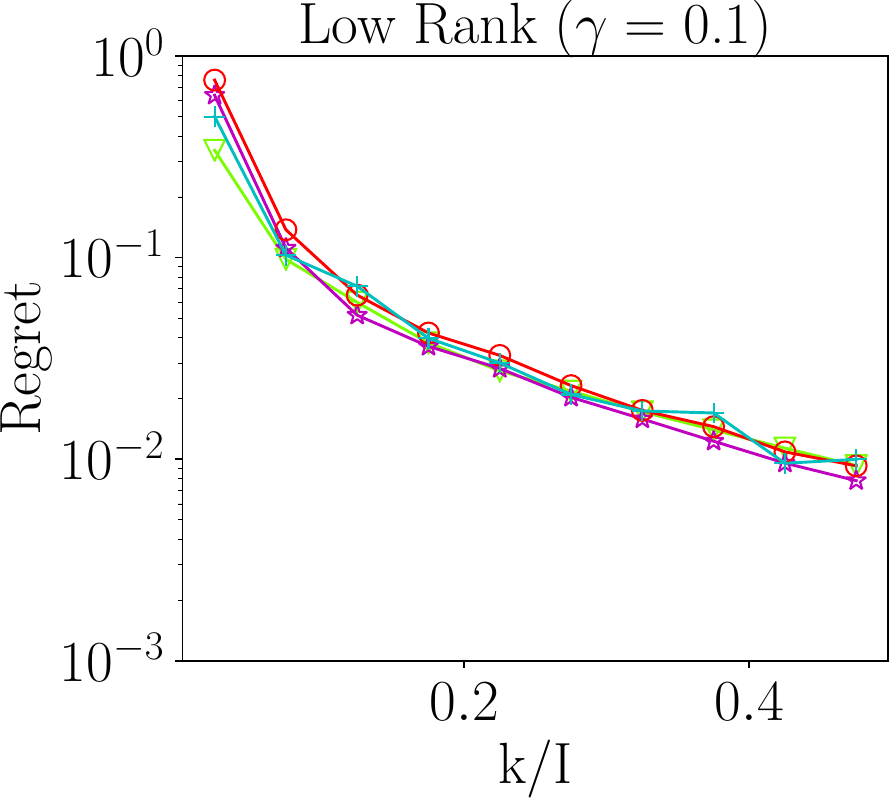}
	\end{subfigure}
	\begin{subfigure}{0.55\textwidth}
		\includegraphics[scale = 0.25]{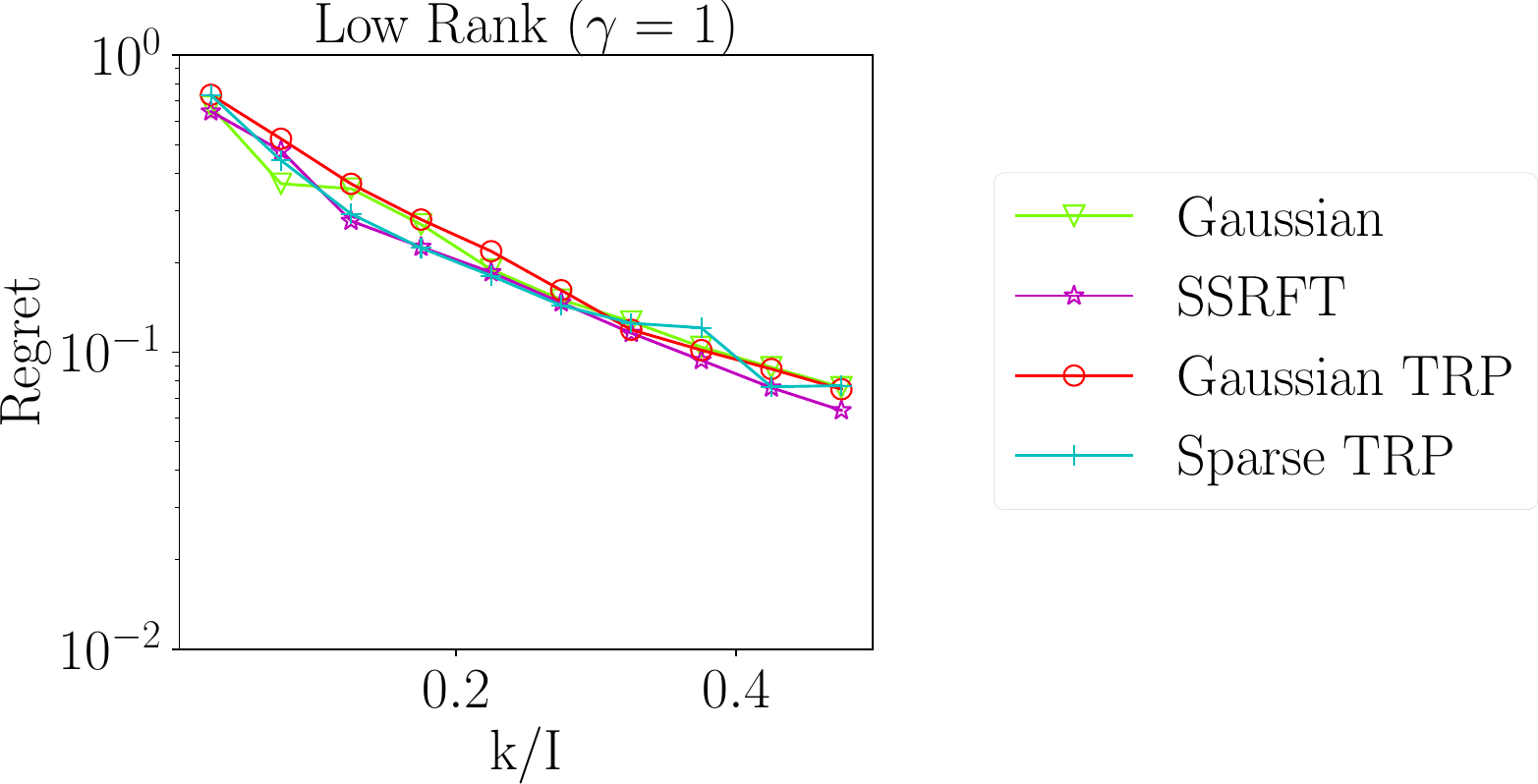}
	\end{subfigure}
	\caption{We approximate 3D synthetic tensors (see \cref{s-synthetic-data}) with $I = 200$,
		using our one-pass algorithm with $r = 5$ and varying $k$ ($s = 2k+1$),
		using a variety of DRMs in the Tucker sketch:
		Gaussian, SSRFT, Gaussian TRP, or Sparse TRP.}
		\label{fig:vary-k-200-app}
\end{figure}

\begin{figure}
	\centering
	\begin{subfigure}{0.3\textwidth}
		\includegraphics[scale = 0.24]{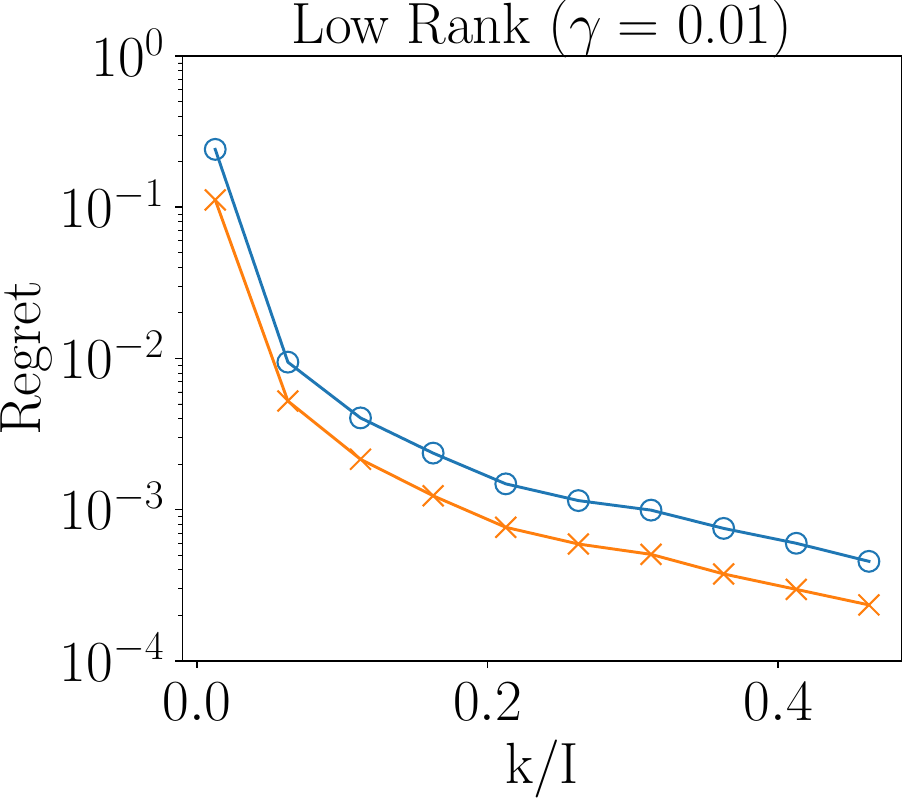}
	\end{subfigure}
	\begin{subfigure}{0.3\textwidth}
		\includegraphics[scale = 0.24]{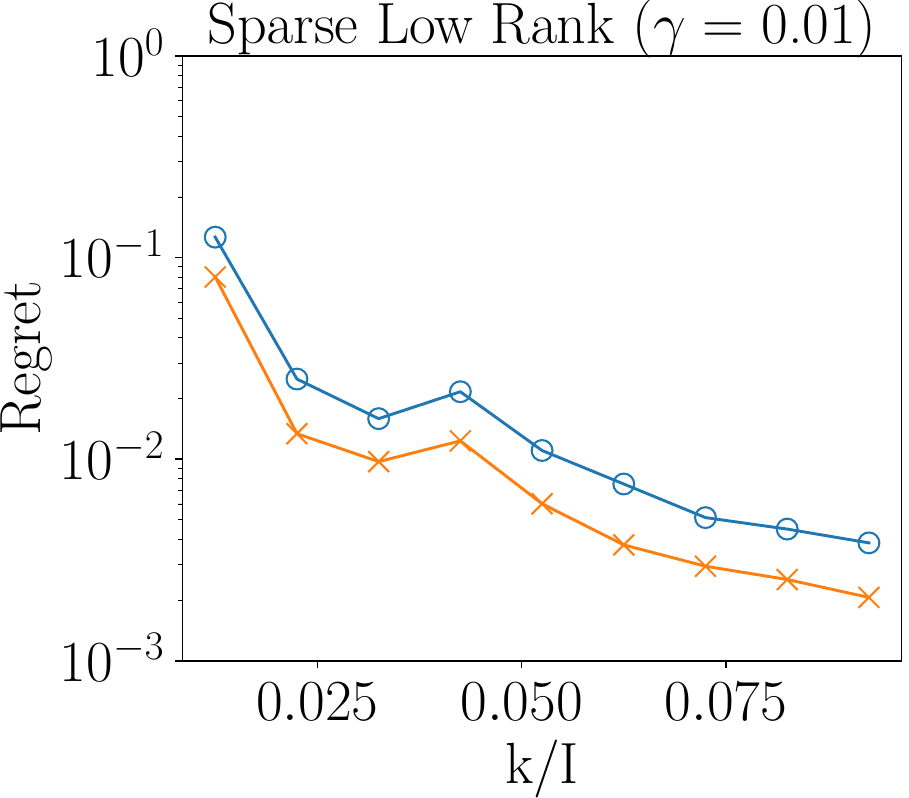}
	\end{subfigure}
	\begin{subfigure}{0.3\textwidth}
		\includegraphics[scale = 0.24]{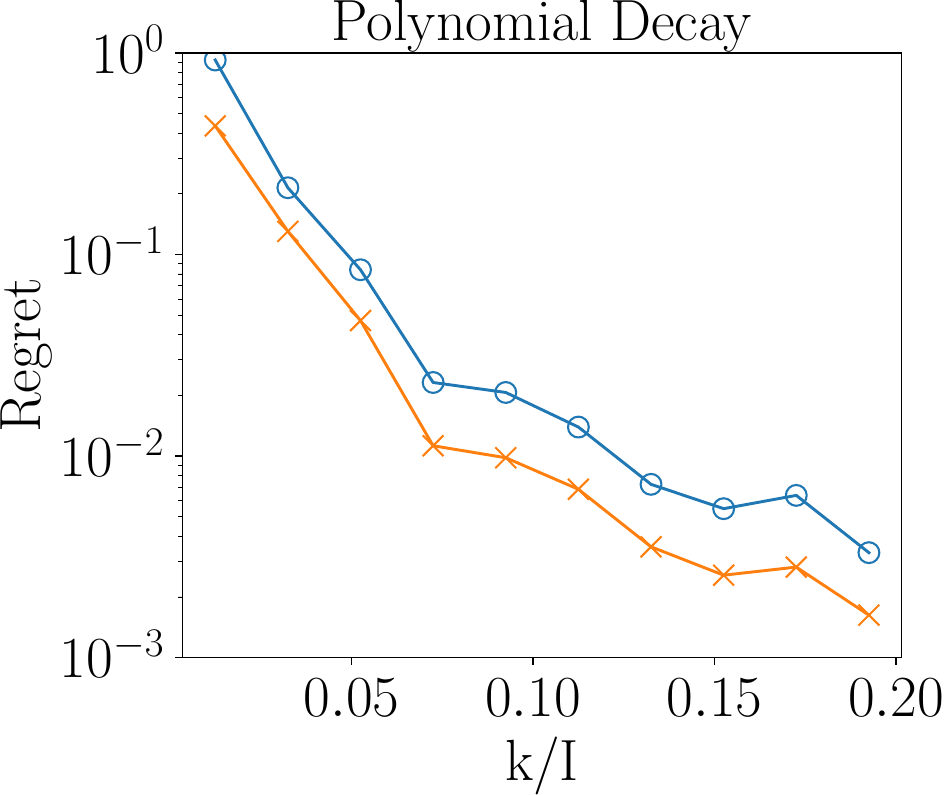}
	\end{subfigure}\\
	\begin{subfigure}{0.3\textwidth}
		\includegraphics[scale = 0.24]{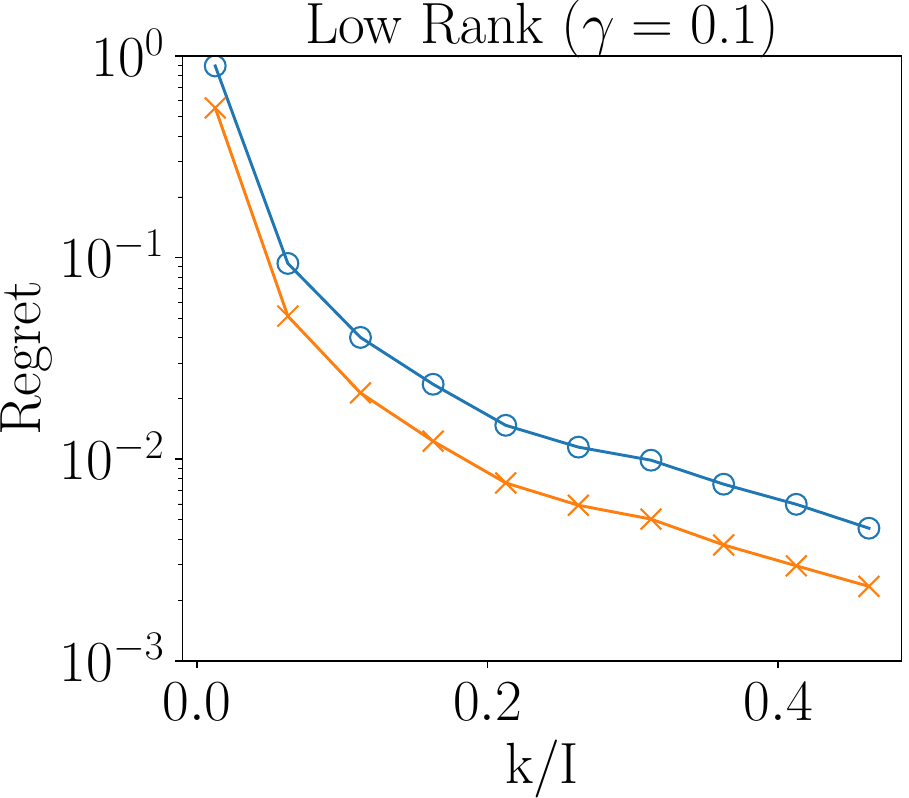}
	\end{subfigure}
	\begin{subfigure}{0.55\textwidth}
		\includegraphics[scale = 0.24]{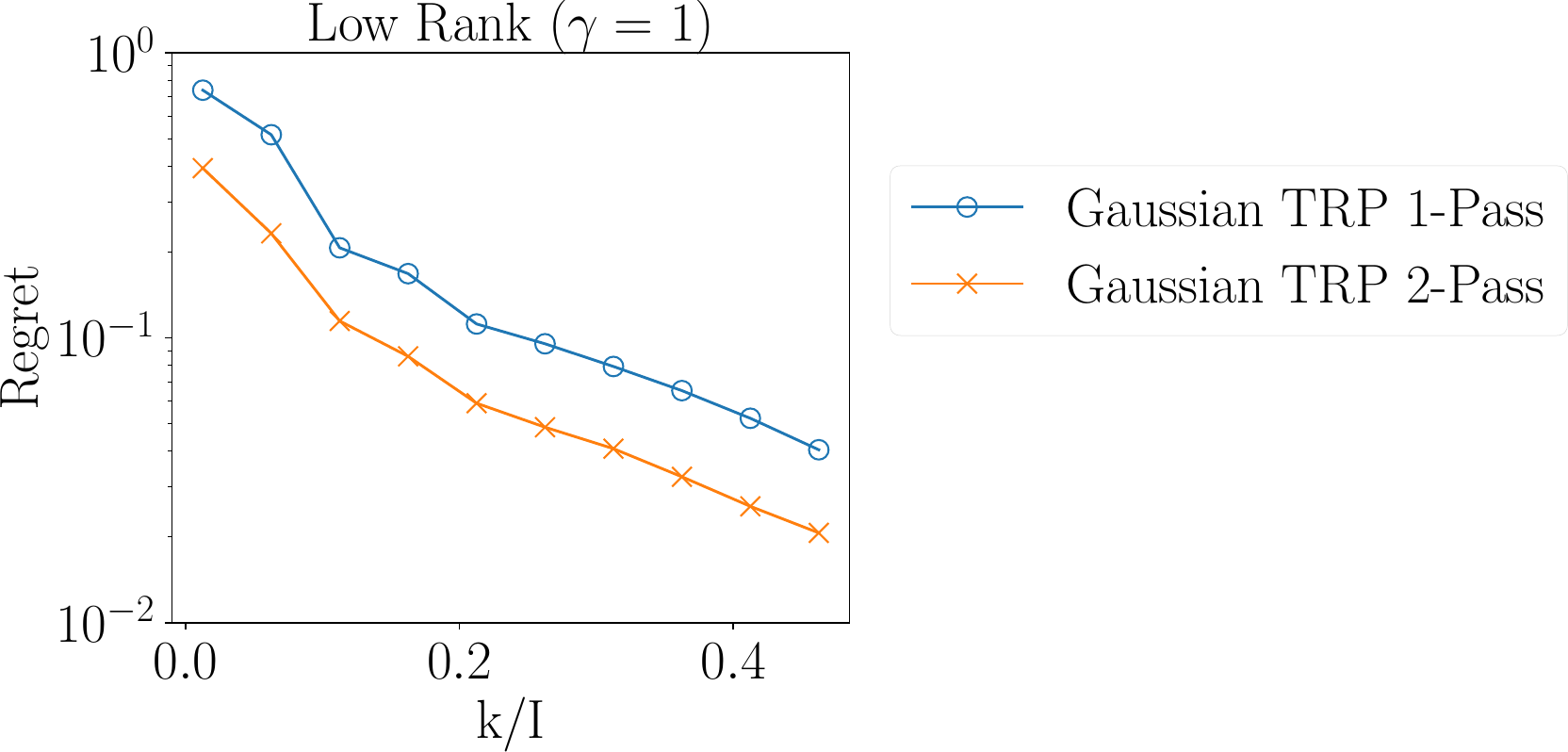}
	\end{subfigure}
	\caption{We approximate 3D synthetic tensors (see \cref{s-synthetic-data}) with $I = 400$,
		using our one-pass and two-pass algorithms with $r = 5$ and varying $k$ ($s = 2k+1$),
		using the Gaussian TRP in the Tucker sketch.}\label{fig:vary-k-400-compare-app}
\end{figure}

\begin{figure}
	\centering
	\begin{subfigure}{0.3\textwidth}
		\includegraphics[scale = 0.25]{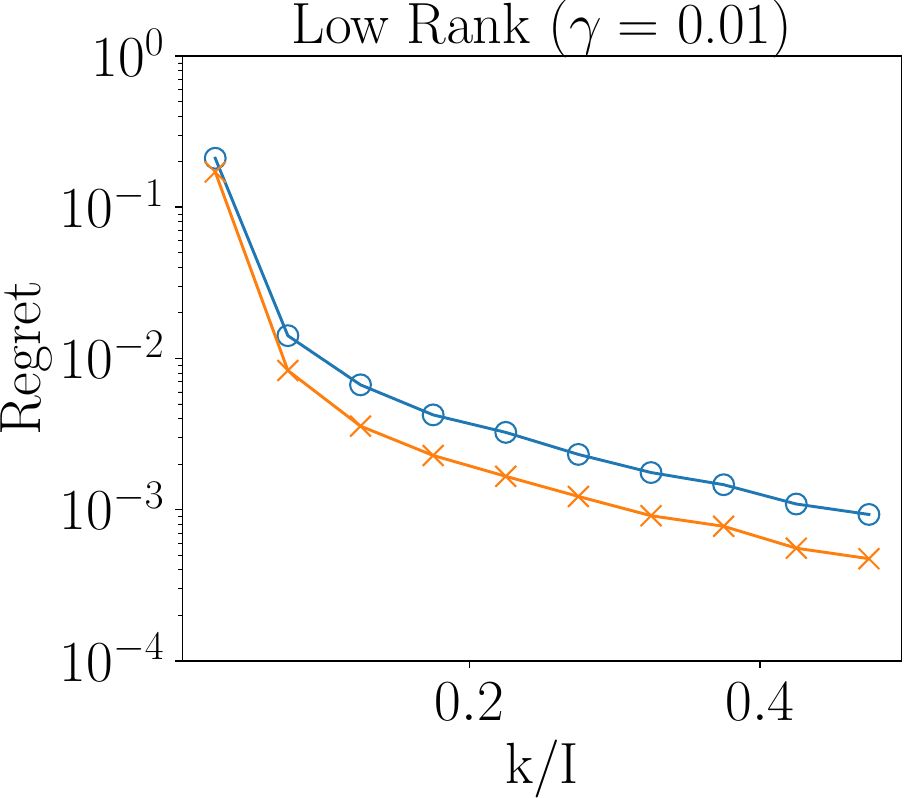}
	\end{subfigure}
	\begin{subfigure}{0.3\textwidth}
		\includegraphics[scale = 0.25]{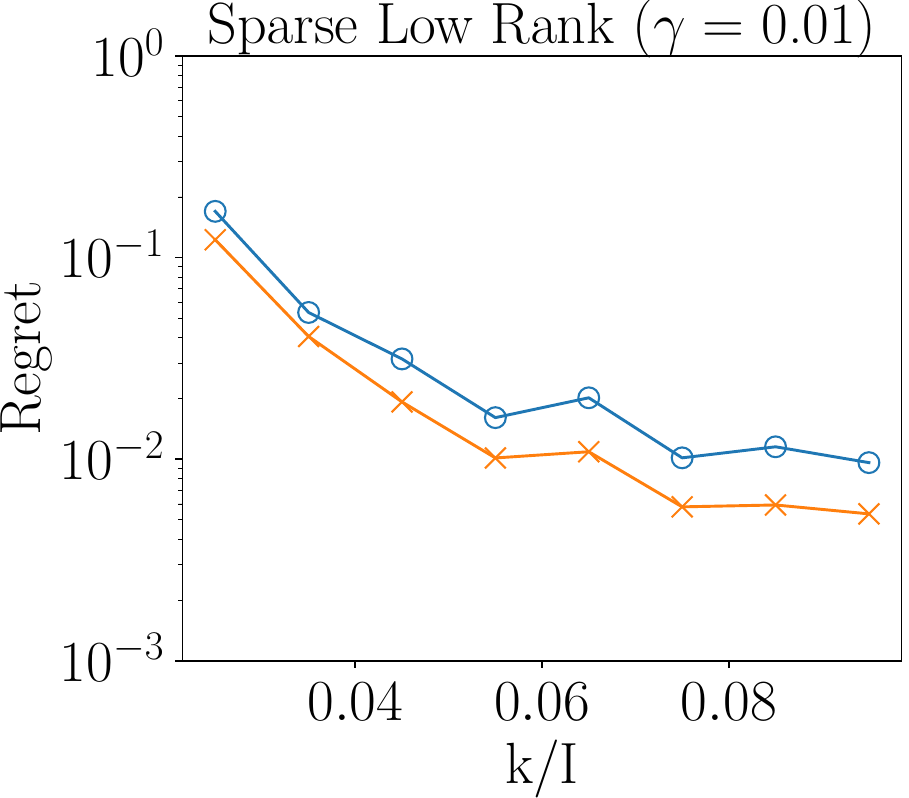}
	\end{subfigure}
	\begin{subfigure}{0.3\textwidth}
		\includegraphics[scale = 0.25]{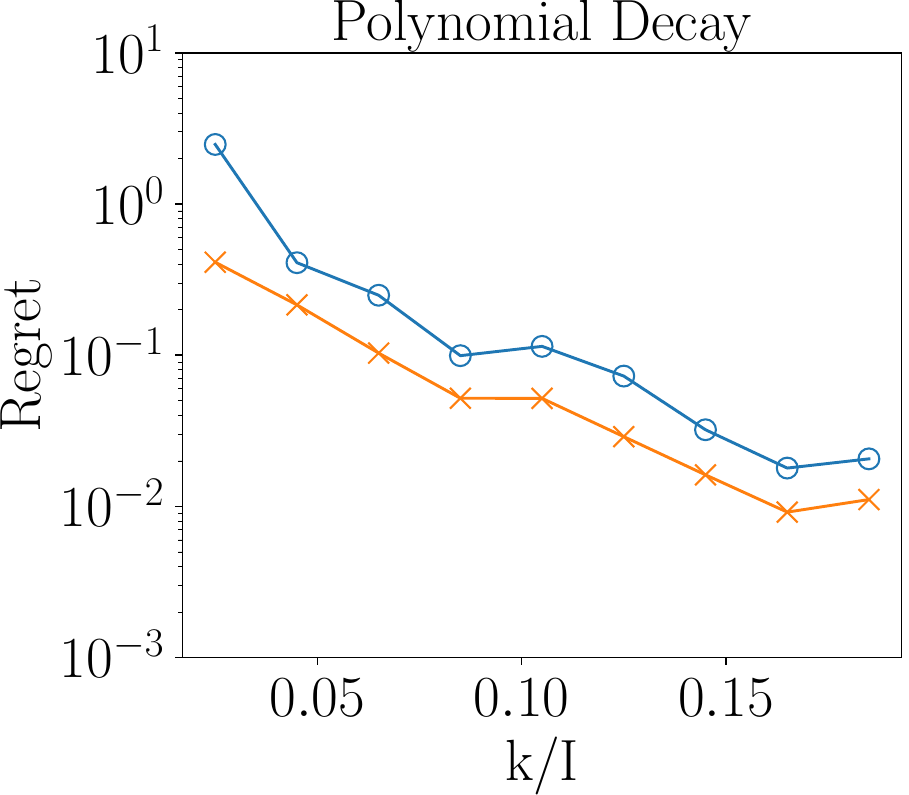}
	\end{subfigure}\\
	\begin{subfigure}{0.3\textwidth}
		\includegraphics[scale = 0.25]{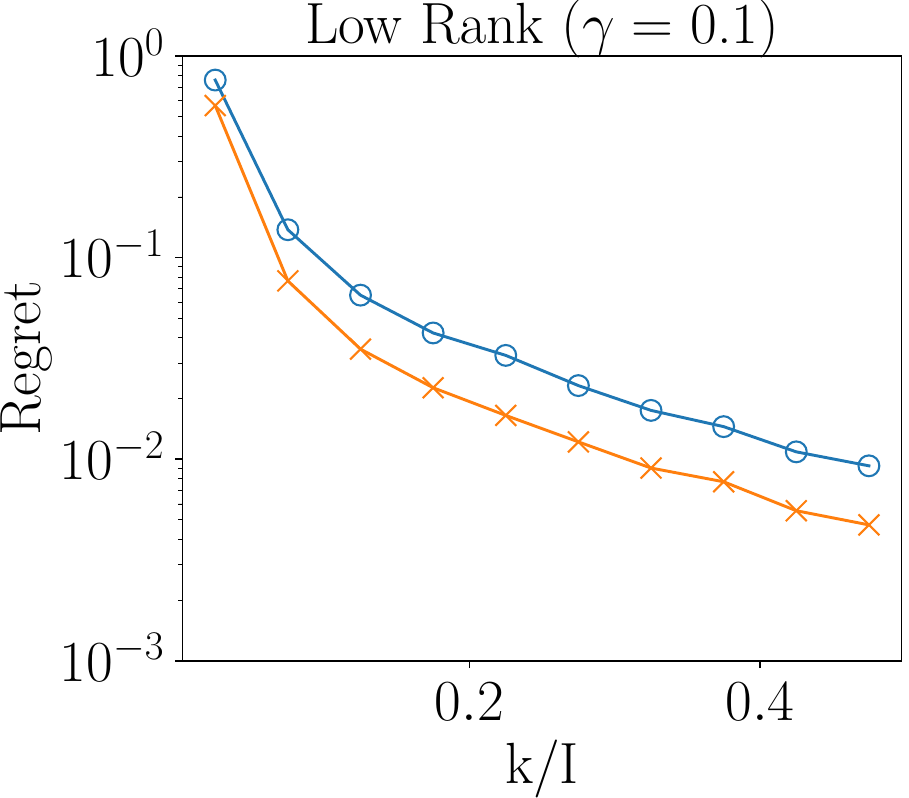}
	\end{subfigure}
	\begin{subfigure}{0.55\textwidth}
		\includegraphics[scale = 0.25]{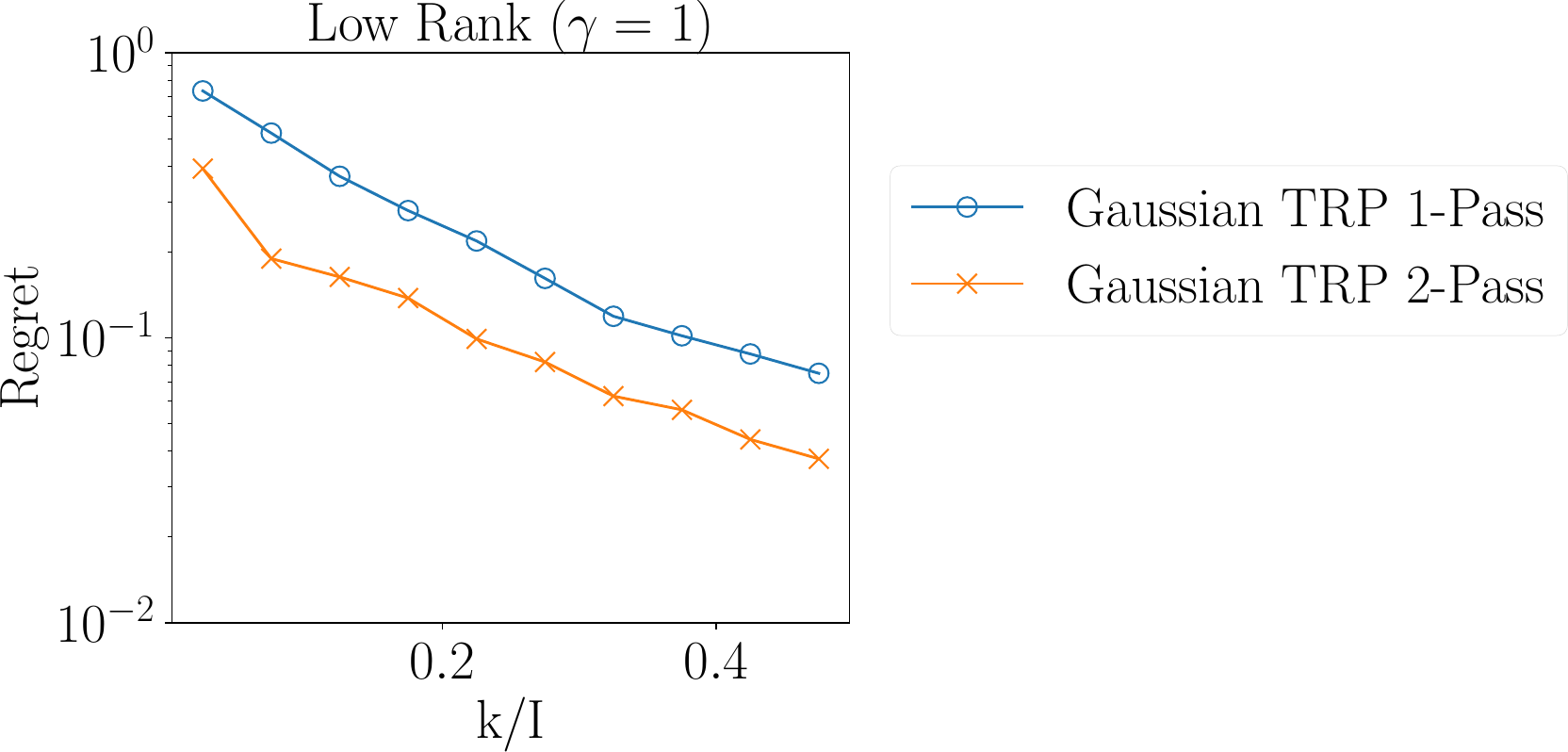}
	\end{subfigure}
	\caption{We approximate 3D synthetic tensors (see \cref{s-synthetic-data}) with $I = 200$,
		using our one-pass and two-pass algorithms with $r = 5$ and varying $k$ ($s = 2k+1$),
		using the Gaussian TRP in the Tucker sketch.}\label{fig:vary-k-200-compare-app}

\end{figure}

This section provides more numerical results on simulated datasets
in \cref{fig:vary-k-400-app},
\cref{fig:vary-k-200-app},
\cref{fig:vary-k-400-compare-app},
and \cref{fig:vary-k-200-compare-app}.


\label{appendix:more_real_data_result}

\begin{figure}[ht]
	\centering
	\includegraphics[height=2.9cm]{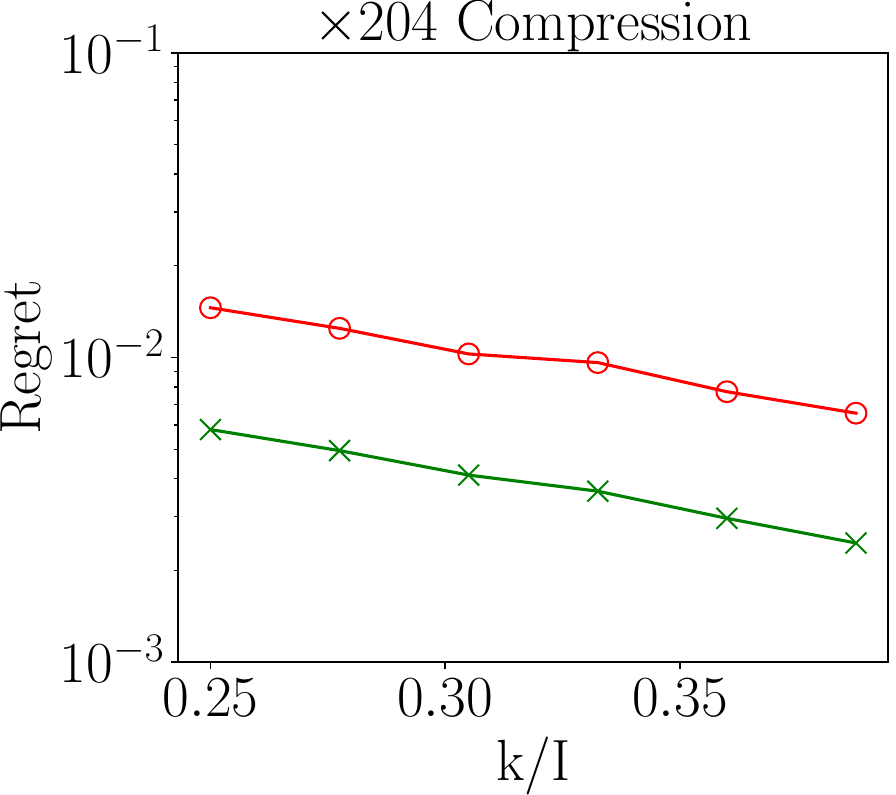}
	\includegraphics[height=2.9cm]{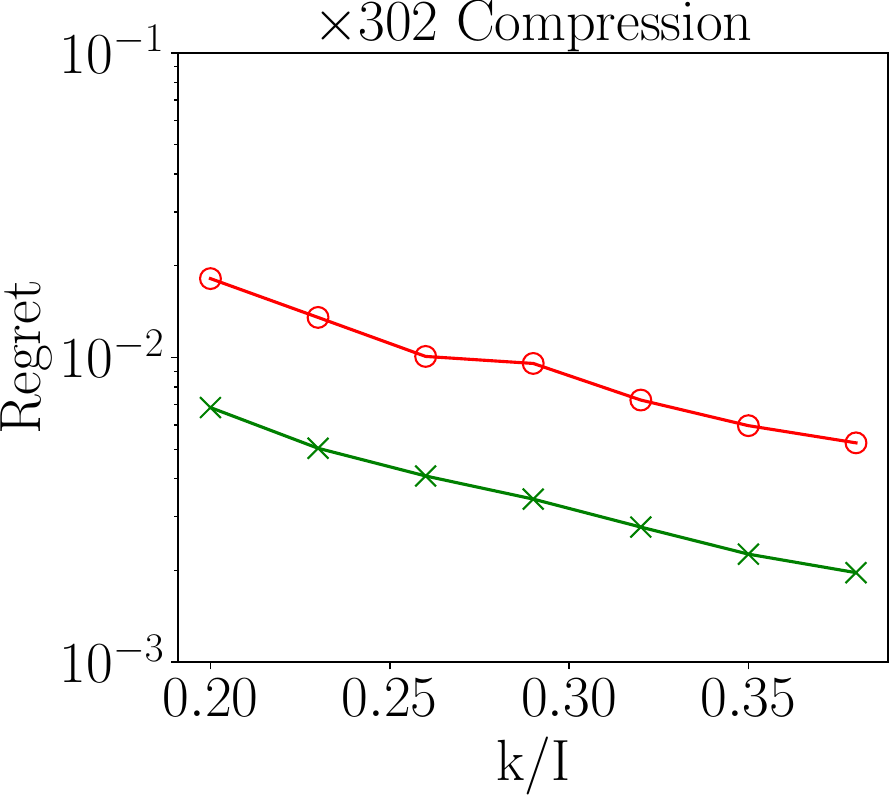}
	\includegraphics[height=2.9cm]{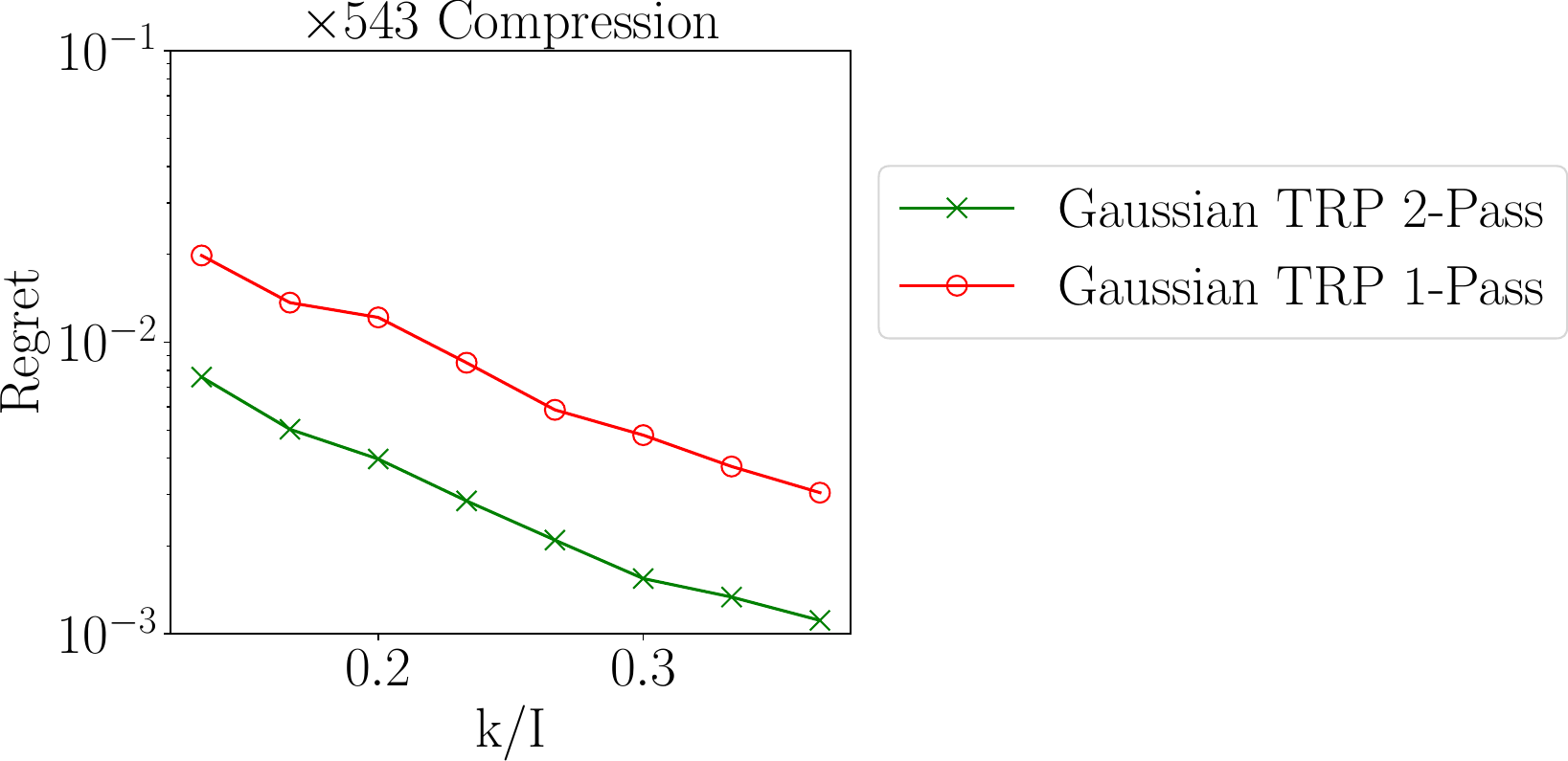}\\
	\textbf{Net Radiative Flux at Surface}\\~\\
	\centering
	\includegraphics[height=2.9cm]{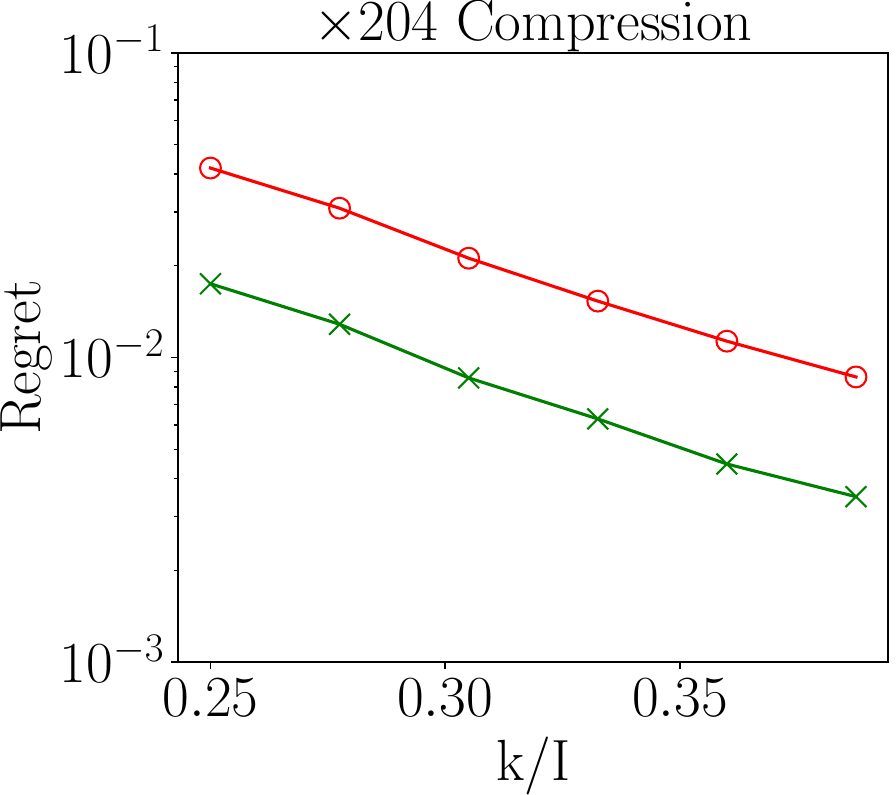}
	\includegraphics[height=2.9cm]{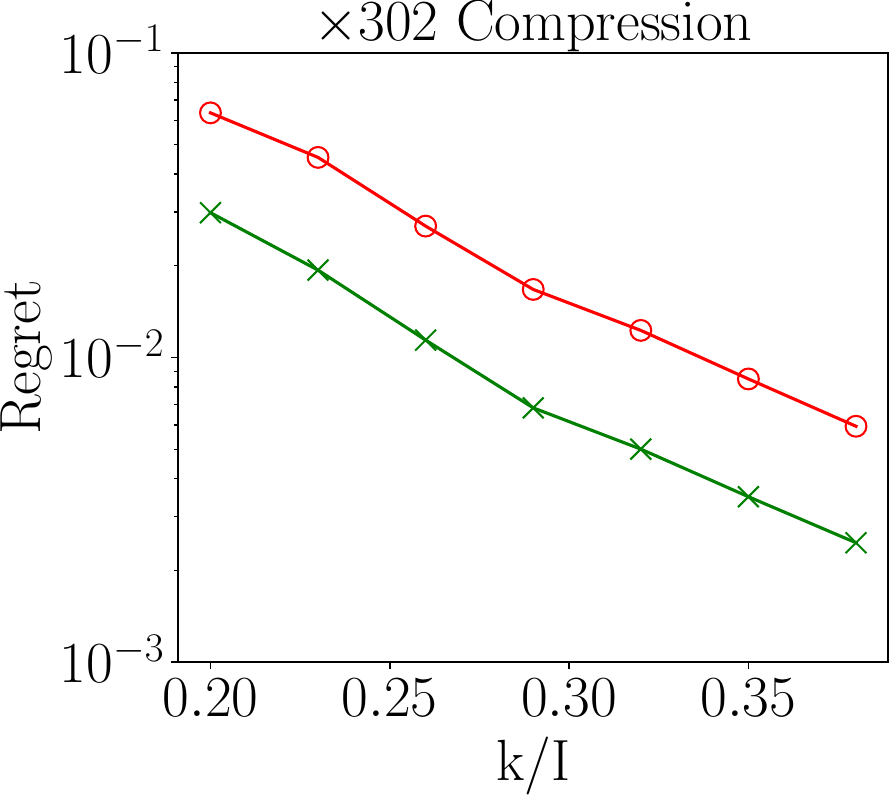}
	\includegraphics[height=2.9cm]{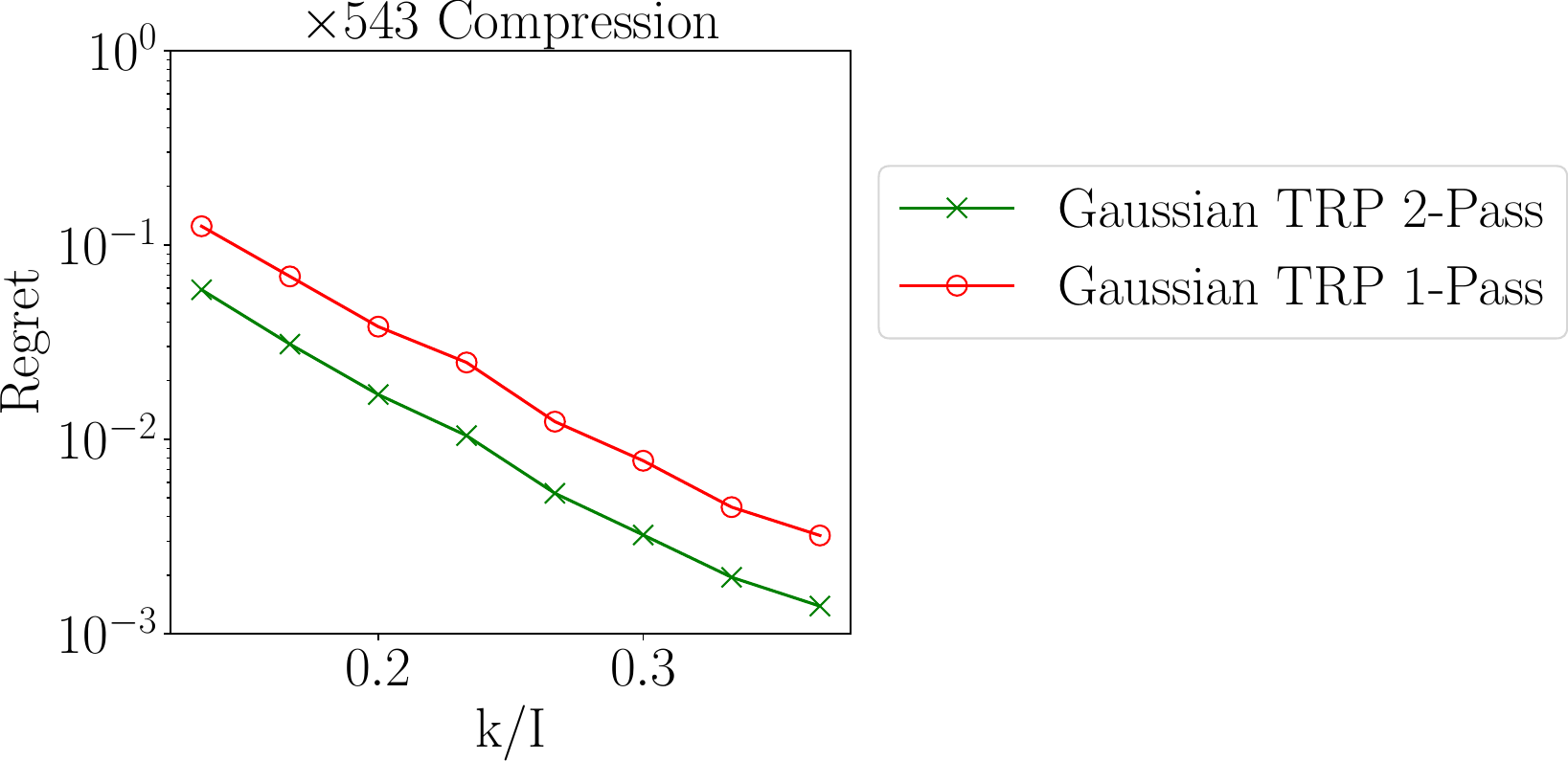}\\
	\textbf{Dust Aerosol Burden}
	\caption{We approximate the net radiative flux and dust aerosol burden data using our one-pass and two-pass algorithms using Gaussian TRP. We compare the performance under different ranks ($r/I = 0.125,0.2,0.067$). The dataset comes from the CESM CAM. The dust aerosol burden measures the amount of aerosol contributed by the dust. The net radiative flux determines the energy received by the earth surface through radiation. } \label{fig:srfrad_burden_dust}
\end{figure}

We also provide more numerical results on real datasets in \cref{fig:srfrad_burden_dust}.

\section{More Algorithms}
This section provides detailed implementations. \par

\begin{algorithm}[ht]
	\caption{Higher order orthogonal iteration (HOOI)
		\cite{de2000multilinear}
		\label{alg:hooi}}
	\textbf{Given:} tensor $\T{X}$, target rank $\V{r} = (r_1, \ldots, r_N)$ \\
	Initialize: compute $\T{X} \approx \llbracket\T{G}; \M{U}_1, \ldots, \M{U}_N \rrbracket$ using HOSVD \\
	Repeat:
	\begin{enumerate}
		\item \emph{Factors.} For each $n \in [N]$,
		\begin{equation}
		\label{eq:factor-update}
		\mathbf{U}_n \leftarrow \argmin_{ \M{U_n} }
		\left\| \llbracket\T{G}; \M{U}_1, \ldots, \M{U}_N \rrbracket - \T{X} \right\|_F^2,
		\end{equation}
		\item \emph{Core.}
		\begin{equation}
		\begin{aligned}
		\label{eq:core_update}
		& \T{G} \leftarrow \argmin_{\T{G}} 
		\left\| \llbracket\T{G}; \M{U}_1, \ldots, \M{U}_N \rrbracket - \T{X} \right\|_F^2.\\
		&  \ie  ~~ \T{G} = \T{X}\times_1 \M{U}_1^\top \times_2 \cdots \times_N \M{U}_N^\top 
		\end{aligned}
		\end{equation}
	\end{enumerate}
	\textbf{Return:} Tucker approximation $\T{X}_{\rm{HOOI}} = \llbracket\T{G}; \M{U}_1, \ldots, \M{U}_N \rrbracket$
\end{algorithm}
Notice the core update \eqref{eq:core_update}
admits the closed form solution
$\T{G} \leftarrow \T{X}\times_1 \mathbf{U}_1^\top \cdots \times_N \mathbf{U}_N^\top$,
which motivates the second step of HOSVD
for a linear sketch appropriate to a streaming setting (\cref{alg:linear_update})
or a distributed setting (\cref{alg:sketch_distributed}).
\label{appendix:more_algorithms}
\begin{algorithm}[th]
	\caption{Linear Update to Sketches}\label{alg:linear_update}
	\begin{algorithmic}[1]
		\Function {SketchLinearUpdate}{$\T{F}, \mathbf{V}_1, \dots, \mathbf{V}_N, \T{H}$; $\theta_1$, $\theta_2$}
		\For{$n = 1, \dots, N$}
		\State $\mathbf{V}_n \leftarrow \theta_1 \mathbf{V}_n + \theta_2 \mathbf{F}^{(n)} \mathbf{\Omega}_n $
		\EndFor
		\State $\T{H} \leftarrow \theta_1 \T{H} + \theta_2 \T{F} \times_1 \mathbf{\Phi}_1 \times \cdots \times_N \mathbf{\Phi}_N $
		\State \Return $(\mathbf{V}_1, \dots, \mathbf{V}_N, \T{H})$
		\EndFunction
	\end{algorithmic}
\end{algorithm}

\begin{algorithm}[th]
\begin{algorithmic}[1]
\caption{Sketching in Distributed Setting}\label{alg:sketch_distributed}
\Require{$\T{X}_i$ is the part of the tensor $\T{X}$ at local machine $i$ and $\T{X} = \sum_{i=1}^m\T{X}_i$.
}
\Function{ComputeSketchDistributed}{$\T{X}_1, \ldots, \T{X}_m$}
\State Send the same random generating environment to every local machine.
\State Generate the same DRM at each local machine.
\For{$i = 1\dots m$}
\State $(\mathbf{V}_1^{(i)}, \cdots,\mathbf{V}_n^{(i)}, \T{H}^{(i)}) \leftarrow$ ComputeSketch($\T{X}_i$)
\EndFor
\For{$j = 1\dots n$}
\State $\mathbf{V}_j\leftarrow \sum_{i=1}^m \mathbf{V}_j^{(i)}$
\EndFor
\State $\T{H} \leftarrow \sum_{i=1}^m \T{H}^{(i)}$
\State \Return $(\mathbf{V}_1, \dots, \mathbf{V}_n, \T{H})$
\EndFunction
\end{algorithmic}
\end{algorithm}

\section{Scrambled Subsampled Randomized Fourier Transform} \label{appendix: ssrft}

In order to reduce the cost of storing the test matrices, in particular, $\mathbf{\Omega}_1, \dots, \mathbf{\Omega}_N$, we can use the Scrambled Subsampled Randomized Fourier Transform (SSRFT). To reduce the dimension of a matrix, $\mathbf{X} \in \mathbb{R}^{m \times n}$, along either the row or the column to size $k$, we define the SSRFT map $\mathbf{\Xi}$ as: 
\begin{equation}
\mathbf{\Xi} = \begin{cases}\mathbf{R}\mathbf{F}^\top \mathbf{\Pi}\mathbf{F}\mathbf{\Pi}^\top \in \mathbb{F}^{k \times m} & \text{(Row linear transform)}\\ 
(\widebar{\mathbf{R}}\widebar{\mathbf{F}}^\top\widebar{\mathbf{\Pi}}\widebar{\mathbf{F}}\widebar{\mathbf{\Pi}}^\top)^\top \in \mathbb{F}^{n \times k} & \text{(Column linear transform)} , \nonumber
\end{cases}
\end{equation}
where $\mathbf{\Pi}, \mathbf{\Pi}' \in \mathbb{R}^{m \times m}, \widebar{\mathbf{\Pi}}, \widebar{\mathbf{\Pi}}' \in \mathbb{R}^{n \times n}$ are signed permutation matrices. That is, the matrix has exactly one non-zero entry, 1 or -1 with equal probability, in each row and column. $\mathbf{F} \in \mathbb{F}^{m \times m}, \mathbf{F} \in \mathbb{F}^{n \times n}$ denote the discrete cosine transform ($\mathbb{F} = \mathbb{R}$) or the discrete fourier transform ($\mathbb{F} = \mathbb{C}$). The matrix $\mathbf{R}, \widebar{\mathbf{R}}$ is the restriction to $k$ coordinates chosen uniformly at random. 

In practice, we implement the SSRFT as in  \cref{alg:ssrft}. It takes only $\mathcal{O}(m)$ or $\mathcal{O}(n)$ bits to store $\mathbf{\Xi}$, compared to $\mathcal{O}(km)$ or $\mathcal{O}(kn)$ for Gaussian or uniform random map. The cost of applying $\mathbf{\Xi}$ to a vector is $\mathcal{O}(n\log n)$ or $\mathcal{O}(m \log m)$ arithmetic operations for fast Fourier transform and $\mathcal{O}(n\log  k)$ or $\mathcal{O}(m \log k)$ for fast cosine transform. Though in practice, SSRFT behaves similarly to the Gaussian random map, its analysis is less comprehensive \cite{boutsidis2013improved,tropp2011improved, ailon2009fast} than the Gaussian case. 

\begin{algorithm}[ht!] 
\begin{algorithmic}[1]
\caption{Scrambled Subsampled Randomized Fourier Transform (Row Linear Transform)}\label{alg:ssrft}
\Require{$\mathbf{X} \in \mathbb{R}^{m \times n}, \mathcal{F} = \mathbb{R}$, \textbf{randperm} creates a random permutaion vector, and \textbf{randsign} creates a random sign vector. \textbf{dct} denotes the discrete cosine transform.}
\Function{SSRFT}{$\mathbf{X}$}
\State \textbf{coords} $\leftarrow$ \textbf{randperm}(m,k) 
\State $\textbf{perm}_{j} \leftarrow \textbf{randperm}(m)$ for $j = 1,2$
\State $\textbf{sgn}_{j} \leftarrow \textbf{randsign}(m)$ for $j = 1,2$
\State $\mathbf{X} \leftarrow \textbf{dct}(\textbf{sgn}_1 \cdot \mathbf{X}[\textbf{perm}_1,:])$   \Comment{elementwise product}
\State $\mathbf{X} \leftarrow \textbf{dct}(\textbf{sgn}_2 \cdot \mathbf{X}[\textbf{perm}_2,:])$
\State \Return $\mathbf{X}[\textbf{coords},:]$
\EndFunction
\end{algorithmic}
\end{algorithm}

\section{TensorSketch} \label{appendix: TensorSketch}
Many authors have developed methods to perform dimension reduction efficiently.
In particular,
\cite{2017arXiv171209473D} proposed a method called \emph{TensorSketch} that aims
to solve least squares problems for which the design matrix has a Kronecker product structure.
\cite{malik2018low} use this technique to compute a one-pass Tucker decomposition.
Here we review the TensorSketch and how it is used in \cite{malik2018low}.

\paragraph{CountSketch} \cite{cormode2008finding} proposed the \textsf{CountSketch} method.
A comprehensive theoretical analysis in the context of low-rank approximation problems appears in \cite{clarkson2017low}.
To compute the sketch $\mathbf{X}\mathbf{\Omega} \in \mathbb{R}^{d \times k}$ for $\mathbf{X} \in \mathbb{R}^{m \times d}$,
\textsf{CountSketch} defines $\mathbf{\Omega} = \mathbf{D}\mathbf{\Phi}$, where
\begin{enumerate}
	\item $\mathbf{D} \in \mathbb{R}^{d \times d}$ is a diagonal matrix with each diagonal entry equal to $(-1,1)$ with probability $(1/2,1/2)$.
	\item $\mathbf{\Phi} \in \mathbb{R}^{d \times k}$ is the matrix form of a hash function.
\end{enumerate}

These two matrices have $2d$ non-zero entries in total
and thus require much less storage than the standard $kd$ entries.
Furthermore, these two matrices can operate on each column of $\mathbf{X}$
at a cost of only $\mathcal{O}(kd)$ arithmetic operations.

\paragraph{TensorSketch}
\cite{malik2018low} proposes to use the CountSketch inside the HOOI method for Tucker decomposition
They apply the sketch to solve least squares problems appearing
in \eqref{eq:factor-update} and \eqref{eq:core_update} in \cref{alg:hooi}.
They use $J_1, J_2$ to denote the reduced dimension.
Using a standard random map, it would require a $J_1$-by-$I_{(-n)}$ random matrix
to solve the problem in \cref{eq:factor-update}
and a $J_2$-by-$\prod_{n = 1}^N I_n$ random matrix to solve the problem in \cref{eq:core_update}.
However, these problems have Kronecker problem structure:
as shown in \cite{malik2018low}, these two stages can be expressed as
\begin{equation}\label{eq: tucker-stage-1}
\text{For } n = 1, \dots, N, \text{update } \mathbf{U}^{(n)}=\underset{\mathbf{U} \in \mathbb{R}^{I_{n} \times R_{n}}}{\arg \min }\left\|\left(\bigotimes_{i=N \atop i \neq n}^{1} \mathbf{U}^{(i)}\right) \mathbf{G}_{(n)}^{\top} \mathbf{U}^{\top}-\mathbf{Y}_{(n)}^{\top}\right\|_{F}^{2}.
\end{equation}

\begin{equation}\label{eq: tucker-stage-2}
\text{Update } \mathcal{G}=\underset{\T{Z} \in \mathbb{R}^{R_{1} \times \cdots \times R_{N}}}{\arg \min } \left\|\left(\bigotimes_{i=N}^{1} \mathbf{U}^{(i)}\right) \vc{\T{Z}}-\vc{\T{Y}}   \right\|_{2}^{2},
\end{equation}
where $\T{Y}$ is the original data.
Here $\forall i \in [n], \mathbf{U}_i$ is the factor matrix, and $\T{G}$ is the core tensor.
The target multilinear rank is $(R_1, \dots, R_N)$.

Following \cite{2017arXiv171209473D},
\cite{malik2018low} proposes to apply TensorSketch to the Kronecker product structure
of the input matrix in the sketch construction,
i.e. $\otimes_{\substack{i = 1\\ i \neq n}}^N \mathbf{U}_i$ in \cref{eq: tucker-stage-1} and $\otimes_{i =1}^N \mathbf{U}_i$ in \cref{eq: tucker-stage-2}.
The TensorSketch method combines the CountSketch of each factor matrix
via the Khatri-Rao product and Fast Fourier Transform.
Consider sketching $\otimes_{i =1}^N \mathbf{U}_i$ in \cref{eq: tucker-stage-2}.
TensorSketch is defined as
\begin{equation}\label{eq:tensorsketch}
\mathbf{\Omega}\mathbf{X}= \text{FFT}^{-1 }\bigg(\odot_{n =1}^N \Big(\text{FFT}\big(\text{CountSketch}^{(n)}(\mathbf{U}^{(n)}) \big)^\top \Big)^\top \bigg)
\end{equation}
By only storing $\text{CountSketch}^{(1)}, \dots, \text{CountSketch}^{(N)}$,
TensorSketch only requires storage $2\sum_{i=1}^N I_n$.
Therefore, the storage cost of the sketch is dominated by the sketch size,
$NR^{n-1}J_1 + J_2R^n \approx NKR^{2n-2}+KR^{2n}$, when $J_1 = KR^{n-1}, J_2 = KR^n$.

\end{appendices}

\end{document}